\documentclass[11 pt]{article}
\usepackage{a4wide, amsmath, amsthm, amssymb, amscd, mathptmx}

\usepackage{amsfonts}
\usepackage[only,llbracket,rrbracket]{stmaryrd}

\usepackage{color}
\usepackage{pstricks}
\usepackage{leftindex}

\usepackage[all]{xy}\CompileMatrices\SelectTips{cm}{12}

\theoremstyle{plain}
\newtheorem{Thm}{\sc Theorem}[section]
\newtheorem{Theorem}[Thm]{\sc Theorem}
\newtheorem{Corollary}[Thm]{\sc Corollary}

\newtheorem*{Corollary*}{\sc Corollary}

\newtheorem{Proposition}[Thm]{\sc Proposition}
\newtheorem*{Proposition*}{\sc Proposition}
\newtheorem{Lemma}[Thm]{\sc Lemma}

\theoremstyle{definition}
\newtheorem{Definition}[Thm]{Definition}

\theoremstyle{remark}
\newtheorem{Remark}[Thm]{Remark}
\newtheorem{Example}[Thm]{Example}
\newtheorem*{Example*}{Example}
\newtheorem*{Remark*}{Remark}

\renewcommand{\AA}{{\mathbb A}}

\newcommand{\GG}{{\mathbb G}}

\newcommand{\LL}{{\mathbb L}}
\newcommand{\MM}{{\mathbb M}}

\newcommand{\ZZ}{{\mathbb Z}}

\newcommand{\PP}{{\mathbb P}}

\newcommand{\RR}{{\mathbb R}}
\newcommand{\TT}{{\mathbb T}}
\newcommand{\VV}{{\mathbb V}}

\newcommand{\cA}{{\mathcal A}}

\newcommand{\cB}{{\mathcal B}}

\newcommand{\cC}{{\mathcal C}}
\newcommand{\cD}{{\mathcal D}}
\newcommand{\cE}{{\mathcal E}}
\newcommand{\cF}{{\mathcal F}}
\newcommand{\cG}{{\mathcal G}}

\newcommand{\cL}{{\mathcal L}}
\newcommand{\cM}{{\mathcal M}}
\newcommand{\cN}{{\mathcal N}}
\newcommand{\cO}{{\mathcal O}}

\newcommand{\cU}{{\mathcal U}}
\newcommand{\cV}{{\mathcal V}}

\newcommand{\Pic}{{\mathop{\rm Pic \, }}}

\newcommand{\Aut}{{\mathop{\operatorname{Aut }\, }}}

\newcommand{\an}{{\mathop{\rm an }}}
\newcommand{\Hom}{{\mathop{{\rm Hom}}}}

\newcommand{\cHom}{{\mathop{{\mathcal H}om}}}

\newcommand{\id}{\mathop{\rm Id}}

\newcommand{\Ext}{{\mathop{{\rm Ext \,}}}}

\newcommand{\cEnd}{{\mathop{{{\mathcal E}nd\,}}}}

\newcommand{\coker}{{\mathop{\rm coker \, }}}

\newcommand{\rk}{{\mathop{\rm rk \,}}}

\newcommand{\Supp}{{\mathop{{\rm Supp \,}}}}
\newcommand{\Spec}{{\mathop{{\rm Spec\, }}}}

\newcommand{\Exc}{{\mathop{{\rm Exc \,}}}}
\newcommand{\Chow}{\mathop{\rm Chow\, }}

\newcommand{\ev}{{\mathop{\rm ev \, }}}

\usepackage{graphicx}

\providecommand{\keywords}[1]{\textbf{\textit{Keywords: }} #1}

\begin{document}

\markboth {\rm }{}

\title{Projective  contact log varieties}
\author{Adrian Langer} \date{\today}

\maketitle


{\footnotesize
	{\noindent \sc Address:}\\
	Institute of Mathematics, University of Warsaw,
	ul.\ Banacha 2, 02-097 Warszawa, Poland\\
	e-mail: {\tt alan@mimuw.edu.pl}
}

\medskip

\begin{abstract}
	We study contact structures on  smooth complex projective varieties with a simple normal crossing divisor, generalizing some well-known results concerning the non-logarithmic case. In particular, we describe the structure of elementary log contractions of such log varieties and we construct the corresponding contact structures.
\end{abstract}

\keywords{complex projective variety; contact structure;  Lie algebroid; log contraction; simple normal crossing divisor}

\section*{Introduction}

Let $X$ be a smooth complex projective variety. A \emph{contact structure} on $X$ is defined as a corank $1$ subbundle $\cF$ of the tangent bundle $T_X$, where  the pairing $\cF\times \cF\to T_X/\cF=\cL$ induced by the Lie bracket is everywhere non-degenerate. In this situation, $X$ must have odd dimension  $2r+1$. Dually, one can consider an $\cL$-valued $1$-form $\theta \in H^0(X, \Omega _X^1\otimes \cL )$. This form defines a contact structure if and only if the expression $\theta \wedge (d\theta) ^{\wedge r}$, which is a well-defined section of $\Omega _X^{2r+1}\otimes \cL^r$, has no zeroes.

Contact structures first emerged in real geometry and later became a significant area of study in complex geometry. In the complex setting, one of the motivations came from the work \cite{LeBrun-Salamon} of C. LeBrun and  S. Salamon, who proved that twistor spaces of compact quaternion-K\"ahler manifolds with positive scalar curvature are contact Fano varieties (see \cite{LeBrun-Salamon}). They also conjectured that every such variety is homogeneous.  Although there has been some recent progress on this conjecture (see \cite{Buczynski-Wisniewski} and \cite{Occhetta-Romano-SolaConde-Wisniewski}), it still remains open.

One of the central classification results for contact complex varieties comes from the combined results of \cite{KPSW} and \cite{Demailly-Contact}. These works show that a projective contact variety is either the projectivization of the  tangent bundle of a smooth projective variety or it is Fano with $b_2=1$.

The main aim of this paper is to study generalization of this result to the logarithmic case. We assume that all varieties are defined over a fixed algebraically closed field $k$ of characteristic $0$ and  we consider contact structures on a projective snc pair, i.e., a pair consisting of a smooth projective variety $X$
and a simple normal crossing divisor $D$. In this case a contact structure is a corank $1$ subbundle $\cF$ of the logarithmic tangent bundle $T_X(-\log D)$ such that the pairing $\cF\times \cF\to T_X(-\log D)/\cF=\cL$ induced by the Lie bracket is everywhere non-degenerate. We usually forget about $\cF$ and we think about the contact structure as the projection $T_X(-\log D)\to \cL$.  Dually, one can  consider an $\cL$-valued logarithmic  $1$-form $\theta \in H^0(X, \Omega _X^1(\log D)\otimes \cL )$. This form is contact if and only if $\theta \wedge (d\theta) ^{\wedge r}$ has no zeroes.  $(X,D)$ is called \emph{contact} if it admits a contact structure (or equivalently a contact form). Such structure is a natural compactification of a contact structure on the quasi-projective variety $U=X\backslash \Supp D$. 

It is natural to expect that, similar to the usual case,  projective contact snc pairs can also be classified. 
It is straightforward to show that the projectivization of the logarithmic tangent bundle of a projective snc pair is contact  (see Proposition \ref{standard-contact-structure}). In this case, the contact form is  given by a logarithmic version of Liouville's $1$-form (cf. \cite[Appendix 4 D]{Arnold-Mathematical-methods} for the classical case). 

One of the key surprises in this theory is the existence of contact  projective snc pairs with $b_2\ge 2$, where the contact structure is not related to such standard contact structures. This phenomenon can be explained as follows. It is well-known in differential geometry and partial differential equations that higher jet bundles of a given vector bundle on a manifold carry a canonical contact form associated with the Cartan distribution (see, e.g., \cite{Saunders-Jet-book}).
This contact structure differs in that it is given by a distribution $\cF\subset T_X$, where $T_X/\cF$ is a vector bundle of possibly higher rank than $1$.  However, in the case of the first jet bundle of a line bundle it gives the usual contact structure (see, e.g., \cite[Example 3.44]{McDuff-Salamon} for the case of the trivial line bundle).  The corresponding contact form is not homogeneous with respect to fiber coordinates, but it extends to a contact form on the natural log variety underlying certain projective bundle compactifying this jet bundle.  This explains how log contact varieties of dimension $(2r+1)$ can arise as projective bundles over $r$-dimensional varieties. More generally,  for any projective snc pair $(Y,B)$ of dimension $d\le r$ and any collection of $(r+1-d)$ line bundles on $Y$, we provide an explicit construction of a contact snc pair of dimension $(2r+1)$, describing it as a projective bundle over $Y$ (see Proposition \ref{non-standard-contact-structure}).
In particular, the classification of projective contact snc pairs is significantly more complex than that of usual contact varieties.   

What is surprising from this point of view is that our construction describes essentially all new contact structures. More precisely, we prove the following theorem, which is a special case of Theorem \ref{structure-of-log-contractions-special-Lie-alg}. 
The $c_1^{\log}(\cM)$ in this theorem denotes the image of the first Chern class $c_1(\cM)\in H^1(Y, \Omega_Y^1)$ in 
$H^1(Y, \Omega_Y^1(\log B))=\Ext^1 (T_{Y}(-\log B), \cO_Y )$.

\begin{Theorem}\label{main1}
Let $(X,D)$ be a contact projective snc pair of dimension $2r+1$. Then $-(K_X+D)$ is not nef and $X$ admits an elementary log contraction $\varphi : X\to Y$ for which $-(K_X+D)$  is $\varphi$-ample. If the restriction of $\varphi$ to $X\backslash \Supp D$ is not an isomorphism then either $(X,D)$ is log Fano with $b_2(X)=1$ or 
$Y$ is smooth and there exists a simple normal crossing divisor $B$ on $Y$
such that  one of the following occurs:
	\begin{enumerate}
		\item  $\dim Y=r+1$  and $(X\to Y)\simeq (\PP( T_Y (-\log B)) \to Y$) and $D=\varphi^{-1}(B)$,
		\item  $\dim Y=r+1-s$ for some $1\le s\le r$ and there exists $s$ line bundles $\cM_1,...,\cM_s$ on ${Y}$ 
		such that $(X\to Y)\simeq (\PP (\cE) \to Y)$, where we have an extension
$$0\to \bigoplus _{i=1}^s\cM_i{\longrightarrow} \cE \longrightarrow\cC \to 0,$$
in which $\cC$ is the extension of $T_{Y}(-\log B)$ by $\cO _Y^{\oplus s}$ given by 
		$-\bigoplus _{i=1}^sc_1^{\log} (\cM_i)$. In this case $D$ is the sum of $\varphi^{-1}(B)$
		and $s$ divisors determined by the cokernels of $\cM_i\to \cE$.
	\end{enumerate}
\end{Theorem}

Our main contribution is description  of log contractions and construction of contact structures in the corresponding cases. The first part of the above theorem is a simple corollary of the result of  F. Touzet \cite{Touzet2016} (see Corollary \ref{Touzet}) generalizing earlier result of J.-P. Demailly \cite{Demailly-Contact} to the log case. In fact, this result shows that $K_X+D$ is not pseudoeffective. So by \cite[Theorem 2]{BDPP2013} there exist on $X$  strongly moving curves $C$ with $(K_X+D)\cdot C<0$. In this case the log Bend and Break conjecture \cite[Conjecture 1.11]{Keel-McKernan} of S. Keel and J. McKernan predicts that the quasi-projective variety $U=X\backslash \Supp D$ is $\AA^1$-uniruled (in the sense of, e.g., \cite[Definition 1.2]{Chen-Zhu}). In case $(X,D)$ admits an elementary log contraction which is not an isomorphism on the open part $U=X\backslash \Supp D$,  this fact follows easily from the above result.

\medskip

We also prove a far reaching generalizations of Theorem \ref{main1} that describe a precise structure of all good log contractions of contact snc pairs. Unlike in the case of contact varieties,  there also appear birational log contractions. In this case the whole exceptional locus is contained in the support of $D$ and the target is often singular. But divisorial contractions can also lead to snc pairs, e.g., we can obtain a contact snc pair by blowing up the fiber of $\PP(T_Y(-\log B))\to Y$ over a point that does not belong to the support of $B$ (see Example \ref{blow-up-works}).

\medskip

 To study log contractions of contact snc pairs, it is essential to explore more general snc pairs with what is called a \emph{semi-contact structure} (see Definition \ref{semi-contact-definition}). A semi-contact structure is defined as a contact structure on some special Lie algebroid generalizing a Picard Lie algebroid in the sense of A. Beilinson and J. Bernstein \cite[Definition 2.1.3]{Beilinson-Bernstein-Jantzen}. This structure enables the description of the restriction of the contact structure on the pair $(X,D)$ to the strata of the stratification of $X$ induced by the components of $D$. A clear advantage of this approach is that it allows for inductive reasoning on the dimension, whereas conventional  contact structures only exist in odd dimensions. 

\medskip

The classification of contact varieties with $b_2\ge 2$ in \cite{KPSW} used Mori's theory, specifically focusing on contractions and rational curves on contact varieties.  Therefore,  it is natural to approach the classification of contact snc pairs using log contractions and the non-proper analogue of rational curves suggested by S. Keel and J. McKernan in \cite{Keel-McKernan} (see also \cite{Chen-Zhu}).  However, many standard results do not directly translate to the log case, and to prove  Theorem \ref{main1}, we must formulate appropriate analogues that are applicable in the log setting. Specifically, we require a relative version of  Keel--McKernan's log bend-and-break \cite[Corollary 5.4]{Keel-McKernan}, a log version of the Ionescu--Wi{\'{s}}niewski inequality, and a log version of the rank computation for  $\Hom(\PP^1, X)$ (see \cite[Chapter II, Proposition 3.10]{Kollar-Rational} for the standard case). Additionally, several results concerning Lie algebroids and generalized Atiyah extensions are extensively used in this paper. 
 With these tools, we attempt to mirror the proof from \cite{KPSW}, though the process becomes quite intricate, leading to various additional challenges that are typically resolved through a careful study of specific loci and the use of dimensional induction.

\medskip

A contact structure can more generally be defined on  a log smooth log variety in the sense of Fontaine--Illusie--Kato (see \cite{Kato-Logarithmic-structures} and \cite{Ogus-log}). In this context, the contact structure is defined by the same data as in teh usual case,  with the only difference being the replacement of a smooth variety with 
a log smooth log variety. These structures also appear as natural compactifications of contact structures on quasi-projective varieties.  In this paper, we  focus on contact snc pairs, but analogous results can also be proven  in the logarithmic case. Unfortunately, this requires  significant additional work, which would obscure the main ideas, so we have chosen to first explain the classical case and address the logarithmic case in a separate paper.

As mentioned earlier,  log contractions of  contact snc pairs (or log smooth log varieties) that are not of fiber type often result in singular varieties. The study of singular contact log varieties and contact log Fano log varieties with $b_2=1$ is also deferred to another paper. We also leave out potential applications of the above theory in complex or real geometry.

\medskip

In this paper we focus on algebraic varieties. This approach  does not restrict the generality of the obtained results, as complex projective manifolds are well-known to be algebraic. It should also be noted that much of this paper could be written in the language of complex manifolds, but we prefer to use the language of algebraic geometry. While this choice might seem unnecessary to some readers  (and it can make the paper more difficult to read), there are mathematical reasons for doing so. One reason is that in some cases, our study leads to singular (projective) varieties, and singularities in complex geometry are not as rigid as those of algebraic varieties. However, the most important reason is that we aim to show that contact  projective (log) varieties have a well-defined local algebraic structure. This might be surprising, as it is well-known that Darboux's theorem, which describes the local structure of a symplectic or contact manifold, does not apply in algebraic local coordinates.

\subsection*{Relation to other work}

Contact structures are well-known to be related to symplectic structures. Rather recently, differential geometers started to study a generalization of symplectic structures under the name of b-symplectic manifolds. These structures first appeared in \cite{Guillemin-Miranda-Pires} and they are described in the language of b-geometry refering to Melrose's calculus for differential operators on real manifolds with (smooth) boundary (see \cite{Melrose-book}).  Complex analogues of symplectic structures in case of complex manifolds with boundary have also been heavily studied under the name of log symplectic manifolds (see, e.g.,  \cite{Cavalcanti-Gualtieri2018}, \cite{Goto2002}, \cite{Goto2016} and \cite{Pym2017}).  
 
Even more recently, b-contact manifolds were defined in \cite{Miranda-Oms} (see also \cite{Bradell-et-al}). The definition is the same as in the case of usual manifolds but it replaces manifolds with b-manifolds (tangent bundle with b-tangent bundle etc.). Our contact snc pairs are an obvious analogue of such structures in complex geometry but it seems they have not been studied till now. Let us also note that our  constructions work in the real case giving a new construction of (compact) b-contact manifolds that compactify $1$-jet bundles (see Appendix B). It seems that this construction has not been noticed before.

\medskip

The structure of the paper is as follows. In Section 1 we generalize Wahl's theorem on ample line bundles of a tangent bundle to the logarithmic case and we recall some results on the logarithmic cone and contraction theorems. In Section 2 we adapt several theorems on rational curves on smooth projective varieties to the log case. In Section 3, we extend contact and symplectic structures to the Lie algebroid case and prove a few auxiliary results. Then in Section 4 we introduce special Lie algebroids that generalize Picard algebroids and we examine their contact structures. Section 5 is dedicated to proving the main result. In Section 6 we present new constructions of contact snc pairs, while in Section 7 we investigate contact projective snc pairs with two elementary log contractions, demonstrating the abundance of such varieties. The appendices provide a classification of semi-contact snc pairs in small dimensions (Appendix A) and demonstrate how analogues of the Section 6 constructions yield new b-contact b-manifolds (Appendix B).

\medskip

\subsection*{Notation}

In the whole paper, we assume that our varieties are defined over an algebraically closed field $k$ of characteristic $0$. A variety is always irreducible. Since all structures in this paper are $k$-linear, by abuse of notation we write $\Omega_X^1$ for the sheaf of  K\"ahler differentials $\Omega_{X/k}$ and $T_X$ for the sheaf  $T_{X/k}=Der _k (\cO_X, \cO_X)$ of $k$-linear derivations of $\cO_X$. For a coherent $\cO_X$-module $\cE$ we write $\cE^*$ for $\cHom _{\cO_X}(\cE, \cO_X)$. All tensor operations on $\cO_X$-modules are performed over $\cO_X$ unless otherwise stated explicitly.

A \emph{vector bundle} on $X$ is a locally free coherent $\cO_X$-module of finite rank.

A \emph{contraction} is a proper surjective map $\varphi: X\to Y$ of normal varieties such that $\varphi_*\cO_X=\cO_Y$ (or equivalently $\varphi$ has connected fibers).

\section{Preliminaries}

In this section we generalize Wahl's characterization of projective spaces \cite{Wahl-Ample-in-tangent} to log varieties, recall the logarithmic cone (and contraction) theorem and prove some auxiliary results on log extremal curves.

\subsection{Logarithmic version of Wahl's theorem}

\begin{Theorem}\label{log-Wahl}
Let $(X, D)$ be a projective snc pair and let $\cL$ be an ample line bundle on $X$. Let $n=\dim X\ge 1$. 
If $H^0(X, T_X(-\log D)\otimes \cL^{-1})\ne 0$ then either $(X, \cL)\simeq (\PP^1, \cO_{\PP^1}(2))$ and $D=0$ or 
$(X, \cL)\simeq (\PP^n, \cO_{\PP^n}(1))$ and $D$ is an arrangement of $s\le n$ hyperplanes in a general position
(i.e., it forms an snc divisor). 
\end{Theorem}

\begin{proof}
The proof is by induction on $n$ with $n=1$ being clear.
Since $T_X(-\log D)\subset T_X$, the first part follows immediately from Wahl's theorem \cite[Theorem 1]{Wahl-Ample-in-tangent}. So we can assume that $(X, \cL)=(\PP^n, \cO_{\PP^n}(1))$ and $D\ne 0$.
For any irreducible component $D_i$ of $D$ we have a short exact sequence
$$0\to T_{D_i}(-\log (D-D_i)|_{D_i})\to T_{\PP^n}(-\log D)|_{D_i}\to \cO _{D_i} (-(D-D_i))\to 0. $$
By assumption we have a nonzero map $\cL\to  T_{\PP^n}(-\log D)$. Its restriction to $D_i$ gives a nonzero
 map $\cL|_{D_i}\to  T_{\PP^n}(-\log D)|_{D_i}$. Since there are no nonzero maps from $\cL|_{D_i}$ to 
 $\cO _{D_i} (-(D-D_i))$, we have an induced non-zero map $\cL|_{D_i}\to T_{D_i}(-\log (D-D_i)|_{D_i})$. Using the induction assumption we see that $D_i$ is a hyperplane and $(D-D_i)|_{D_i}$ consists of at most $(n-1)$ hyperplanes in $D_i\simeq \PP^{n-1}$. This easily implies the required assertion. 
\end{proof}

\begin{Remark}
	Let $D$ be an arrangement of $1\le s\le n$ hyperplanes in a general position in $\PP^n$. Then
 by \cite[Proposition 2.10]{Dolgachev-Kapranov} we have
$$T_{\PP^n}(-\log D)\simeq \cO_{\PP^{n}}^{\oplus (s-1)}\oplus \cO_{\PP^{n}}(1)^{\oplus (n-s+1)}.$$
This shows that the above theorem is optimal. Note also that $D$ is projectively equivalent
to the standard hyperplane arrangement given by $x_1...x_s=0$. In particular, intersection of any subset of hyperplanes in $D$ is connected.
\end{Remark}

The above theorem and the main result of M. Andreatta and J. Wi{\'{s}}niewski \cite[Theorem]{Andreatta-Wisniewski2001} immediately imply the following corollary:

\begin{Corollary}
	Let $(X, D)$ be a projective snc pair with $n=\dim X\ge 1$ and let $\cE$ be an ample vector bundle contained in $T_X(-\log D)$. Then we have the following possibilities for $(X,D)$ and $\cE$:
	\begin{enumerate}
		\item  $(X,D)\simeq (\PP^n ,0)$ and $\cE \simeq \cO_{\PP^n}(1)^{\oplus m}$ for some $m\le n$,
		\item  $(X,D)\simeq (\PP^n ,0)$ and $\cE\simeq T_{\PP^n}$,
		\item $X\simeq \PP^n$,  $D$ is an arrangement of $1\le s\le n$ hyperplanes and $\cE\simeq \cO_{\PP^{n}}(1)^{\oplus m}$ for some $m\le {(n-s+1)}$.
	\end{enumerate}
\end{Corollary}

\subsection{Logarithmic cone theorem}

Let $X$ be a smooth projective variety.
Let $Z_1(X)_{\RR}$ be the real vector space generated by integral curves
on $X$, and let $N_1(X)$ be the space  $Z_1(X)_{\RR}$ modulo numerical equivalence. We denote by
$ N^1(X)$ the space of real Cartier divisors modulo numerical equivalence. We also set
$\rho(X):=\dim_{\RR}N^1(X).$

Let $NE(X)\subset N_1(X)$ be the convex cone generated by classes
of effective $1$-cycles and let $\overline{NE}(X)$ be its closure. For
$L\in N^1(X)$, we use the notation
$ \overline{NE}(X)_{L\geq0}
:=
\bigl\{
z\in\overline{NE}(X)\mid L\cdot z\geq0
\bigr\}$
and
$N_1(X)_{L<0}
:=
\bigl\{
z\in N_1(X)\mid L\cdot z<0
\bigr\}.$

A subcone $F\subset\overline{NE}(X)$ is called an \emph{extremal face}
if, whenever $z_1,z_2\in\overline{NE}(X)$ and $z_1+z_2\in F$, one has
$z_1,z_2\in F$. It is called \emph{$(K_X+D)$-negative} if
$(K_X+D)\cdot z<0$
for every nonzero $z\in F$. A one-dimensional extremal face is called an
\emph{extremal ray}.  
A $(K_X+D)$-negative extremal face $F$ is called a \emph{$(K_X+D)$-extremal face}.

A nef Cartier divisor $H_F$ is called a \emph{supporting divisor} of an
extremal face $F$ if
\[
F=
\bigl\{
z\in\overline{NE}(X)\mid H_F\cdot z=0
\bigr\}.
\]

We need the following version of the cone theorem that follows from the cone theorem for dlt pairs (see \cite[Theorems 3.25(2) and 3.35]{Kollar-Mori-book} and \cite[Chapter III, Corollary 1.4]{Kollar-Rational}).

\begin{Theorem}\label{cone-theorem}
	Let $(X,D)$ be a projective snc pair. Then: 	
	\begin{enumerate}
		\item There are countably many rational curves $C_j\subset X$ such that $(K_X+D)\cdot C_j<0$ and
		$$\overline{NE} (X)=\overline{NE} (X) _{K_X+D\ge 0}+\sum _j \RR^{+}[C_j].$$
		The rays $\RR^{+}[C_j]\subset N_1(X)$ are locally discrete in $N_1(X)_{K_X+D<0}$ and every extremal ray $R$  contained in $N_1(X)_{K_X+D<0}$ is of the form $R=\RR^{+}[C_j]$ for some $j$.
		\item For any $(K_X+D)$-extremal face $F$ 
		there is a unique contraction $\varphi_F: X\to Y$ to projective $Y$ such that an irreducible curve $C\subset X$ is mapped by $\varphi _F$ to a point if and only if $[C]\in F$. 
		\item Every $(K_X+D)$-negative extremal face $F$ admits a supporting
		divisor $H_F$. Moreover, $$mH_F-(K_X+D)$$
		is ample for all sufficiently large $m$.  
	\end{enumerate}
\end{Theorem}

A contraction $\varphi : X\to Y$ is called a \emph{good log contraction} of $(X,D)$ if $-(K_X+D)$ is $\varphi$-ample. Any good log contraction is of the form $\varphi_F$ for some extremal face $F$ (see \cite[Lemma 3.2.5]{Kawamata-Matsuda-Matsuki}). 
A $1$-dimensional $(K_X+D)$-extremal face is called a \emph{ $(K_X+D)$-extremal ray}. If $R$ is a $(K_X+D)$-extremal ray then the good log contraction $\varphi_R$ is called \emph{elementary}.

Let $(X,D)$ be a projective snc pair such that $(K_X+D)$ is not nef. Then either $\rho (X)=1$ and $-(K_X+D)$ is ample (the log Fano case) or $\rho (X)\ge 2$ and there exists a good elementary log contraction $\varphi : X\to Y$ with $\dim Y>0$.

\medskip

If $D=\sum _{i\in I}D_i$ for some smooth irreducible $D_i$ then we write $D_J:=\bigcap _{j\in J} D_j$ for $J\subset I$. In particular, we set $D_{\emptyset}=X$. For a subset  $J\subset I$ we denote by $D^{J}$ the restriction of the divisor 
$\sum _{i\in I\backslash J}D_i$ to $D_J$. 
Let us note the following complement to the cone theorem:

\begin{Proposition}\label{cone-theorem-complement}
	Let $D=\sum _{i\in I}D_i$ and let us fix some $J\subset I$ (possibly $\emptyset$).
	If $F$ is a $(K_X+D)$-extremal face and  $C\subset D_J $ is a curve such that $[C]\in F$
	then we can find a rational curve $C'\subset D_J$ such that $[C']\in F$ and
	$$0< -(K_X+D)\cdot C'\le \dim D_J+1.$$	
	Moreover, if  $C\not \subset D_i$ for $i\in I\backslash J$ then we can assume that $C\cap C'\ne \emptyset$ and 
	$C'\not \subset D_i$ for $i\in I\backslash J$.
\end{Proposition}

\begin{proof}
Enlarging $J$ if necessary we can assume that $J$ is the set of all $j\in I$ for which $C\subset D_j$. Then  we have
	$$0>(K_X+D)\cdot C=(K_{D_J}+D^{J})\cdot C\ge K_{D_J}\cdot C.$$ 
	Let $H_F$ be a supporting divisor for $F$. 	 By \cite[Chapter II, Theorem 5.8]{Kollar-Rational} for any closed point $x\in C$ we can find a rational curve $C'$ passing through $x$ and such that $$H_F\cdot C'\le 2 \dim D_J \frac{H_F\cdot C}{-K_{D_J}\cdot C}=0.$$
	Then $H_F\cdot C'=0$ and hence $[C']\in F$. This shows that we can assume that $C$ is rational.
	By \cite[Chapter II, Theorem 5.7]{Kollar-Rational} for any $x\in C\subset D_{J}$ there is an effective $1$-cycle $\sum a_iC_{i}$ in $D_J$ whose irreducible components are rational curves such that $[C]$ is effectively algebraically equivalent to $\sum a_{i}C_{i}$,   $-K_{D_{J}}\cdot C_{i}\le \dim D_J+1$ for all $i$, and $x\in \Supp \sum a_{i}C_{i}$. Since $0=H_F\cdot C=\sum a_i H_{F}\cdot C_{i}\ge a_i H_{F}\cdot C_{i}\ge 0$, we see that
	$H_{F}\cdot C_{i}=0$ and hence $[C_i]\in F$. Since $x\in \Supp \sum a_{i}C_{i}$ we can also take $i$ so that $x\in C_i$. Choosing $x$ so that $x\not \in D _{j}$ for $j \in I\backslash J$, we obtain
	$$0<-(K_X+D)\cdot C_i=-(K_{D_J}+D^{J})\cdot C_i\le -K_{D_J}\cdot C_i \le \dim D_J+1.$$ 
\end{proof}

\medskip

For a  good log contraction $\varphi $ we define its \emph{length} 
$l_D(\varphi)$ as the minimum of $-(K_X+D)\cdot C$ taken among all rational curves $C$ contracted by $\varphi$. If the minimum is attained for $C$  then we call $C$ an \emph{extremal rational curve}.
By  Proposition \ref{cone-theorem-complement} we have $1\le l_D(\varphi)\le \dim X+1.$

Let $X$ be a smooth projective variety and let  $ \Hom (\PP^1, X)$ be the scheme parametrizing morphisms from $\PP^1$ to $X$. This scheme comes with the evaluation morphism $F: \PP^1\times \Hom (\PP^1, X) \to X$. For any family $V\subset \Hom (\PP^1, X)$ we write $F_V$ for the restriction of $F$ to $\PP^1\times V$. We denote by ${\rm locus} (V)$ the image of $F_V$.
The family $V$ is called \emph{unsplit} if curves, obtained as images of morphisms from this family, cannot be deformed to a $1$-cycle consisting of more than one curve (counting multplicities). We refer to \cite[Chapter IV, Definition 2.1]{Kollar-Rational} for a precise definition.

\begin{Lemma}\label{bound-on-length}
If  $ l_D(\varphi)=\dim X+1$ then  $X\simeq \PP^n$  and $D=0$.
In particular, if $D\ne 0$ then $l_D(\varphi)\le \dim X$. 
\end{Lemma}

\begin{proof}
	Assume that  $l_D(\varphi)=\dim X+1$. Since any good log contraction factors through an elementary good log contraction, we can assume that $\varphi$ is elementary. 
	
Then for any extremal rational curve $C$ contracted by $\varphi$ we have $-(K_X+D)\cdot C=\dim X+1$. So $C\not \subset \Supp D$ by Proposition \ref{cone-theorem-complement}. In particular, $-K_X\cdot C\ge -(K_X+D)\cdot C=\dim X+1$. Let $f: \PP^1\to C$ be the normalization and let $V$ be an irreducible component of $\Hom (\PP^1, X)$ containing $[f]$. 
The family $V$ is unsplit. Otherwise, a reducible degeneration would have components whose classes lie in the same extremal ray, and one component would have $-(K_X+D)$-degree strictly smaller than $l_D(\varphi)$.

By Ionescu--Wi\'{s}niewski's inequality (see, e.g., \cite[Chapter IV, Corollary 2.6]{Kollar-Rational}), $-K_X\cdot C=\dim X+1$ and  $\overline{{\rm locus} (V)}=X$. If $D=0$ then $\varphi$ contracts $X$ to a point  and hence $X$ is Fano. Then $X\simeq \PP^n$ by  \cite{Cho-Miyaoka-Shepherd-Barron} and \cite[Theorem 1.1]{Kebekus-CMSB}. If $D\ne 0$ then
$D\cdot C=0$ and  there exists $[g]\in V$ such that the intersection of $C'=g(\PP^1)$ with $D$ is non-empty. 
	Since $D\cdot C'=D\cdot C=0$,  $C'$ is contained in $ \Supp D$. But $\RR^{+}[C']=\RR^{+}[C]$, so we get a contradiction with Proposition \ref{cone-theorem-complement}.	
\end{proof}

The proof of the following lemma follows standard arguments used, e.g., in the proof of \cite[Theorem 1.6]{Kawamata1997}.

\begin{Lemma}\label{reduction-to-log-contractions-on-strata}
	Let $(X, D)$ be a projective snc pair and let $\varphi: X\to Y$ be a good log contraction of $(X, D)$. Let $Z$ be an irreducible component of $D$ and let $D^Z$ be the restriction of $(D-Z)$ to $Z$. Then $T=\varphi(Z)$ is normal and $\varphi|_Z: Z\to T$ is a good log contraction of $(Z, D^Z)$.
\end{Lemma}

\begin{proof}
	Since $-(K_X+D)$ is $\varphi$-ample, for small $\epsilon>0$ the divisor
	$-(K_X+D)+\epsilon (D-Z)$ is also $\varphi$-ample. Since
	$$-Z= K_X+ (-(K_X+D)+\epsilon (D-Z))+(1-\epsilon) (D-Z),$$
	the relative version of the Kawamata-Viehweg vanishing theorem implies that 
	$R^i\varphi _*(\cO_X (-Z))=0$ 	for all $i>0$.  
	In particular, the map 	$\cO_Y=\varphi _*\cO_X \to \varphi_*\cO_{Z}$ is surjective.
	So if $T=\varphi(Z)$ then $\cO_T\to  \varphi_*\cO_{Z}$ is surjective. Since it is also injective, it is an isomorphism and hence $T$ is normal. To see that  $\varphi|_Z$ is a good log contraction it is sufficient to note that
	$K_Z+D^Z= (K_X+D)|_Z$. 
\end{proof}

By induction the above lemma implies the following corollary.

\begin{Corollary} \label{normality-of-image}
	Let $\varphi: X\to Y$ be a good log contraction of a projective snc pair  $(X, D)$ with $D=\sum _{i\in I}D_i$. Let us fix $J\subset I$ and an irreducible component $Z$ of $D_J$. Then $T=\varphi(Z)$ is normal and $\varphi|_Z: Z\to T$ is a good log contraction of $(Z, D^J|_Z)$.
\end{Corollary}

\section{Rational curves on snc pairs}\label{section:rational-curves}

 In this section we generalize several results on rational curves on smooth projective varieties to projective snc pairs.
 
\subsection{Log Bend-and-Break}

Let $X$ be a smooth projective variety. For any closed point $x\in X$ we consider the subscheme $$\Hom (\PP^1, X; 0\to x)$$ of  $\Hom (\PP^1, X)$ parametrizing $g: \PP^1\to X$  mapping $0\in \PP^1 $ to $x\in X$. For a closed point $x\in {\rm locus} (V)$, we write $V^x$ for the family $V\cap  \Hom (\PP^1, X; 0\to x)$ and ${\rm locus} (V, x)$ for the locus of $V^x$.

\medskip

Now let $(X,D)$ be a  projective snc pair and let us set $U:= X\backslash \Supp D$. Let 
$\{D_i\}_{i\in I}$ be the set of irreducible components of $D$ and let $G=\{G_i\}_{i\in I}$ be a collection of effective Cartier divisors on $\PP^1$. Then we can consider 
the closed subscheme  $\Hom (\PP^1, X; G\subset D)$ of $\Hom (\PP^1, X)$  parametrizing morphisms $g: \PP^1\to X$ such that for every irreducible component $D_i$ of $D$ we have $G_i\subset g^*D_i$ (see \cite[5.1]{Keel-McKernan}).
For any $x\in U$,  $\Hom (\PP^1, X; G\subset D, 0\to x)$ denotes the scheme-theoretic intersection of the schemes $\Hom (\PP^1, X; G\subset D)$ and $ \Hom (\PP^1, X; 0\to x)$ inside  $ \Hom (\PP^1, X)$.

Let $f: \PP^1\to X$ be a fixed morphism such that $f(\PP^1)\not\subset \Supp D$ and  $G_i\subset f^*D_i$ for all $i\in I$.  Let $W$ be an irreducible component of $\Hom (\PP^1, X; G\subset D)$ containing $[f]$. Then
\cite[Proposition 5.3]{Keel-McKernan} implies that
$$\dim _{[f]}W  \ge \chi (\PP^1, f^*T_X)-\deg G= -\deg f^*K_X-\deg G+\dim X.$$
Moreover, if $G_i=f^*D_i$ for all $i\in I$ then
$$ T_{[f]}\Hom (\PP^1, X; G\subset D)=H^0(\PP^1, f^*T_X(-\log D)).$$

\medskip

The proof of the following proposition follows that of \cite[Corollary 5.4]{Keel-McKernan}
and its improvement \cite[Proposition 6.4]{Chen-Zhu}.

\begin{Proposition}\label{Keel-McKernan}
	Assume that $K_X+D$ is not nef and let $C$ be an extremal rational curve for some $(K_X+D)$-extremal face $F$.  Let $J\subset I$ be the set of all $j\in I$ for which  $C\subset D_j$ and let  $U_J=D_{J}\backslash \Supp D^{J}$.	 Then for every closed point $x\in C\cap U_J$ there exists a morphism $g: \PP^1\to D_J$ such that $x\in g(\PP^1)$, $[g(\PP^1)]\in F$ and $g^*(D^{J})$ is supported in at most one point.  
\end{Proposition}

\begin{proof}
	Let $f: \PP^1\to C$ be the normalization of $C$ and let $V$ be the connected component of 
	$\Hom (\PP^1, D_J)$ containing $[f]$. 
	The \emph{size} of a $0$-dimensional scheme $G$ of $\PP^1$ is  the maximum length of a subscheme contained in $G$, which is supported at a single point. Let $g: \PP^1\to D_J$ be a morphism with $[g]\in V^x$. Since $x\in C':=g(\PP^1)$, we have 	$C'\not \subset \Supp D^{J}$ and $g^*(D^{J})$ is a well defined Cartier divisor. Let us choose $g$  for which $g^*(D^{J})$ has the greatest size that we denote by $s$. We claim that $g$ satisfies the required assertions.

	Let $G\subset g^*(D^{J})$ be a subscheme of length $s$ supported at a single point.
	Since it is clear that $[C']\in F$, we need only to check the last condition. Assume that $g^*D^{J}$ is supported in $m\ge 2$ points.
	Let $W$ be an irreducible component of $\Hom  (\PP^1, D_J; G\subset D^{J}, 0\to x)$ containing $[g]$,
	and let $\tilde W$ be an irreducible component of $\Hom  (\PP^1, D_J)$ containing $W$. 
	Since $D^{J}\cdot C'\ge (m-1)+s$ and $W$ is an irreducible component of the fibre of $\tilde W\to D_J$, $h\to h(0)$, over $x$,  we have
	\begin{align*}
	\dim W&\ge \dim \tilde W-\dim D_J\ge{\chi (\PP^1, g^*T_{D_J})-\deg G-\dim D_J=-K_{D_J}\cdot C'-s}\\
	&\ge  {
	-(K_{D_J}+D^{J})\cdot C'+m-1= -(K_X+D)\cdot C+m-1\ge 2}.	
	\end{align*}
	Note that $C'$ is an extremal rational curve, so for large $m\gg 0$ it has minimal degree with respect to ample divisors of the form $-(K_X+D)+mH_F$. Hence we get a contradiction in the same way as in the proof  of \cite[Corollary 5.4]{Keel-McKernan} (or  \cite[Proposition 6.4]{Chen-Zhu}).
\end{proof}

\begin{Remark}
	Note that $D_J$ need not be irreducible although every connected component of $D_J$ is smooth and irreducible. 
\end{Remark}

\subsection{Loci of rational curves on snc pairs}

As pointed out in the proof of \cite[Proposition 2.9]{KPSW}, the following proposition due to P. Ionescu and J. Wi{\'{s}}niewski follows from \cite[Chapter IV, Proposition 2.5]{Kollar-Rational} (note that unlike claimed in \cite{Kollar-Rational}, the proof does not give any information if we do not assume that $V$ is generically unsplit).

\begin{Proposition}\label{easy-log-Ionescu-Wisniewski}
	Let $X$ be a smooth projective variety and let $C$ be a rational curve in $X$. Let $f: \PP^1 \to C$ be the normalization of $C$ and  let $V$ be a closed irreducible subset of $\Hom (\PP^1, X)$, which contains $[f]$ and which is invariant under the action of $\Aut (\PP^1)$.
	If the family $V$ is unsplit then for any closed point  $x\in {\rm locus} (V)$ we have
	$$\dim V\le \dim  {\rm locus} (V)+ \dim  {\rm locus} (V,x)+1.$$
\end{Proposition}

\medskip

Let $(X,D)$ be a projective snc pair with $D\ne 0$ and let $U= X\backslash \Supp D$ as before.
For any family $V\subset \Hom (\PP^1, X)$ we set ${\rm locus} _U(V):={\rm locus} (V)\cap U$ and $ {\rm locus} _U(V,x):= {\rm locus} (V,x)\cap U$. 

If $f: \PP^1\to X$ is a fixed morphism such that $f(\PP^1)\not\subset \Supp D$ then
we can consider $f^*D$ as a collection of  effective Cartier divisors on $\PP^1$. If $g$ with $g(\PP^1)\not\subset \Supp D$ is in the same connected component of $\Hom (\PP^1, X; f^*D\subset D)$ as $f$ then $f^*D_i= g^*D_i$ for all $i\in I$. This follows from the fact that  $\deg f^*D=\deg g^*D$ and $f^*D_i\le g^*D_i$ for $i\in I$. 

\medskip 

The following proposition is a logarithmic analogue of the above proposition.

\begin{Proposition}\label{log-Ionescu-Wisniewski}
	Let $C$ be a rational curve in $X$ and let $f: \PP^1 \to C$ be the normalization of $C$.
	Assume that $f^{-1}(U)\simeq \AA^1$
	and let $W$  be a closed irreducible subset of $\Hom (\PP^1, X; f^*D\subset D)$ containing $[f]$ and invariant under the action of the subgroup of $\Aut (\PP^1)$ fixing the support of $f^*D$.  If the family $W$ is unsplit  then  for  any closed point  $x\in {\rm locus} _U(W)$ we have
	$$\dim W \le \dim  {\rm locus} _U(W)+ \dim  {\rm locus} _U(W,x).$$
\end{Proposition}

\begin{proof}
	Let $\tilde{W}= \Aut (\PP^1)\cdot W\subset \Hom (\PP^1, X)$. Since $C$ is not contained in the support of $D$, $W$ is not $\Aut (\PP^1)$-invariant.  Hence
	$\dim \tilde W= \dim W+1$. Let us consider the map $F_{\tilde W}: \PP^1\times \tilde{W}\to X$. For any $x\in {\rm locus} _U(W)$ there exists $p\in \PP^1$ and $[g]\in \tilde W$ such that $g(p)=x$. Taking an automorphism $h$ of $\PP^1$ such that $h(0)=p$, we see that $x= gh(0)$ and $[gh]\in \tilde W$. So already the image of the map $\varphi: \tilde W\to X$, $[g]\to g(0)$ contains $ {\rm locus} _U(W)$. 
	In particular, $(\varphi ^{-1}(x))_{red}= \{[g]\in \tilde  W: g(0)=x \} =\tilde  W^x$ so
	we have
	$$\dim \tilde W^x\ge \dim\tilde  W- \dim {\rm locus} _U(W).$$
	Similarly,  the image of the map $\psi: \tilde  W^x\to X$, $[g]\to g(1),$ contains $ {\rm locus} _U(W,x)\backslash \{x\}$. So for any closed point $y\in {\rm locus} _U(W,x)\backslash \{x\}$ we have
	$$\dim  \{[g]\in \tilde W: g(0)=x, g(1)=y \}\ge \dim \tilde  W^x- \dim {\rm locus} _U(W,x).$$
	If $\dim  \{[g]\in \tilde  W: g(0)=x, g(1)=y \}\ge 2$ then $f$ moves, fixing $0$ and $1$, in a $2$-dimensional family and so its image must move. Since $W$ is unsplit, the family $\tilde W$ is also unsplit and we get a contradiction with  
	Mori's Bend and Break (see \cite[Chapter II, Corollary 5.5.2]{Kollar-Rational}). Therefore
	$$\dim \tilde W =\dim W+1\le \dim {\rm locus} _U(W)+\dim \tilde  W^x\le \dim  {\rm locus} _U(W)+ \dim  {\rm locus} _U(W,x)+1.$$
\end{proof}

\medskip

Every vector bundle $\cE$ on $\PP^1$ splits into a direct sum of the form $\bigoplus \cO_{\PP^1}(a_i)$. We set $\rk ^{+}(\cE):=
\#\{i: a_i>0\}$.

The following proposition is an analogue of \cite[Chapter II, Proposition 3.10]{Kollar-Rational} in the logarithmic case.

\begin{Proposition}\label{log-bound-on-rank}
	Assume that for $f: \PP^1\to X$, the divisor $f^*D$ is supported in at most one point. 
	Let us  fix a closed point $x\in U$ and let $W$ be an irreducible component of 
	$\Hom (\PP^1, X; f^*D\subset D)$ containing $[f]$.
	Then for any closed point $p\in \PP^1\backslash (\{0\}\cup \Supp f^*D)$
	and any  $g: \PP^1\to X$ such that $[g]\in W^x$ we have
	$$\rk dF_{W^x}(p, [g])=\rk ^{+}(g^*T_X(-\log D)).$$
	In particular, we have 
	$$\dim {\rm locus} _U(W,x)\le \max _{\{g:\, [g]\in W^x  \}}  \rk ^{+} \left(g^*T_X(-\log D) \right).$$
\end{Proposition}

\begin{proof}
	Let $I_0\subset \cO_{\PP^1}$ be the ideal sheaf of $\{0\}\subset \PP^1$.
	As in \cite[II.3]{Kollar-Rational} we have
	$$T_{\PP^1\times W^x}\otimes k(p, [g])=T_{\PP^1}\otimes k(p)+H^0(\PP^1,g^*T_X(-\log D) \otimes I_{0})$$	
	and $ dF_{W^x}(p, [g])=dg(p)+\phi (p, g)$, where $dg(p): T_{\PP^1}\otimes k(p)\to g^*T_X\otimes k(p)=T_X\otimes k(g(p))$
	is the differential of $g$ at $p$ and 
	$$\phi (p, g): H^0(\PP^1,g^*T_X(-\log D) \otimes I_{0})\to g^*T_X(-\log D) \otimes I_{0}\otimes k(p)=T_X\otimes k(g(p))$$
	is the evaluation map. Let us set $B=(f^*D)_{red}$ and recall that $g^*D=f^*D$. A standard local computation shows that the differential of $g$ induces the map $T_{\PP^1}(-B)=T_{\PP^1}(-\log B)\to g^*T_X(-\log D)$. So the first claim follows from surjectivity of the evaluation map 
	$$H^0(\PP^1, T_{\PP^1}(-B)\otimes I_0)\to T_{\PP^1}(-B)\otimes I_0\otimes k(p)=T_{\PP^1}\otimes k(p).$$ 
	
	The second part follows from \cite[Chapter III, Proposition 10.6]{Har77}.
\end{proof}

\section{Symplectic and contact structures on smooth Lie algebroids}

In this section we revise symplectic and contact structures in the Lie algeborid setup. We also recall and revise Atiyah extensions in the log case and prove some auxiliary results.

\subsection{Lie algebroids}

Let $X$ be a $k$-variety and let
$\LL = (L, [\cdot, \cdot]_L, \alpha: L\to T_{X})$ be a $k$-Lie algebroid on $X$ such that $L$ is a vector bundle (in this case we say that $\LL$ is \emph{smooth}). Let us recall that $\LL$ is given by a sheaf of $k$-Lie algebras structure on $L$ 
and the anchor map $\alpha$, which is a map of both $\cO_X$-modules and  sheaves of $k$-Lie algebras, and which satisfies Leibniz's rule 
$$[x,fy]=(\alpha (x))(f) y+f[x,y]$$
for all local sections $f\in \cO_X$ and $x, y\in L$ (see, e.g., \cite[Definition 2.1]{Langer-Lie-algebroids}).
In the following we usually omit $k$ in the notation and often write $\LL$ not only for a Lie algebroid but also for the underlying vector bundle $L$.
 
For any Lie algebroid $\LL$ we can consider the canonical de Rham complex 
$(\Omega_{\LL}^{\bullet} , d^{\bullet} _{\LL})$, where $\Omega_{\LL}^m=\bigwedge ^m (L^*)$ (see \cite[2.1]{Langer-Lie-algebroids}). The differential  in this complex is induced by 
$$d_{\LL}\in {\rm Der}_{k}(\cO_X, \Omega_{\LL}^1)= \Hom _{X}(\Omega_{X}^1, \Omega_{\LL}^1),$$ 
corresponding to the dual of the anchor map $\alpha: L\to T_{X}$.  Note that the de Rham complex is a complex of $\cO_X$-modules but the maps in this complex are only $k$-linear.

Let us recall that a map of Lie algebroids is an $\cO_X$-linear map of underlying vector bundles that commutes with anchors and Lie brackets.
It is useful to use the following definition that is modelled on the Lie algebra case:

\begin{Definition} 
	An \emph{abelian extension} of a Lie algebroid $\LL$ by $\LL''$ is a short exact sequence of Lie algebroids
	$$0\to \LL''\to \LL'\to \LL\to 0$$
	in which $\LL''$ is a trivial Lie algebroid on a trivial bundle.
An abelian extension is called \emph{central} if $\LL''$ is generated by global sections whose images belong to the center of the Lie algebra  $H^0(X, \LL')$.
\end{Definition}

Both abelian and central extensions appear naturally in many places in this paper. Note that the word ``abelian''  (or ``central'') can be somewhat misleading as it does not mean that $\LL''$ is abelian (central) in $\LL'$.

\subsection{Symplectic structure}

Let us fix   a $k$-variety $X$ equipped with a smooth Lie algebroid $\LL$. Let $\cL$ be
a line bundle on $X$. 
We say that  a $k$-linear morphism  $\nabla : \cL \to \Omega^1_{\LL} \otimes _{\cO_{X}}\cL$ is an \emph{$\LL$-connection} if it satisfies the Leibniz rule
$$\nabla (fe)=f\nabla (e)+d_{\LL}(f)\otimes e$$
for $e\in \cL$ and $f\in \cO_X$. Then we can define $k$-linear maps $\nabla_i : \Omega^i_{\LL} \otimes _{\cO_{X}}\cL \to \Omega^{i+1}_{\LL} \otimes _{\cO_{X}}\cL  $ by setting 
$$\nabla _i (\omega\otimes e)=d_{\LL}(\omega)\otimes e +(-1)^i\omega\wedge \nabla (e)$$
for $e\in \cL$ and $\omega\in \Omega^i_{\LL} $. If $\nabla_1\circ \nabla =0$ then we say that $\nabla$ is \emph{flat}. If $\nabla$ is flat then the  associated sequence  
$$\cL \mathop{\longrightarrow}^{\nabla} \Omega^1_{\LL} \otimes _{\cO_{X}}\cL\mathop{\longrightarrow}^{\nabla _1} \Omega^2_{\LL} \otimes _{\cO_{X}}\cL\mathop{\longrightarrow}^{\nabla _2} \cdots $$
is a complex. By abuse of notation in the following we write $\nabla$ instead of $\nabla _i$.

Let us also  recall that for a line bundle $\cL$ we have  a canonical bijection between  flat $\LL$-connections on $\cL$, and $\LL$-module structures on $\cL$ compatible with the $\cO_X$-module structure (see, e.g., \cite[Lemma 3.2]{Langer-Lie-algebroids}).

\begin{Definition} 
	Let $\cL$ be a line bundle equipped with a flat $\LL$-connection $\nabla$.
	An \emph{$\cL$-valued symplectic form} on $\LL$ is a section  $\omega\in H^0(X, \Omega_{\LL}^2\otimes _{\cO_X} \cL)$  such that
	\begin{enumerate}
		\item $\omega$ is \emph{closed}, i.e.,  $\nabla (\omega)=0$,
		\item $\omega$ is \emph{nondegenerate}, i.e., the map 
		$\hat\omega:  {\LL} \to \Omega^1_{\LL} \otimes _{\cO_{X}}\cL=\cHom_{\cO_X} (\LL, \cL)$, $x\to \omega (x\wedge  \cdot)$ is an isomorphism.
	\end{enumerate}	
\end{Definition}

The second condition in the above definition implies that the rank of $L$ is even. If $\rk L =2r$ then 
$\omega ^{\wedge r}$ induces an isomorphism 
$ \cO_{X}{\stackrel{\simeq}{\longrightarrow}}\Omega ^{2r}_{\LL}\otimes _{\cO_{X}}\cL^{\otimes r}.$

The above definition generalizes holomorphic symplectic manifolds  that correspond to smooth varieties with an $\cO_X$-valued symplectic form on the canonical Lie algebroid structure on the tangent bundle, where $\cO_X$ is equipped with the canonical connection $d: \cO_X\to \Omega_X^1$.

\medskip

Let $(X,D)$ be an snc pair. The logarithmic tangent bundle  $T_X(-\log D)$ carries a canonical Lie algebroid structure with canonical inclusion   $\alpha: T_X(-\log D)\to T_X$ as the anchor map and with the Lie bracket induced from the Lie bracket on $T_{X}=Der _k (\cO_X, \cO_X)$. By abuse of notation we denote this Lie algebroid again by $T_X(-\log D)$ and we call it the canonical Lie algebroid.
The following definition is a reformulation of the standard definition (see, e.g., \cite[Definition 1.1]{Goto2002}).

\begin{Definition}
A \emph{log symplectic form} on $(X,D)$	is an $\cO_X$-valued symplectic form on the canonical Lie algebroid $T_X(-\log D)$, where $\cO_X$ is equipped with the canonical logarithmic connection $d: \cO_X\to \Omega_X^1 (\log D)$. A pair $(X,D)$ is called \emph{log symplectic} if it admits a log symplectic form.
\end{Definition}

\subsection{Contact structure}

Let us fix   a $k$-variety $X$ equipped with a smooth Lie algebroid $\LL$.
Let $\cL$ be a line bundle on $X$ and let  $ \theta : L\to \cL$ be an $\cO_X$-linear map with kernel $\cF$. 
Although the Lie bracket $[\cdot, \cdot]_{\LL}: {\bigwedge}^2_k L \to L$ is not $\cO_X$-linear, the composition
$${\bigwedge}^2\cF \longrightarrow  {\bigwedge}^2 _kL   \mathop {\longrightarrow}^{[\cdot, \cdot]_{\LL}} L\mathop{\longrightarrow}^{\theta} \cL , \quad x\wedge y\to \theta([x,y]_{\LL})$$
is $\cO_X$-linear (see, e.g., \cite[Lemma 4.4]{Langer-Lie-algebroids}).

\begin{Definition}
An \emph{$\cL$-valued contact structure} on $\LL$ is a surjective $\cO_X$-linear map $ \theta : L\to \cL$
such that ${\bigwedge}^{2}\cF\to \cL,$ $x\wedge y\to \theta ([x,y]_{\LL})$  is nondegenerate, i.e.,  the map $\hat\theta: \cF \to \cF^*\otimes _{\cO_{X}}\cL=\cHom_{\cO_X} (\cF, \cL)$, $x\to \theta ([x, \cdot]_{\LL})$ is an isomorphism.
\end{Definition}

Note that our definition implies that $\rk \cF\ge 2$.
Since  $\hat\theta$ is a non-degenerate skew-symmetric form on $\cF$ with values in $\cL$ (i.e., $\hat\theta=-\hat\theta^*\otimes \id _{\cL}$),  the rank of $\cF$ is even and  ${\bigwedge}^{2}\cF\to \cL$ induces an isomorphism $\det \cF= {\bigwedge}^{2r}\cF{\stackrel{\simeq}{\longrightarrow}}\cL^{\otimes r}$, where $2r$ is the rank of $\cF$. In particular, $\det L\simeq  \cL ^{r+1}$. Since $\cF\subset L$ is not closed under the Lie bracket on $L$, there is no natural Lie algebroid structure on $\cF$ and we cannot say that $\hat \theta$  is a symplectic form. This problem will be resolved in Proposition \ref{general-LeBrun-lemma}, where we show that there exists a natural symplectic form compatible with $\hat \theta$ on some other Lie algebroid.

\begin{Lemma}
For any $\cO_X$-linear map $\omega: \cL^{-1}\to \Omega_{\LL}^1$ with cokernel being a vector bundle $\cF^*$, the composition
$$\varphi: \cL^{-1}\mathop{\longrightarrow}^{\omega}\Omega_{\LL}^1\mathop {\longrightarrow}^{d_{\LL}}\Omega_{\LL}^2\to {\bigwedge}^{2}(\cF^*)=({\bigwedge}^{2}\cF)^*$$
is $\cO_X$-linear.  If $\omega=\theta^*$ then up to a sign $\varphi$ coincides with the dual to ${\bigwedge}^{2}\cF\to \cL,$ $x\wedge y\to \theta ([x,y]_{L})$.	
\end{Lemma}

\begin{proof}
To prove that $\varphi$ is $\cO_X$-linear, we must show that for any local section $s \in \cL^{-1}$ and any local function $f \in \cO_X$ we have
$$\varphi(fs) = f\varphi(s).$$
Let $\eta = \omega(s) \in \Omega^1_{\LL}$. Let us evaluate the 2-form $\varphi(fs)$ on two local sections $x, y \in \cF$, using the standard formula for the Lie algebroid exterior derivative and the Leibniz rule:
\begin{eqnarray*}
	d_{\LL}(f\eta)(x\wedge y)& =& \alpha_x((f\eta)(y)) - \alpha_y((f\eta)(x)) - (f\eta)([x,y]_{\LL})\\
	& =& \big(\alpha_x(f)\eta(y) + f\alpha_x(\eta(y))\big) - \big(\alpha_y(f)\eta(x) + f\alpha_y(\eta(x))\big) - f\eta([x,y]_{\LL})\\
	&=& f \cdot d_{\LL}\eta(x\wedge y) + \alpha_x(f)\eta(y) - \alpha_y(f)\eta(x).
\end{eqnarray*}
Because $\cF^*$ is the cokernel of $\omega$, the 1-form $\eta = \omega(s)$ vanishes on any section of $\cF$. 
Therefore  $\eta(x) =\eta(y) = 0$ and we get
$$d_{\LL}(f\eta)(x\wedge y) = f \cdot d_{\LL}\eta(x\wedge y).$$

To see the second assertion  note that for any $x, y\in \cF$, $s \in \cL^{-1}$
and $\eta= \omega(s)\in \Omega_{\LL}^1$ we have
$$(d_{\LL}\eta )(x\wedge y)=\alpha_x(\eta (y))-\alpha_y(\eta (x))-\eta ([x,y]_{{\LL}})= -\eta ([x,y]_{{\LL}}).$$
\end{proof}
 
\medskip
 
The above lemma allows us to consider 
$$\varphi^{\wedge r}: \cL ^{-r}\to {\bigwedge}^{2r}(\cF^*)=\det (\cF^*)= \det \Omega_{\LL}^1 \otimes \cL.$$
If $\omega=\theta^*$ then $\theta$ is nondegenerate if and  only if the map $\varphi^{\wedge r}$ is an isomorphism.
Let us also note that although $d_{\LL}\omega: \cL^{-1}\to  \Omega_{\LL}^2$ is not $\cO_X$-linear,
$\omega\wedge (d_{\LL}\omega)^{\wedge r}$ gives a well defined $\cO_X$-linear map
from $\cL ^{-(r+1)}$ to $\det \Omega_{\LL}^1$. This map can be identified with $\varphi^{\wedge r}\otimes \id _{\cL^{-1}}$.  A standard computation proving these assertions is left to the reader.

 \begin{Definition}
 An {$\cL$-valued form} $\omega\in H^0(X, \Omega_{\LL}^1\otimes _{\cO_X} \cL)=\Hom\, (\cL^{-1}, \Omega_{\LL}^1)$ is called \emph{contact} if the induced $\cO_X$-linear map $\omega\wedge (d_{\LL}\omega)^{\wedge r}: \cL ^{-(r+1)}\to \det \Omega_{\LL}^1$ is an isomorphism. 
 \end{Definition}
 
 Note that the above remarks show that giving an $\cL$-valued contact structure on $\LL$ is equivalent to giving an 
 {$\cL$-valued contact form}.

\medskip

Sometimes it is also convenient to consider contact structures defined only at the generic point of $X$:

\begin{Definition}
A \emph{generically $\cL$-valued contact structure}  on $\LL$ is a non-zero  $\cO_X$-linear map $ \theta :\LL \to \cL$ for which there exists a non-empty Zariski open subset $U\subset X$ such that $\theta|_U$ is an $\cL|_U$-valued contact structure on $\LL|_U$.
\end{Definition}

If $\theta: \LL\to \cL$ is a generically $\cL$-valued contact structure and $\omega: \cL^{-1}\to \Omega _{\LL}^1$ is its dual  then $\omega\wedge (d\omega)^{\wedge r}$ defines a non-zero section of $\det \Omega_{\LL}^1\otimes \cL ^{r+1}$. In particular, $\det \Omega_{\LL}^1=\cL^{-(r+1)}(B)$ for some effective divisor $B$. $\omega$ defines an $\cL$-valued contact structure on $\LL$ if and only if $B=0$.

\subsection{Jet bundles and the Atiyah extension for Lie algebroids}\label{Atiyah-classes}

Let us recall Atiyah classes defined by Lie algebroids after \cite[4.2.5 and 8.1]{Calaque-Van-den-Bergh}.
If $\LL$ is a Lie algebroid then we define the \emph{first $\LL$-jet $k$-algebra bundle}  $J^1_{\LL}$ as 
$\cHom _{\cO_X}({L}\oplus \cO_X, \cO_X)=\Omega_{\LL}^1\oplus \cO_X$ with a natural commutative $k$-algebra structure given by $(\eta _1,f_1)\cdot (\eta_2, f_2)=(f_2\eta_1+f_1\eta_2, f_1f_2)$. $J^1_{\LL}$ comes with 
an $\cO_X$-$\cO_X$-bimodule structure, where the left $\cO_X$-module structure is
given by the monomorphism
$$\alpha_1: \cO_X\to J^1_{\LL}, \quad f\to (0,f)$$
and the right $\cO_X$-module structure is
given by the monomorphism
$$\alpha_2: \cO_X\to J^1_{\LL}, \quad f\to (d_{\LL}(f),f).$$

Since $\Omega_{\LL}^1$ is an ideal in $J^1_{\LL}$, we have the short exact sequence 
$$0\to \Omega_{\LL}^1\longrightarrow J^1_{\LL}\longrightarrow \cO_X\to 0$$ 
of $J^1_{\LL}$-modules. In general, this sequence  is not split although it is split both as a sequence of left $\cO_X$-modules and as a sequence of right $\cO_X$-modules.

For a vector bundle $\cE$ on $X$, we define its \emph{first $\LL$-jet bundle } as $J^1_{\LL}(\cE)= J^1_{\LL}\otimes _{\cO_X}\cE$. It has a natural structure of a left $J^1_{\LL}$-module and the above short exact sequence induces the short exact sequence
$$0\to \cE \otimes _{\cO_X} \Omega_{\LL}^1\simeq \Omega_{\LL}^1\otimes _{\cO_X}\cE\longrightarrow J^1_{\LL}(\cE) \longrightarrow \cE\to 0$$  
of left $J^1_{\LL}$-modules (hence also of left $\cO_X$-modules). 
This sequence is called the \emph{$\LL$-Atiyah extension} of $\cE$ and it
defines so called \emph{$\LL$-Atiyah class} of $\cE$:
$$\alpha_{\cE} \in \Ext ^1_X (\cE, \Omega_{\LL}^1\otimes _{\cO_X}\cE)= H ^1 (X, \Omega_{\LL}^1\otimes _{\cO_X}\cEnd \cE ).$$ 
This sequence is split as a sequence of sheaves of abelian groups (or even of sheaves of $k$-vector spaces) with the canonical splitting $j^1: \cE \to  J^1_{\LL}(\cE) $
induced by $\alpha_2$. Existence of an $\cO_X$-linear splitting of this sequence (treated as a sequence of left $\cO_X$-modules) is equivalent to existence of an $\LL$-connection $\cE\to \Omega_{\LL}^1\otimes _{\cO_X} \cE$.

\medskip

If $\cE=\cL$  is a line bundle, tensoring the {$\LL$-Atiyah extension} of $\cL$ by $\cL^{-1}$ (on the left) gives the extension
$$0\to \Omega_{\LL}^1\longrightarrow \cA_{\LL}(\cL) \longrightarrow \cO_X \to 0$$ 
in which $ \cA_{\LL}(\cL)= \cL^{-1} \otimes _{\cO_X} J^1_{\LL}(\cL)$. Sometimes we also call this extension the {$\LL$-Atiyah extension} of $\cL$.
By construction  we have a canonical pushout diagram
$$\xymatrix{
	0\ar[r]&\Omega_X^1\ar[r]\ar[d]^{\alpha^*}& \cA_{T_X}(\cL)\ar[d]\ar[r]& \cO_X\ar[r]\ar@{=}[d]&0\\
	0\ar[r]&\Omega_{\LL}^1\ar[r]& \cA_{\LL}(\cL)\ar[r]& \cO_X\ar[r]&0\\
}$$
in which $\alpha^*$ is the dual of the anchor map $\alpha: L\to T_{X}$.
The $T_X$-Atiyah class of $\cL$ can be identified with the first Chern class $c_1(\cL)\in H^1(X,\Omega_{X}^1)$
(sometimes people use a different normalization in the complex case).
The $\LL$-Atiyah  class  is therefore the image of $c_1(\cL)$ under the canonical map $H^1(X,\Omega_{X}^1)\to H^1(X,\Omega_{\LL}^1)=\Ext^1(\cO_X, \Omega_{\LL}^1)$ induced by $d_{\LL}$. We denote the corresponding class by $c_1^{\LL}(\cL)$.

In the following we need also a generalization of the Atiyah extension for a collection $\underline{\cL}=(\cL_1,...,\cL_s)$ of line bundles on $X$. Then we define the \emph{multi-$\LL$-Atiyah extension of $\underline{\cL}$}
	$$0\to \Omega_{\LL}^1\to  \cA_{\LL} ({\underline{\cL}})\to \cO_X^{\oplus s}\to 0$$
	as the pushout of the direct sum of Atiyah extensions for line bundles $\cL_1,...,\cL_s$
	$$0\to (\Omega_{\LL}^1)^{\oplus s}\to \cA_{\LL}(\cL_1)\oplus ...\oplus\cA_{\LL}(\cL_s) \to \cO_X^{\oplus s}\to 0$$
	by the codiagonal map $(\Omega_{\LL}^1)^{\oplus s}\to \Omega_{\LL}^1$.
Clearly, the extension class of the multi-Atiyah extension in $ \Ext^1(\cO_X^{\oplus s},\Omega_{\LL}^1)= H^1(X,\Omega_{\LL}^1)^{\oplus s}$ is given by  $\oplus_{i=1}^sc_1^{\LL}(\cL_i)$.

\medskip

\subsection{Differential operators for Lie algebroids}\label{Atiyah-classes2}

Let $\cD _{\LL}$ be the universal enveloping algebra of the Lie algebroid $\LL$ and let $\cE$ be a vector bundle on $X$.
Then we define the sheaf $\cD_{\LL}(\cE )$  of $\LL$-differential operators acting on $\cE ^{*}$ as a subsheaf of $\cEnd _k (\cE^{*})$ equal to $\bigcup _{i\ge 0 } \cD_{\LL}^{\le i}(\cE )$, where
$ \cD_{\LL}^{\le -1}(\cE ):= 0$ and
$$\cD_{\LL}^{\le i}(\cE ):= \{ x\in \cEnd _k (\cE^*) : \,   xf-fx\in \cD_{\LL}^{\le i-1}(\cE ) \hbox{ for all }f\in \cO_X \}  $$
for $i\ge 0$.
Then $\cD_{\LL}(\cO_X )=\cD_{\LL}$ and we have a canonical isomorphism of sheaves of rings
$$\cE^{*}\otimes _{\cO_X} \cD_{\LL}\otimes _{\cO_X}\cE \simeq \cD_{\LL}(\cE ) $$
given by  
$$(s\otimes x\otimes u) (t)= x(\langle u, t\rangle )s,$$
where $s,t \in \cE^{*}$, $u\in \cE$ and $x\in \cD_{\LL}$.

Note that $\cD_{\LL}(\cE )$ is a sheaf of rings with the ring structure induced from 
$\cEnd _k (\cE^*)$ and $\cD_{\LL}^{\le 0}(\cE )\simeq \cEnd _{\cO_X} (\cE^*)$ is a sheaf of subrings of $\cD_{\LL}(\cE )$ containing $\cO_X$.

If $\cL$ is a line bundle then $\cD_{\LL}^{\le 0}(\cL )\simeq \cO_X$ and 
$\cD_{\LL}^{\le 1}(\cL )$ has a canonical Lie algebroid structure, defined by the Lie bracket coming from the associative algebra structure on $\cD_{\LL}(\cL )$
and the anchor map $$\cD_{\LL}^{\le 1}(\cL ) \longrightarrow \cD_{\LL}^{\le 1}(\cL )/ \cD_{\LL}^{\le 0}(\cL ) \simeq \LL \to T_X.$$
Moreover, we have a central extension of Lie algebroids
$$0\to \cO_X \longrightarrow \cD_{\LL}^{\le 1}(\cL ) \longrightarrow \LL \to 0,$$
since the image of $1\in H^0 (X, \cO_X)$  is central, acting as the identity on $\cL^{-1}$. This extension can be identified with the dual to the {$\LL$-Atiyah extension} of $\cL$ and its extension class in $\Ext ^1(\LL, \cO_X)= H^1(X, \Omega ^1_{\LL})$ is given by $-c_1^{\LL}(\cL)$.
In particular, we have $\cD_{\LL}^{\le 1}(\cL )= (\cA_{\LL}(\cL))^*$,
i.e., $\cA_{\LL}(\cL)= \Omega_{\cD_{\LL}^{\le 1}(\cL )}^1$.

If $\underline{\cL}=(\cL_1,...,\cL_s)$ is a collection of line bundles then we can construct
$$0\to \cO_X ^{\oplus s} \longrightarrow \cD_{\LL}^{\le 1}(\underline{\cL} ) \longrightarrow \LL \to 0$$ 
as the pullback of the direct sum of
$$0\to \cO_X ^{\oplus s}\longrightarrow \cD_{\LL}^{\le 1}({\cL _1} )\oplus ...\oplus \cD_{\LL}^{\le 1}({\cL _s} )\longrightarrow \LL ^{\oplus s}\to 0$$ 
by the diagonal map $\LL\to \LL ^{\oplus s}$. The class of this extension is given by $-\oplus_{i=1}^sc_1^{\LL}(\cL_i)$.
Note that $\cD_{\LL}^{\le 1}(\underline{\cL} )$ has a natural structure of a Lie algebroid for which the above extension is a central extension of Lie algebroids. Under the identification $\cD_{\LL}^{\le 1}(\underline{\cL} )=(\cA_{\LL} ({\underline{\cL}}))^*$, the anchor map corresponds to the dual of 
$\Omega_X^1\to \Omega_{\LL}^1\to \cA_{\LL} ({\underline{\cL}})$.

\subsection{Relation between contact and symplectic structures}

Below we describe generalization of the classical process of symplectification of a contact manifold (see, e.g., \cite[Appendix 4 E]{Arnold-Mathematical-methods})  to arbitrary contact structures on Lie algebroids.
We follow the approach of \cite[Proposition 2.4]{LeBrun} and \cite[Lemma 1.4]{Beauville-Fano-contact}
(see also \cite[2.1 and 2.2]{KPSW}).  Classically this gives a relation between contact structures on a variety and $\GG_m$-invariant symplectic structures on the associated principal $\GG_m$-bundle. The Lie algebroid point of view makes this a relation between contact and symplectic structures on the same variety but for different Lie algebroids. 

First, let us make the following easy observation.

\begin{Lemma} 
	A line bundle $\cL$ has a canonical flat $\cD_{\LL}^{\le 1}(\cL )$-connection $\nabla ^{\sf can}$.
\end{Lemma}

\begin{proof}
	By construction we have a canonical map
	$$j^1:  \cL\to J^1_{\LL}(\cL)\simeq \cA_{\LL}(\cL) \otimes _{\cO_X} \cL =\Omega_{\cD_{\LL}^{\le 1} (\cL )}^1  \otimes _{\cO_X} \cL ,$$
	giving a $\cD_{\LL}^{\le 1}(\cL )$-connection $\nabla ^{\sf can}$. The connection is flat as $\cL$ is a line bundle.
	Alternatively, we can obtain the required connection  $\nabla ^{\sf can}$ considering the action of $\LL$-differential operators on $\cL$, dual to the canonical action on  $\cL^{-1}$. 	
\end{proof}

\begin{Proposition}\label{general-LeBrun-lemma}
Let $\LL$ be a (smooth) $k$-Lie algebroid on a $k$-variety $X$ and let $\cL$ be a fixed line bundle on $X$. There exists a natural bijection between $\cL$-valued contact structures on  $\LL$ and $\cL$-valued symplectic forms on the Lie algebroid $\cD_{\LL}^{\le 1}(\cL )$, where $\cL$ is equipped with the canonical flat $\cD_{\LL}^{\le 1}(\cL )$-connection $\nabla ^{\sf can}$.
\end{Proposition}

\begin{proof}
Let $\MM:=\cD_{\LL}^{\le 1}(\cL)$ and let $\pi:\MM\to\LL$ be the
canonical map. Let
\( e=\mathrm{id}_{\cL^{-1}}\in H^0(X,\MM) \)
be the canonical section spanning $\ker(\pi)\simeq\cO_X$. The section
$e$ is central in $\MM$ and acts as the identity on $\cL^{-1}$.
Consequently, it acts as minus the identity on $\cL$.
Therefore, for every \(
\omega\in H^0(X,\Omega_{\MM}^j\otimes_{\cO_X}\cL), \)
the covariant Lie derivative along $e$ satisfies
	\[ 
\mathbf{L}_e \omega = \nabla^{\sf can}_e \circ \omega - \sum_{i=1}^k \omega(\dots, [e, \cdot]_{\MM}, \dots) =- \omega .
\]	
 So Cartan's magic formula 
$\mathbf{L}_e = \nabla^{\sf can} \circ \iota_e + \iota_e \circ \nabla^{\sf can}$ implies that
\begin{equation}\label{Cartan}
	-\omega =  \nabla^{\sf can} (\iota_e \omega) + \iota_e( \nabla^{\sf can} \omega)
\end{equation}
for all  $\omega \in H^0 (X, \Omega^j_{\MM}\otimes _{\cO_X}\cL )$.

\emph{Step 1.}
Let $\omega \in H^0(X, \Omega^2_{\MM }\otimes _{\cO_X}\cL )=\Hom (\bigwedge ^2\MM, \cL)$ be an $\cL$-valued symplectic form on $\MM$.  	
Because $\nabla^{\sf can} \omega = 0$, equation (\ref{Cartan}) reduces to
	\[
	\omega = \nabla^{\sf can} \tilde{\theta},\] where \(\tilde{\theta} = -\iota_e \omega \in H^0(X, \Omega^1_{\MM }\otimes _{\cO_X}\cL )=\Hom (\MM, \cL).
	\)
Skew-symmetry implies $\tilde{\theta}(e) =- \omega(e, e) = 0$, so $\tilde{\theta}$ vanishes on $\ker(\pi)$. Thus $\tilde{\theta}$ uniquely descends to a section of $ \Omega^1_{\LL}\otimes _{\cO_X}\cL$, representing an $\cO_X$-linear map $\theta: \LL \to \cL$. Nondegeneracy of $\omega$ implies $\tilde{\theta}_p \neq 0$ at every point $p \in X$, ensuring that $\theta$ is surjective.
Let $\cF := \ker(\theta) \subset \LL$ and  $\tilde \cF : = \pi^{-1}(\cF) = \ker(\tilde{\theta}) \subset \MM$. For any $x, y \in \cF$ with lifts $\tilde{x}, \tilde{y} \in \tilde \cF$ we have
		\[
		\omega (\tilde{x}, \tilde{y}) = (\nabla^{\sf can}\tilde{\theta})(\tilde{x}, \tilde{y}) = \nabla^{\sf can}_{\tilde{x}}(\tilde{\theta}(\tilde{y})) - \nabla^{\sf can}_{\tilde{y}}(\tilde{\theta}(\tilde{x})) - \tilde{\theta}([\tilde{x}, \tilde{y}]_{\MM}) = -\theta([x, y]_{\LL}).
		\]
		Because $\omega (e, \tilde \cF) = \tilde{\theta}(\tilde \cF) = 0$, the radical of $\omega|_{\tilde \cF}$ is precisely $\cO_X e$. Nondegeneracy of $\omega$ on $\MM$ forces the induced bilinear form $\omega|_{\tilde \cF}$ to be nondegenerate on $\tilde \cF / (\cO_X \cdot e) \cong \cF$. Thus, $x \wedge y \mapsto \theta([x, y]_{\LL})$ is nondegenerate on $\cF$, making $\theta$ a contact structure.

\emph{Step 2.}
	Given an $\cL$-valued contact structure $\theta: \LL \to \cL$, set $\tilde{\theta} = \pi^* \theta ^*= \theta\circ\pi\in H^0(X,\Omega^1_{\MM }\otimes _{\cO_X}\cL )$ and define
	\[
	\omega := \nabla^{\sf can} \tilde{\theta} \in H^0(X, \Omega^2_{\MM }\otimes _{\cO_X}\cL )
	\]
 Flatness of $\nabla^{\sf can}$ implies that  $\nabla^{\sf can}\omega = 0$,
 so $\omega$ is closed. 
 
 To prove that $\omega$ is nondegenerate, let $\cF := \ker(\theta)$ and  $\tilde \cF : = \pi^{-1}(\cF) = \ker(\tilde{\theta})$.
 Since $e$ is central in $\MM$, $\tilde\theta(e)=0$, and
 $\nabla^{\sf can}_e$ acts as minus the identity on $\cL$, we have, for
 every local section $v$ of $\MM$,
 \[
 \omega(e,v)
 =
 \nabla^{\sf can}_e\bigl(\tilde\theta(v)\bigr)
 -\nabla^{\sf can}_v\bigl(\tilde\theta(e)\bigr)
 -\tilde\theta([e,v]_{\MM})
 =
 -\tilde\theta(v).
 \]
 
 Suppose that a local section $v$ of $\MM$ belongs to the radical of
 $\omega$. Pairing $v$ with $e$, we obtain
 \[
 0=\omega(e,v)=-\tilde\theta(v),
 \]
 and hence $v\in\tilde{\cF}$. So for every local section
 $w\in\tilde{\cF}$ we have
 $$
 	0=\omega(v,w)
 	=
 	\nabla^{\sf can}_v\bigl(\tilde\theta(w)\bigr)
 	-\nabla^{\sf can}_w\bigl(\tilde\theta(v)\bigr)
 	-\tilde\theta([v,w]_{\MM})
 	=
 	-\theta\bigl([\pi(v),\pi(w)]_{\LL}\bigr).
 $$
 Since $\pi:\tilde{\cF}\to\cF$ is surjective and the pairing
 \(
 \cF\times\cF\longrightarrow\cL,
 \) \(
 (x,y)\longmapsto\theta([x,y]_{\LL}),
 \)
 is nondegenerate, it follows that $\pi(v)=0$. Consequently,
 $v=f e$ for some local function $f$.
 
 Finally, for every local section $w$ of $\MM$,
 \[
 0=\omega(v,w)
 =f\,\omega(e,w)
 =-f\,\widetilde\theta(w).
 \]
 The map $\tilde\theta=\theta\circ\pi:\MM\to\cL$ is surjective, so
 locally one can choose $w$ such that $\tilde\theta(w)$ generates
 $\cL$. Hence $f=0$, and therefore $v=0$. Thus $\omega$ is
 nondegenerate and hence it is an $\cL$-valued symplectic form on $\MM$.
 
\emph{Step 3.}
Let us check that the two constructions are mutually inverse. If
$\tilde\theta=\pi^*\theta^*=\theta\circ\pi$, then
$\iota_e\tilde\theta=0$ and $\mathbf L_e\tilde\theta
=-\tilde\theta$. Hence Cartan's formula gives
\[
-\iota_e\bigl(\nabla^{\sf can}\widetilde\theta\bigr)
=
-\mathbf L_e\widetilde\theta
+
\nabla^{\sf can}(\iota_e\widetilde\theta)
=
\widetilde\theta.
\]
Conversely, if $\omega$ is closed, equation~\eqref{Cartan} gives
\[
\nabla^{\sf can}(-\iota_e\omega)=\omega.
\]
Thus the mutually inverse correspondences are determined by
\[
\theta\circ\pi=-\iota_e\omega,
\qquad
\omega=\nabla^{\sf can}(\pi^*\theta^*).
\]	
\end{proof}

\medskip

\begin{Remark}
	As in the proof of \cite[Proposition 2.4]{LeBrun} we have the following
	commutative diagram:
	$$\xymatrix{
		\cF\ar[d]_{\hat \theta}^{\simeq} \ar[r] &\LL&\cD_{\LL}^{\le 1}(\cL )\ar[l]\ar[d]_{\hat \omega_{\cL}}^{\simeq }\\
		\cF^*\otimes \cL&\Omega^1_{\LL}\otimes \cL\ar[l]\ar[r]&\Omega ^1_{\cD_{\LL}^{\le 1}(\cL )}\otimes \cL\\
	}
	$$
	This diagram shows that the induced map $\Omega^1_{\LL}\otimes \cL\to \LL$ (using inverses of the isomorphisms) is compatible with both contact form on $\LL$ and symplectic form on $\cD_{\LL}^{\le 1}(\cL )$.
\end{Remark}

\medskip

A relation between the above point of view and the  classical one (see, e.g., \cite[Lemma 1.4]{Beauville-Fano-contact}) can be explained as follows. Let
$\pi:\widetilde X:=\Spec _X (\bigoplus _{m\in \ZZ} \cL^{m})\to X$
be the principal $\mathbb G_m$-bundle associated with $\cL^{-1}$, and
let $\alpha:\LL\longrightarrow T_X$ denote the anchor map. Let us recall that one defines the \emph{the inverse-image Lie algebroid} as
\[
\pi^{!}\LL
:=
T_{\widetilde X}\times_{\pi^*T_X}\pi^*\LL
=
\left\{
(v,\ell)\in T_{\widetilde X}\oplus \pi^*\LL
\ \middle|\
d\pi(v)=\pi^*\alpha(\ell)
\right\}.
\]
The anchor map is given by the first projection
$\tilde\alpha:\pi^{!}\LL\longrightarrow T_{\tilde X}$,
$(v,\ell)\longmapsto v$.
To describe the Lie bracket, write two local sections of $\pi^{!}\LL$ as
\[
u=
\left(v,\sum_i f_i\pi^*\ell_i\right),
\qquad
u'=
\left(w,\sum_j g_j\pi^*m_j\right),
\]
where $\ell_i,m_j$ are local sections of $\LL$.
Then we define the Lie bracket on $\pi^{!}\LL$ by
\[
	[u,u']_{\pi^{!}\LL}
	:=
	\Bigg(
	[v,w],\,
	\sum_{i,j}f_ig_j\pi^*[\ell_i,m_j]_{\LL}\\
	+\sum_jv(g_j)\pi^*m_j
	-\sum_iw(f_i)\pi^*\ell_i
	\Bigg).
\]
This definition is independent of the chosen presentations of the
local sections.

Since $\pi$ is smooth, there is a canonical short exact sequence of vector bundles
\[
0\longrightarrow T_{\tilde X/X}
\longrightarrow \pi^{!}\LL
\longrightarrow \pi^*\LL
\longrightarrow0.
\]
For a principal $\mathbb G_m$-bundle, the relative tangent bundle is
canonically trivialized by the Euler vector field $\varepsilon$.
Consequently, the above sequence can be written as
\[
0\longrightarrow
\cO_{\tilde X}
\longrightarrow \pi^{!}\LL
\longrightarrow \pi^*\LL
\longrightarrow 0,
\]
where the first non-zero map maps $1$ to $-\epsilon$.
Since $\cD_{\LL}^{\leq1}(\cL)$ acts on $\cL^{-1}$, every local section
$v$ of $\cD_{\LL}^{\leq1}(\cL)$ induces a
$\mathbb G_m$-invariant vector field $\tilde v$ on $\tilde X$.
This gives a canonical isomorphism
\[
\Phi:
\pi^*\cD_{\LL}^{\leq1}(\cL)
\mathop{\longrightarrow}^{\simeq}
\pi^!\LL,
\qquad
\pi^*v\longmapsto
\bigl(\tilde v,\pi^*\sigma(v)\bigr),
\]
where $\sigma:\cD_{\LL}^{\leq1}(\cL)\longrightarrow\LL$ is the canonical projection. 
The pullback of the canonical inclusion $\cO_X\to\cD_{\LL}^{\leq1}(\cL)$ coincides with
$\cO_{\tilde X} \to\pi^!\LL$ mapping $1$ to $-\epsilon$.

The canonical flat $\cD_{\LL}^{\leq1}(\cL)$-connection on $\cL^{-1}$ pulls back, through
$\Phi$, to a flat $\pi^!\LL$-connection on $\pi^*\cL^{-1}$. By the projection formula $\pi_*(\pi^* \cL^{-1})\simeq \bigoplus _{m\in \ZZ} \cL^{m}$, so $\pi^*\cL^{-1}$ has a canonical section  corresponding to $1\in H^0(X, \cL^{0}\simeq \cO_X)$. This section trivializes $\pi^*\cL^{-1}$ and it
is flat for the pullback connection. Consequently, $\cO_{\tilde X}$ has a canonical flat $\pi^!\LL$-connection corresponding to the pullback of $\nabla ^{\sf can}$. In the following a \emph{symplectic form on $\pi^!\LL$} is 
an $\cO_{\tilde X}$-valued symplectic form, where $\cO_{\tilde X}$ is equipped with this $\pi^!\LL$-connection.
Such form is called \emph{$\GG_m$-equivariant} if it is homogeneous of weight $1$ (cf. \cite[1.3]{Beauville-Fano-contact}).

\begin{Corollary}\label{general-LeBrun-lemma-v2} 
There exists a natural bijection between $\cL$-valued contact structures on  $\LL$ and $\GG_m$-equivariant 
symplectic forms on the Lie algebroid $\pi^!\LL$.
\end{Corollary}

\begin{proof} 
It is sufficient to note that $\GG_m$-equivariant symplectic forms on the Lie algebroid $\pi^!\LL$
correspond to $\cL$-valued symplectic forms on the Lie algebroid $\cD_{\LL}^{\le 1}(\cL )$ and use Proposition \ref{general-LeBrun-lemma}.
\end{proof}

\medskip

\begin{Corollary}\label{general-LeBrun-lemma-v3} 
Let $(X,D)$ be an snc pair. Let $\cL$ be a line bundle on $X$ and let $\pi:\widetilde X\to X$ be the principal $\mathbb G_m$-bundle associated with $\cL^{-1}$. Let us set \(\LL:=\cD_{T_X(-\log D)}^{\leq1}(\cL). \)
Then we have canonical bijections between:
\begin{enumerate}
	\item $\cL$-valued contact structures on $(X,D)$,
	\item $\cL$-valued symplectic forms on $\LL$, where $\cL$ is equipped with the canonical flat ${\LL}$-connection,
	\item $\GG_m$-equivariant  log symplectic forms on $(\tilde X,\pi ^{-1}(D))$.
\end{enumerate}
\end{Corollary}

\begin{proof}
Bijection between (1) and (2) follows from Proposition \ref{general-LeBrun-lemma}.
Since $\pi$ is smooth there is a canonical isomorphism of $\mathbb G_m$-equivariant Lie algebroids
\[
T_{\tilde X}\bigl(-\log \pi^{-1}(D)\bigr)
\simeq \pi^{!}T_X(-\log D).
\]
So Corollary \ref{general-LeBrun-lemma-v2} gives a bijection between (1) and (3).
\end{proof}

\subsection{Atiyah extensions on snc pairs}\label{subsection-Atiyah-diff}

Let $(X,D)$ be a projective snc pair. 
If $\cL$ is a line bundle on $X$ then we usually write $\cA_{\cL}$ for the Atiyah bundle $ \cA_{T_X}(\cL)$
and $\cB_{\cL}$ for the logarithmic Atiyah bundle $ \cA_{T_X(-\log D)}(\cL)$. 
In the following we write $c_1^{\log} (\cL)$ for  $c_1^{T_X(-\log D)} (\cL)$, i.e., for
the image of $c_1(\cL)$ under the canonical map $H^1(X,\Omega_X^1)\to H^1(X, \Omega_{X}^1(\log D))$. 
Similarly, if  $\underline{\cL}$ is  a collection of line bundles  then we write
$\cA_{\underline{\cL}}$ for the {multi-Atiyah bundle of $\underline{\cL}$} and 
 $\cB_{\underline{\cL}}$ for the logarithmic  multi-Atiyah bundle.

Let $(Y, B)$ be  another projective snc pair and let $\varphi : X\to Y$ be a morphism.
The canonical surjection $\Omega_X^1\to \Omega ^1_{X/Y}$ allows us to consider 
the relative  Atiyah extension
$$0\to \Omega_{X/Y}^1\to \cA_{\cL/Y}\to \cO_X\to 0$$
defined by the pushout of the standard Atiyah extension. The class of this extension in $H^1(X, \Omega_{X/Y}^1)=\Ext ^1(\cO_X, \Omega_{X/Y}^1)$ is the image of $c_1(\cL)$ under the induced map $H^1(X, \Omega_{X}^1)\to H^1(X, \Omega_{X/Y}^1)$.

Let us write $D=\sum _{i\in I}D_i$ as a sum of smooth irreducible $D_i$. Then  for any $J\subset I$
we write $D_J$ for $\bigcap _{j\in J} D_j$.
Let us assume that $\varphi$ is smooth and the restriction of $\varphi$ to each $D_J$ is also smooth. 
Then we can decompose $D$ into the horizontal part $D^h$ corresponding to components that dominate $Y$ and the vertical part $D^v$ that consists of the remaining components. Let us assume that $D^v=\varphi ^{-1}(B)$.

Then $\Omega ^1_{X/Y}(\log D^h)=\Omega ^1_{X}(\log D^h)/ \varphi^*\Omega^1_Y$ is locally free and it fits into the short exact sequence
$$0\to \Omega ^1_{X/Y}\to \Omega ^1_{X/Y}(\log D^h) \to {\bigoplus}\cO_{D^h_i}\to 0 $$
in which $D^h_i$ are the irreducible components of $D^h$. Let us also note that $\Omega ^1_{X/Y}(\log D^h)$
is canonically isomorphic to $\Omega ^1_{X}(\log D)/ \varphi^*\Omega^1_Y(\log B)$. 
So as above we can define the relative logarithmic Atiyah extension of $(X,D)/(Y,B)$ as the extension
$$0\to\Omega ^1_{X/Y}(\log D^h) \to \cB_{\cL/Y}\to \cO_X\to 0,$$
whose class in  $H^1(X, \Omega ^1_{X/Y}(\log D^h))$ is the image of $c_1^{\log} (\cL)$  under the canonical map $$H^1(X,\Omega ^1_{X}(\log D))\to H^1(X, \Omega ^1_{X/Y}(\log D^h)).$$ 
This extension can be identified with the relative logarithmic Atiyah extension of $\cL$ for $(X,D^h)/Y$.
By construction we have the following commutative diagram:
$$\xymatrix{
	&0\ar[d]&0\ar[d]&0\ar[d]&\\
	0\ar[r]&\Omega_{X/Y}^1\ar[r]\ar[d]& \cA_{\cL/Y}\ar[d]^{\alpha}\ar[r]& \cO_X\ar[r]\ar@{=}[d]&0\\
	0\ar[r]& \Omega ^1_{X/Y}(\log D^h)\ar[d]\ar[r]& \cB_{\cL/Y} \ar[r]\ar[d]&  \cO_X \ar[r]\ar[d]&0\\
	0\ar[r] &{\bigoplus}\cO_{D^h_i}\ar[r]^{\simeq}\ar[d]&\coker \alpha \ar[r]\ar[d]&0&\\
	&0&0
	&&
}$$	
Note also that our definitions imply that the following diagram is commutative:
$$
\xymatrix{
	&0\ar[d]&0\ar[d]&0\ar[d]&\\
	0\ar[r]&\varphi ^*\Omega^1_Y(\log B)\ar@{=}[r]\ar[d]& \varphi ^*\Omega^1_Y(\log B)\ar[d]\ar[r]& 0\ar[r]\ar[d]&0\\
	0\ar[r]&\Omega ^1_{X}(\log D)\ar[d]\ar[r]& \cB_{\cL} \ar[r]\ar[d]&  \cO_X \ar[r]\ar@{=}[d]&0\\
	0\ar[r] & \Omega ^1_{X/Y}(\log D^h)\ar[r]\ar[d]&\cB_{\cL/Y}\ar[r]\ar[d]&\cO_X\ar[r]\ar[d]&0\\
	&0&0&0&\\
}
$$
Therefore we get the following lemma:  

\begin{Lemma}\label{seq:A-B-D^h}
	Under the above assumptions there exist canonical short exact sequences:
	$$0\to\cA_{\cL/Y}\to \cB_{\cL/Y}\to {\bigoplus}\cO_{D^h_i}\to 0$$
and 
$$0\to\varphi ^*\Omega^1_Y(\log B)\to \cB_{\cL}\to \cB_{\cL/Y}\to 0.$$
\end{Lemma}

\medskip

For future reference we also note the following easy lemma:

\begin{Lemma}\label{Atiyah-bundles-sequence}
	Let  
	$$0\to\cO_{X}^{\oplus m}\to \LL\to  T_{X}(-\log D)\to 0 $$	
	be a short exact sequence of Lie algebroids
	and let $\cL$ be a line bundle on $X$. Then we have an induced short exact sequence
	$$0\to \cB_{\cL} \to  \cA_{\LL}(\cL)\to \cO_{X}^{\oplus m}\to 0.$$
\end{Lemma}

\begin{proof}
	The assertion follows from the following commutative diagram
	$$\xymatrix{
		&0\ar[d]&0\ar[d]&0\ar[d]&\\
		0\ar[r]& \Omega_{X}^1(\log D)\ar[r]\ar[d]& \Omega_{\LL}^1\ar[d]\ar[r]& \cO_X^{\oplus m}\ar[r]\ar[d]^{\simeq}&0\\
		0\ar[r]&\cB_{\cL}\ar[d]\ar[r]^-{\beta}& \cA_{\LL}(\cL) \ar[r]\ar[d]&  \coker \beta \ar[r]\ar[d]&0\\
		0\ar[r]&\cO_X\ar@{=}[r]\ar[d]&\cO_X\ar[r]\ar[d]&0&\\
		&0&0&&\\
	}$$	
	(cf. Subsection \ref{Atiyah-classes}).
\end{proof}

\subsection{Hyperplane subbundles in projective bundles}\label{hyperplane-subbundles}

Let us consider the following situation. Let
$Y$ be a smooth projective variety and let $\pi: X=\PP(\cE)\to Y$ be a projectivization of a vector bundle $\cE$ of rank $m+1$. Let us set  $\cL:= \cO_{\PP (\cE)}(1)$. In $X$ we consider a smooth divisor $D$, which is also a projective bundle over $Y$. Assume that over some point $y\in Y$ we have $(X_y, D_y)\simeq (\PP^m, H)$, where $H$ is a hyperplane.
Then $\cL(-D)|_{X_y}\simeq \cO_{X_y}$ and since $\Pic \PP(\cE)\simeq \Pic Y\times \ZZ$, there exists some line bundle $\cM$ on $Y$ such that  
$$ \cL(-D)\simeq \pi^* \cM.$$
Let us set $\cG:=\pi_* (\cL|_D)$. This is a vector bundle and $(D, \cL|_D)\simeq (\PP(\cG), \cO_{\PP(\cG)} (1))$.
By the projection formula we have $\pi_{*}(\cL(-D))=\cM$ and  $R^1\pi_{*}(\cL(-D))=0$. 
Pushing down the standard short exact sequence
$$0\to \cL(-D)\to \cL \to \cL|_{D}\to 0,$$
we get the short exact sequence
$$0\to \cM\to \cE \to \cG\to 0.$$
In particular, we have $\det \cG \simeq \det \cE\otimes \cM^{-1}$.

\begin{Lemma}\label{relative-cotangent-bundle}
	$\Omega_{X/Y}^1(\log D) $ is canonically isomorphic to $\pi^*\cG(-1)$.
\end{Lemma}

\begin{proof}
	Since $\pi$ is a projective bundle, we have
	$\pi_{*}(\Omega _{X/Y}^1(1))=0$ and $R^1\pi_{*}(\Omega _{Z/T}^1(1))=0$. 
	We have a short exact sequence
	$$0\to \Omega _{X/Y}^1(1) \to (\Omega _{X/Y}^1 (\log D))(1)\to \cO_{D}(1)\to 0,$$
	where the first non-zero map is the canonical inclusion and the second one is the relative Poincar\'e residue. 
	Pushing the above sequence down to $Y$,  we get a canonical isomorphism 
	$$ \pi_{*}((\Omega _{X/Y}^1 (\log {D}))(1))\mathop{\longrightarrow}^{\simeq}\pi_{*} (\cO_{D}(1))=\cG.$$ 
	Let us recall that if $H$ is a hyperplane in $\PP^{m}$ then $(\Omega_{\PP^{m}}^1(\log H))(1)\simeq \cO_{\PP^{m}}^{\oplus m}$. So the relative evaluation map
	$$ \pi^*\pi_{*}((\Omega _{X/Y}^1 (\log {D}))(1))\longrightarrow (\Omega _{X/Y}^1 (\log {D}))(1)$$
	is an isomorphism on the fibers of $\pi$ and hence it is an isomorphism. This shows that 
	the relative logarithmic cotangent bundle $\Omega _{X/Y}^1 (\log {D})$ is canonically isomorphic to $\pi^*\cG(-1)$.
\end{proof}

\medskip

In the following lemma we use the notation from Subsection \ref{subsection-Atiyah-diff}.

\begin{Lemma}\label{exact-sequence-without-contact}
	We have a canonical short exact sequence
	$$0\to \pi^* \cA_{{\cM}} \to \cB_{\cL}\to  \pi^*\cG (-1)\to 0.$$ 
\end{Lemma}

\begin{proof}
	In the short exact sequence
	$$0\to \Omega_{X}^1\to \Omega_X^1(\log {D})\to \cO_{D}\to 0,$$
	the connecting map $H^0(X,\cO_{D})\to H^1(X,\Omega_X^1)$ maps $1$ to $c_1(\cO_X({D}))$. So the image of $c_1(\cO_X({D}))=c_1(\cL)- c_1(\pi^*\cM)$ in $H^1(X,\Omega_X^1(\log {D}))$  vanishes. This shows that $c_1^{\log}(\cL)= c_1^{\log}(\pi^*\cM)$ and hence
  the logarithmic Atiyah extensions of $\cL$ and $\pi^*\cM$ for the pair $(X, D)$ are isomorphic. 
	By naturality of the first Chern class $\pi^*c_1(\cM) \in H^1(X,\pi^*\Omega_Y^1)$ is mapped to $c_1(\pi^*\cM) \in H^1(X,\Omega_X^1)$. Thus  we get the following commutative diagram
	$$\xymatrix{
		0\ar[r]&\pi^*\Omega_Y^1\ar[r]\ar[d]& \pi^*\cA  _{\cM}\ar[d]\ar[r]& \pi^* \cO_Y\ar[r]\ar@{=}[d]&0\\
		0\ar[r]&\Omega_X^1(\log D)\ar[r]& \cB _{\cL}\ar[r]& \cO_X\ar[r]&0\\
	}$$
	Hence by Lemma \ref{relative-cotangent-bundle} we have the short exact sequence
	$$0\to  \pi^* \cA_{{\cM}} \to \cB_{\cL}\to  \Omega_{X/Y}^1(\log {D}) \simeq  \pi^*\cG(-1)\to 0.$$
\end{proof}

\begin{Remark}
Let us point out that the image of 
$\pi^*\cA  _{\cM}\to \cB _{\cL}\to \cB _{\cL/Y}$ is isomorphic to $\cO_X$, so it gives a canonical splitting of the relative logarithmic Atiyah extension
$$0\to \Omega_{X/Y}^1(\log {D}) \to \cB_{\cL/Y}\to \cO_X\to 0.$$
Since $\Omega_{X/Y}^1(\log {D}) \simeq  \pi^*\cG(-1)$, this splitting can be also given by the relative evaluation map
$\cO_X\simeq \pi^*\pi_* \cB_{\cL/Y}\to  \cB_{\cL/Y}.$
This follows from the fact that $\pi^*\cA  _{\cM}\to \cB _{\cL/Y}$ factors through  $\pi^*\cA  _{\cM}\to \pi^*\pi_* \cB_{\cL/Y}$.
\end{Remark}

\section{Contact structures on special Lie algebroids}

In this section we define a semi-contact structure on an snc pair and classify such structures in case $b_2=1$. Then we study the Kodaira dimension of semi-contact snc pairs and prove some result on existence of lines.

\subsection{Contact snc pairs and special Lie algebroids}

Let $(X,D)$ be an snc pair.

\begin{Definition} 
	Let $\cL$ be a line bundle on $X$.
	A \emph{(generically) $\cL$-valued contact structure}  on $(X,D)$ is a (generically) $\cL$-valued contact structure on the canonical Lie algebroid $T_X(-\log D)$. We say that $(X,D, \cL )$ is 
	\emph{(generically)  contact} if $(X,D)$ is an snc pair, and
there exists a (generically) $\cL$-valued contact structure on $(X, D)$.
Sometimes we say that $(X,D )$ is \emph{(generically)  contact} if there exists a line bundle $\cL$ on $X$ such that  $(X,D, \cL )$ is  {(generically)  contact}.
\end{Definition}

Note that if $\theta$ is a generically $\cL$-valued contact structure on $(X,D)$
then there exists a  non-empty Zariski open subset $V\subset X$ such that $\theta|_V$ is an $\cL|_V$-valued contact structure on $V$ (it is sufficient to take $V=U\backslash \Supp D$ for $U$ in the definition). So $(X,D)$ is generically contact if and only if $X$ is generically contact. Clearly, a smooth projective variety which is birational to a generically contact variety is also generically contact.
  
Note that \cite[Section 1]{Buczynski-Kapustka-Kapustka} uses a different definition of a generically contact variety in which the map $T_X\to \cL$ is surjective. This notion strongly restricts possible varieties and it depends on the choice of a birational model.
 
\medskip
 
Let $(X, D)$ be a projective snc pair. We have on $X$ a stratification induced by the components of the divisor $D$.
 To describe a stratum, we write $D=\sum _{i\in I} D_i$ as a sum of irreducible divisors and we fix $J\subset I$.
 Then a stratum $Z$ is an irreducible component of $D_J$. In the following we write $D^Z$ for the divisor $D^J|_Z$. To make things more precise we write $\imath $ for the closed embedding $Z\hookrightarrow X$. 
 An important  role in this paper is played by the following Lie algebroid structure induced on $Z$.

 \begin{Lemma}\label{special-Lie-algebroid-structure}
 	Let $m= \dim X-\dim Z$. Then we have a canonical Lie algebroid structure on $ \imath ^*T_{X}(-\log D)$ and a canonical central extension 	
 	$$0\to\cO_{Z}^{\oplus m}\to \imath ^* T_{X}(-\log D)\to  T_{Z}(-\log D^Z)\to 0
 	$$
 	of Lie algebroids. 
 \end{Lemma}
 
 \begin{proof}
 	Exactness of the above sequence as a sequence of $\cO_Z$-modules follows from standard short exact sequences for logarithmic bundles by taking duals. For example, the first non-zero map in the sequence is dual to the direct sum of restrictions to $Z$  of Poincar\'e residues $\Omega ^1_X(\log D)\to \cO_{D_i}$ for $i\in J$.  The anchor map for $\imath ^* T_{X}(-\log D)$ is given by the composition $\imath ^* T_{X}(-\log D)\to  T_{Z}(-\log D^Z)\hookrightarrow T_Z$. The Lie bracket  for 
 	$\imath ^* T_{X}(-\log D)=\cO_Z \otimes _{\imath ^{-1}\cO_X}\imath ^{-1} T_{X}(-\log D)$
 	is given by the formula
 	$$[f_1 \otimes \delta_1 , f_2 \otimes \delta_2]= 1\otimes [\tilde f_1\delta_1, \tilde f_2\delta_2],$$
 	where $f_1,f_2\in \cO_Z$, $\delta_1, \delta_2\in T_{X}(-\log D)$ and $\tilde f_i\in \cO_X$ is a lift of $f_i$ for $i=1,2$. It is easy to check that the formula does not depend on the choice of lifts (here it is convenient to interpret the Lie bracket as the standard bracket in the ring of differential operators and use the fact that if $f\in I_Z$ then $\delta (f)\in I_Z$ for $\delta\in T_{X}(-\log D)$). This gives rise to the required Lie algebroid structure.
 	It is also clear from the formula, that the map $\imath ^* T_{X}(-\log D)\to  T_{Z}(-\log D^Z)$ is a map of Lie algebroids.
 	So the anchor map on the kernel is trivial and hence the corresponding Lie algebra structure is $\cO_Z$-linear.
 	Now let us choose local (algebraic) coordinates $(x_1,...,x_n)$ so that $Z$ is given by $x_1=...=x_m=0$. Then $T_X(-\log D)$
 	has a local basis $x_1\frac{\partial}{\partial x_1}, ...,x_m\frac{\partial}{\partial x_m}, \frac{\partial}{\partial x_{m+1}}, ..., \frac{\partial}{\partial x_n}$ and the kernel of $\imath ^* T_{X}(-\log D)\to  T_{Z}(-\log D^Z)$ 
 	is generated as an $\cO_Z$-module by $1\otimes \left( x_i\frac{\partial}{\partial x_i}\right)$ for $i=1...,m$. Since $\left[ \frac{\partial}{\partial x_i},  \frac{\partial}{\partial x_j}\right]=0$, we have
 	$$1\otimes \left[x_i\frac{\partial}{\partial x_i}, \sum _{j\le m}f_j x_j\frac{\partial}{\partial x_j}+\sum _{j> m} f_j\frac{\partial}{\partial x_j}\right]=  \sum _{j\le m}x_i\frac{\partial f_j}{\partial x_i}\otimes x_j\frac{\partial}{\partial x_j}+\sum _{j> m} x_i\frac{\partial f_j}{\partial x_i} \otimes \frac{\partial}{\partial x_j}= 0$$ for all $1\le i\le m$.
So  the Lie bracket on the kernel is trivial and  $\left\{1\otimes \left( x_i\frac{\partial}{\partial x_i}\right)\right\}_{i=1,..,m}$
define central sections of  $\imath ^* T_{X}(-\log D)$.  
 \end{proof}

\begin{Remark}
	The above lemma allows us to upgrade $\imath: Z\hookrightarrow X$ to a morphism of D-schemes
	$$(Z, \cD_{T_{Z}(-\log D^Z)})\to (X, \cD_{T_{X}(-\log D)})$$
	in the sense of \cite[1.4.2]{Beilinson-Bernstein-Jantzen}.
\end{Remark}

\begin{Remark}
One can easily check that the class of the extension  from Lemma \ref{special-Lie-algebroid-structure}
is given by $$\bigoplus _{i\in J}c_1^{\log}(\cO_Z(D_i))\in \Ext ^1( T_{Z}(-\log D^Z), \cO_{Z}^{\oplus m})=H^1(Z, \Omega_Z^1(\log D^Z))^{\oplus m}.$$ 
\end{Remark}

\medskip
 
 Note that if $(X,D, \cL)$ is a  contact projective log variety then for any $Z$ as in the lemma, the restriction 
 $\imath ^* T_{X}(-\log D)\to \imath ^*\cL$ is a  contact structure on the above Lie algebroid $ \imath ^*T_{X}(-\log D)$. This allows us to study logarithmic contact structures by induction on the dimension of the underlying variety. It also motivates the following definition of  contact structures on some  Lie algebroids generalizing Picard Lie algebroids (see \cite[Definition 2.1.3]{Beilinson-Bernstein-Jantzen}).

\begin{Definition}\label{semi-contact-definition}
	Let $(X,D)$ be a projective snc pair. 
	\begin{enumerate}
		\item  A \emph{special Lie algebroid} on $(X,D)$ is a Lie algebroid $\LL$ together with an abelian extension 
		$$0\to \LL' \to \LL \to T_X(-\log D)\to 0$$
		of Lie algebroids. In the following we often omit the extension in the notation.
		\item  A {special Lie algebroid} $\LL$ on $(X,D)$ is called \emph{super-special} if the corresponding extension is central.
		\item 	A \emph{special (generically) contact Lie algebroid} on $(X,D)$ is a triple $(\LL, \cL, \LL\to \cL)$ consisting of a special  Lie algebroid $\LL$, a line bundle $\cL$ and a (generically) $\cL$-valued contact structure $\LL\to \cL$. In the following we omit the contact structure in the notation.
		\item 	We say that $(X, D)$ is  \emph{(generically) semi-contact of rank $(2r+1)$}  if there exists a  special (generically) contact Lie algebroid  $(\LL, \cL)$  on $(X, D)$ with $\rk \LL=2r+1$.
	\end{enumerate}
\end{Definition}

The proof of the following lemma is similar to that of Lemma \ref{special-Lie-algebroid-structure} and we leave it to the reader.

\begin{Lemma}\label{special-Lie-algebroid-structure-v2}
Let $(\LL, \cL)$ be a special contact Lie algebroid on 	 $(X,D)$. Then for any  irreducible stratum $Z$ of the stratification induced by $D$, the restriction  $(\LL|_Z, \cL|_Z)$ is a special contact Lie algebroid on $(Z, D^Z)$.
\end{Lemma}

\subsection{Bounds on existence of a semi-contact structure}

The following result shows that varieties admitting a semi-contact structure are quite rare.
 
\begin{Proposition}\label{dimension-contact-Lie-algebroid}
	Let $(\LL, \cL)$ be a special contact Lie algebroid of rank $(2r+1)$ on $(X,D)$. If $K_X+D$ is not nef then $\dim X\ge r$ and if equality holds then $D=0$. If moreover $b_2(X)=1$ then one of the following cases occurs:
	\begin{enumerate}
		\item $\dim X=2r+1$ and $(X,D)$ is a contact log Fano variety.
		\item 	For some $0\le s\le r$ we have $(X, \cL)\simeq (\PP^{r+s} , \cO_{\PP^{r+s}}(1)) $, $D$ is an arrangement of $s$ hyperplanes in a general position and $\cD_{\LL}^{\le 1}(\cL )\simeq \cO_{\PP^{r+s}}^{\oplus (r+1)}\oplus \cO_{\PP^{r+s}}(1)^{\oplus (r+1)}$. Moreover, we have
		either $\LL\simeq \cO_{\PP^{r+s}}^{\oplus r}\oplus \cO_{\PP^{r+s}}(1)^{\oplus (r+1)}$ or $s=0$ and  $\LL\simeq \cO_{\PP^r}^{\oplus (r+1)}\oplus T_{\PP^r}$. 
	\end{enumerate}
\end{Proposition}

\begin{proof} 
Let us write $\dim X=r+s$ for some $s\le r+1$. By assumption there exists a short exact sequence of Lie algebroids
\begin{equation}\label{seq-1}
0\to\cO_{X}^{\oplus (r+1-s)}\to \LL\to  T_{X}(-\log D)\to 0 	
\end{equation}	
	So existence of an $\cL$-valued contact structure on $\LL$ implies that $-(K_X+D)=(r+1)c_1(\cL)$. If $K_X+D$ is not nef then there exists a good log contraction $\varphi$ of $(X, D)$. So by Lemma \ref{bound-on-length} we can find a rational curve $C$ such that $$0<-(K_X+D)\cdot C=(r+1)\cL\cdot C\le \dim X+1=r+s+1.$$ This shows that
	$s\ge 0$. If $s=0$ then the same lemma implies that $D=0$. 

Now assume that $b_2(X)=1$. Then $-(K_X+D)$ is ample and hence $\cL$ is ample.
For small $\epsilon>0$ the divisor  $-(K_X+D)+\epsilon D$ is also ample.
So by the Kawamata--Viehweg vanishing theorem we have
\begin{equation}\label{van-1}
	H^j(X, \cO_X)= H^j(X, K_X+\lceil -(K_X+D)+\epsilon D \rceil )=0
\end{equation}
for all $j>0$.
For any irreducible component $D_i$ of $D$, $-(K_{D_i}+D^i)=-(K_X+D)|_{D_i}$ is ample, so we also have 
\begin{equation}\label{van-1.5}
H^j(D_i, \cO_{D_i})=0
\end{equation}
for $j>0$ (note that the above argument for vanishing did not use $b_2=1$).
Since $b_2(X)=1$, we have $h^1(X, \Omega _X^1)=1$. Let us consider the Atiyah extension of $\cL$
$$0\to \Omega_X^1 \longrightarrow \cA _{\cL}\longrightarrow \cO_{X} \to 0.$$ 	
Since $c_1(\cL)\in H^1(X, \Omega_X^1)$ is non-zero, the sequence is non-split. Hence the connecting homomorphism $H^0(X, \cO_X)\to H^1(X, \Omega_X^1)$ is non-zero.  Since $H^1(X, \cO_X)=0$ (see (\ref{van-1})), this gives 
\begin{equation}\label{van-2}
	H^1(X,\cA _{\cL} )=0. 
\end{equation}
Let $$0\to \Omega_X^1 (\log D)\longrightarrow \cB _{\cL}\longrightarrow \cO_{X} \to 0	$$
be the logarithmic Atiyah extension of $\cL$. The short exact sequence 
\begin{equation}\label{seq-1.5}
0\to \Omega_X^1 \longrightarrow \Omega_X^1(\log D) \longrightarrow \bigoplus \cO_{D_i} \to 0	
\end{equation}
induces the  short exact sequence
$$0\to \cA_{\cL} \longrightarrow \cB_{\cL} \longrightarrow \bigoplus \cO_{D_i} \to 0.$$
The long exact cohomology sequence for this sequence and (\ref{van-2}) show that 
\begin{equation}\label{van-3}
	H^1(X,\cB _{\cL} )=0.
\end{equation}
In the following we will also use  the short exact sequence
\begin{equation}\label{seq-2}
	0\to \cO_X \longrightarrow \cD_{\LL}^{\le 1}(\cL ) \longrightarrow \LL \to 0.
\end{equation}
Now we divide the proof into two cases depending on the extension (\ref{seq-1}).

\medskip

\emph{Case 1.} If the extension (\ref{seq-1}) does not split then $s\le r$ and there exists  a commutative diagram 
$$\xymatrix{
	0\ar[r]&\cO_{X}^{\oplus (r+1-s)}\ar[r]\ar@{>>}[d]& \LL\ar[r]\ar[d]& T_{X}(-\log D)\ar[r]\ar@{=}[d]& 0\\
	0\ar[r]&\cO_{X}\ar[r]& \cB_{\cL}^*\ar[r]& T_{X}(-\log D)\ar[r]& 0\\
}$$
The middle column of this diagram induces a short exact sequence
\begin{equation*}
	0\to\cO_{X}^{\oplus (r-s)} \to  \LL\to \cB_{\cL}^*\to 0.
\end{equation*}
Since the corresponding Ext group vanishes by (\ref{van-3}), this sequence splits and $\LL \simeq \cO_{X}^{\oplus (r-s)}\oplus \cB_{\cL}^*$.
Since $\Ext^1(\LL, \cO_X)=H^1(X,\cO_{X}^{\oplus (r-s)}\oplus \cB_{\cL})=0$ by (\ref{van-1}) and (\ref{van-3}),
sequence (\ref{seq-2}) splits.  Therefore  
$\cD_{\LL}^{\le 1}(\cL )  \simeq \cO_{X}^{\oplus (r-s+1)}\oplus \cB_{\cL}^*.$
By Proposition \ref{general-LeBrun-lemma} there exists an isomorphism $\cD_{\LL}^{\le 1}(\cL )\simeq \cL\otimes (\cD_{\LL}^{\le 1}(\cL )) ^*$.  
So 
$$\cO_{X}^{\oplus (r-s+1)}\oplus \cB_{\cL}^*\simeq \cL^{\oplus (r-s+1)}\oplus (\cL\otimes \cB_{\cL}).$$
By the Krull--Schmidt theorem, a decomposition of a vector bundle into a direct sum of indecomposable bundles is unique up to reordering of the factors (see \cite{Atiyah-Krull-Schmidt}). In particular, $\cL^{\oplus (r-s+1)}$ is a direct factor of 
$\cB_{\cL}^*$. But we have a short exact sequence 
$$0\to \cO_{X} \longrightarrow \cB ^*_{\cL}\longrightarrow T_{X}(-\log D) \to 0.$$
So if $s<r+1$ then we get a non-zero map $\cL\to T_{X}(-\log D)$. Since $\cL$ is ample and $-(K_X+D)=(r+1)c_1(\cL)$, Theorem \ref{log-Wahl} shows that $X\simeq \PP^{r+s}$, $\cL \simeq \cO_{\PP^{r+s}}(1)$ and  $D$ is an arrangement of $s$ general hyperplanes.
If $s=0$ then $\cB ^*_{\cL}\simeq \cO_{\PP^{r}}(1)^{\oplus (r+1)}$.
If $s>0$ then $T_{X}(-\log D)\simeq \cO_{\PP^{r+s}}^{\oplus (s-1)}\oplus \cO_{\PP^{r+s}}(1)^{\oplus (r+1)}$
by \cite[Proposition 2.10]{Dolgachev-Kapranov}. In both cases we have
$$\cD_{\LL}^{\le 1}(\cL )\simeq \cO_{\PP^{r+s}}^{\oplus (r+1)}\oplus \cO_{\PP^{r+s}}(1)^{\oplus (r+1)}.$$
 	
\medskip

\emph{Case 2.} If the extension (\ref{seq-1}) splits then $\LL\simeq \cO_{X}^{\oplus (r-s+1)}\oplus  T_{X}(-\log D)$. If the sequence (\ref{seq-2}) splits  then 
$\cD_{\LL}^{\le 1}(\cL )  \simeq \cO_{X}^{\oplus (r-s+2)}\oplus T_{X}(-\log D).$
So we have an isomorphism 
$$\cO_{X}^{\oplus (r-s+2)}\oplus  T_{X}(-\log D)\simeq \cL^{\oplus (r-s+2)}\oplus (\cL\otimes \Omega_X^1(\log D)).$$
Then $\cL^{\oplus (r-s+2)}$ is a direct factor of $T_{X}(-\log D)$. As above Theorem \ref{log-Wahl}  implies that  $X\simeq \PP^{r+s}$, $\cL \simeq \cO_{\PP^{r+s}}(1)$ and  $D$ is an arrangement of $s$ general hyperplanes. But in this case $T_{X}(-\log D)$ contains  only $(r-s+1)$ copies of $\cL$, a contradiction.
It follows that  the sequence  (\ref{seq-2}) does not split.
However, by (\ref{van-1}) we have 
$$\text{Ext}^1(\LL, \cO_X)=H^1(X, \LL^*)\simeq H^1(X, \cO_X)^{\oplus(r-s+1)} \oplus H^1(X, \Omega_X^1(\log D))=H^1(X, \Omega_X^1(\log D)).$$
By (\ref{van-1.5}) the canonical map  $H^1(X, \Omega _X^1)\to H^1(X, \Omega _X^1(\log D) )$, obtained from sequence (\ref{seq-1.5}), is surjective and the only non-zero class in $ H^1(X, \Omega _X^1(\log D) )$ (up to scalar) is equal to the image of $c_1(\cL)$. So $\cD_{\LL}^{\le 1}(\cL )  \simeq \cO_{X}^{\oplus (r-s+1)}\oplus \cB_{\cL}^*.$ As above this implies that if $s<r+1$ then $X\simeq \PP^{r+s}$,  $D$ is an arrangement of $s$ general hyperplanes and $\cD_{\LL}^{\le 1}(\cL )\simeq \cO_{\PP^{r+s}}^{\oplus (r+1)}\oplus \cO_{\PP^{r+s}}(1)^{\oplus (r+1)}.$	
\end{proof}

\medskip

Existence of semi-contact structure in case (2) is described in the following example:

\begin{Example}\label{example-P-{r+s}}
Let us first consider $X=\PP^{2r+1}$ with an arrangement $D$ of $(r+1)$ hyperplanes $H_1,...,H_{r+1}$ in a general position.
Then $(X,D)$ has a contact structure given by projection 
$ T_{X}(-\log D)\simeq \cO_X^{\oplus r} \oplus\cO_X(1)^{\oplus (r+1)}\to \cO_X(1)$ onto one of the factors.
For every $0\le s\le r$, restricting to $H_1\cap ...\cap H_{r+1-s}\simeq \PP ^{r+s}$ and using Lemma \ref{special-Lie-algebroid-structure}, we get a semi-contact structure on 
$\PP^{r+s}$ with an arrangement of $s$ hyperplanes in a general position. 
\end{Example}

\medskip

Note that in the proof of Proposition \ref{dimension-contact-Lie-algebroid} we did not use the Lie algebra structure on $\LL$. But we use it to show the following slightly stronger form of the first assertion.

\begin{Lemma}\label{rank-estimate}
If $(X, D)$ is a semi-contact snc pair of rank $(2r+1)$ then $r\le \dim X$.
\end{Lemma}

\begin{proof}
Let $(\LL, \cL)$ be a special contact Lie algebroid of rank $(2r+1)$ on $(X,D)$
with $\cF=\ker (\theta: \LL \to \cL)$. By definition we have an abelian extension 
$$0\to \LL' \to \LL \to T_X(-\log D)\to 0$$
of Lie algebroids. 

We claim that the map $\gamma: \LL'\cap \cF \to \cHom _{\cO_X}(\LL, \cL)$ mapping $x\in \LL'\cap \cF$ to
the map $y\to \theta([x,y]_{\LL})$ is $\cO_X$-linear.
Let us first note that for any  $x\in \LL'$ the map $\LL\to \LL$ sending $y$ to $[x,y]_{\LL}$ is $\cO_X$-linear. Indeed, for any $x\in \LL'$ and $y\in \LL$ 
the Lie algebroid Leibniz rule gives
$$[x, fy]_{\LL} = \alpha _x(f) y+f[x, y]_{\LL} =f[x, y]_{\LL}.$$
This implies that the map  $y\to \theta([x,y]_{\LL})$ is also $\cO_X$-linear and hence $\gamma$ is well-defined.

Now we need to check that for any $f \in \cO_X$ and $x \in \LL' \cap \cF$
we have $\gamma (fx)=f \gamma (x)$. 
For any $y\in \LL$ we have
$$\gamma (fx)(y)=\theta([fx,y]_{\LL}).$$ 
Using the skew-symmetry of the Lie bracket and the standard Leibniz rule with respect to the anchor map $\alpha: \LL \to T_X(-\log D)$, we can expand the bracket:
$$[fx, y]_{\LL} = -[y, fx]_{\LL} = -\left(f[y, x]_{\LL} + \alpha_y(f)x\right) = f[x, y]_{\LL} - \alpha_y(f)x.$$
Now, apply the $\cO_X$-linear contact form $\theta$:
$$\theta([fx, y]_{\LL}) = f\theta([x, y]_{\LL}) - \alpha_y(f)\theta(x).$$
Since  $x$ lies in $\LL' \cap \cF$, we  have $\theta(x) = 0$, which with the above implies $\gamma (fx)(y)=f \gamma (x)(y)$.

Now let us recall that by definition of a contact structure, the map $\hat\theta: \cF \to \cHom_{\cO_X} (\cF, \cL)$, $x\to \theta ([x, \cdot]_L)$ is an isomorphism.
This implies that the map $\hat\theta|_{\LL'\cap \cF}$ is injective. But this map is the composition of $\gamma$ with the canonical projection 
$ \cHom _{\cO_X}(\LL, \cL)\to \cHom_{\cO_X} (\cF, \cL)$, so $\gamma$ is also injective.
Since $[x, \LL']_{\LL}=0$ and $\LL/\LL'\simeq T_X(-\log D)$, this implies that
the composition of $\gamma$ with the canonical projection 
$ \cHom _{\cO_X}(\LL, \cL)\to \cHom_{\cO_X} ( T_X(-\log D), \cL)$ is also injective.
Therefore ${\LL'\cap \cF}$ has rank $\le \dim X$.
But ${\LL'\cap \cF}$  is the kernel of the map $\alpha\oplus \theta: \LL\to T_X(-\log D)\oplus \cL$, so it has rank $\ge (2r+1)-(\dim X+1)=2r-\dim X$. This shows that $r\le \dim X$. 
\end{proof}

\subsection{Kodaira dimension of semi-contact snc pairs}

In the following we write  $\kappa (\cM)$ for the Iitaka dimension of $\cM$ (see \cite[Definition 2.1.3]{Laz04}).
The proof of the following lemma follows that of \cite[Proposition 2]{Druel-Contact}.

\begin{Lemma}	\label{Druel's-lemma}
	 If  $(X,D, \cL )$ is a generically contact projective snc pair then $\kappa (\cL ^{-1})=-\infty$.
	In particular, if $(X,D)$ is contact then $\kappa (K_X+D)=-\infty$.
\end{Lemma}

\begin{proof}
	Assume that $\kappa (\cL ^{-1})\ge 0$. By \cite[Proposition 4.1.6]{Laz04} there exists a smooth projective  variety $\bar X$ and a generically finite morphism $\pi: \bar X\to X$ such that $h^0(\bar X, \bar \cL^{-1})>0$, where $\bar \cL=\pi^*\cL$. Passing to a log resolution, we can also assume that $\bar D:= (\pi^*D)_{\sf red}$ is a simple normal crossing divisor.
	By assumption there exists a non-empty Zariski open subset $U\subset X$ and  $\omega\in H^0(X,\Omega _X^1(\log D)\otimes _{\cO_X} \cL)$ such that $\omega|_U$ is an $\cL|_U$-valued contact form on $U$. There also exists a non-empty Zariski open subset $V\subset U\backslash \Supp D$ such that $\bar V=\pi^{-1}(V)\to V$ is finite \'etale. Taking any non-zero section $s\in H^0(\bar X, \bar \cL^{-1})$ we can find a non-empty Zariski open subset $W\subset \bar V$ such that $s$ has no zeroes on $W$. Note that by the logarithmic ramification formula, we have a natural injective $\cO_{\bar X}$-linear map $\pi^*\Omega_X^1(\log D)\hookrightarrow \Omega^1_{\bar X} (\log \bar D)$.
	Composing it with $\pi^*\omega$ we obtain
	$\bar \omega : \bar \cL^{-1}\to \Omega^1_{\bar X} (\log \bar D)$. Let us set $\eta:=\bar\omega(s)\in	H^0(\bar X,\Omega_{\bar X}^1(\log\bar D)).$ 
	Since $\bar \omega|_W=\pi^* (\omega |_V)|_W$ gives an $\bar \cL|_{W}$-valued contact form, we have $d (\eta|_W)\ne 0$ in
	$H^0(W, \Omega^2_{\bar X}|_W)=H^0(W, \Omega^2_{\bar X}(\log \bar D)|_W)$. But by degeneracy of the logarithmic Hodge to de Rham spectral sequence on   $(\bar X, \bar D)$,	all global logarithmic $1$-forms on $(\bar X, \bar D)$ are closed. Hence  $d\eta =0\in H^0(\bar X, \Omega^2_{\bar X} (\log \bar D))$. So  $d (\eta|_W)=d\eta |_W$ also vanishes, a contradiction.  
\end{proof}

Let us recall that a line bundle $\cM$ on a smooth projective variety $X$ is \emph{pseudoeffective} if and only if $c_1(\cM)\in N^1(X)_{\RR}$ lies in the closure of the cone generated by classes of effective divisors.
The following  theorem was proven by J.-P. Demailly in \cite[Corollary 1]{Demailly-Contact} in the non-logarithmic case  and then generalized by F. Touzet in \cite[Theorem 5 and Corollary 7.6]{Touzet2016} to the logarithmic case.

\begin{Theorem}
	Let $(X,D)$ be a projective snc pair and let	$\cL$ be a line bundle on $X$. If $\cL^{-1}$ is pseudo-effective and $\omega\in H^0(X, \Omega_X^1(\log D)\otimes \cL)$ is a logarithmic $1$-form with values in $\cL$ then $\omega\wedge d\omega=0$. 
\end{Theorem}

This theorem implies a much stronger form of Lemma \ref{Druel's-lemma}:

\begin{Corollary}\label{Touzet}
	Let $(X, D)$ be  a  projective snc pair.
	If $(X,D, \cL )$ is generically contact then  $\cL^{-1}$ is not pseudo-effective. In particular,
	$\kappa (\cL ^{-1})=-\infty$. If $(X,D)$ is a contact snc pair then $K_X+D$ is not pseudoeffective. In particular, $X$ is uniruled. 
\end{Corollary}

\medskip

The following lemma deals with semi-contact, non-contact snc pairs.

\begin{Lemma}\label{special-contact-Kodaira}
	Let $(\LL, \cL)$ be a special generically contact Lie algebroid of rank $(2r+1)$ on a projective snc pair $(X,D)$. If $\dim X<2r+1$ then $H^0(X, \cL)\ne 0$. 
\end{Lemma}

\begin{proof}
By assumption, we have a short exact sequence
$$	0\to\cO_{X}^{\oplus s}\to \LL\to  T_{X}(-\log D)\to 0 ,$$
in which $s\ge 1$. If $H^0(X, \cL)= 0$ then the composition $\cO_{X}^{\oplus s}\to \LL\to \cL$ vanishes. So  $ \LL\to \cL$ induces a generically surjective map $T_X(-\log D)\to \cL$. Passing to the dual  we see that the contact form
$\omega: \cL^{-1}\to \Omega_{\LL}^1$ factors through $\cL^{-1}\to \Omega_X^1(\log D)$. But then
$\omega\wedge (d_{\LL}\omega)^{\wedge r}: \cL ^{-(r+1)}\to \det \Omega_{\LL}^1$ 
factors through $ \cL ^{-(r+1)}\to  \Omega_X^{2r+1}(\log D)=0$ and we get a contradicton.
\end{proof}

\begin{Corollary}\label{semi-contact-Kodaira}
Let $(X, D)$ be  a  projective snc pair.
If $(X,D, \LL, \cL )$ is semi-contact and  $K_X+D$ is pseudoeffective then $\dim X<\rk \LL$ and
$\cL \simeq \cO_{X}$ (in particular,  $K_X+D\sim 0$).
\end{Corollary}

\begin{proof}
By Corollary \ref{Touzet}, $\dim X<\rk \LL$. So Lemma \ref{special-contact-Kodaira} implies that
$H^0(X, \cL)\ne 0$. Since the cone of pseudo-effective divisors is convex and $\cL^{-1}$ is pseudo-effective, the class of $\cL $ in $N^1(X)$ vanishes. Since $\cL$ has a non-zero section this gives $\cL \simeq \cO_{X}$.
\end{proof}

\medskip

The following example shows that semi-contact snc pairs with $K_X+D\sim 0$ really occur.

\begin{Example}\label{example-semi-abelian-semi-contact}
Let $(X, D)$ be a projective snc pair of dimension $n$ and assume that a semi-abelian variety $G$ acts on $X$ with $X\backslash \Supp D$ as an open orbit.  Since  $G$ is a commutative $k$-group, the Lie algebra $\mathfrak{g}$ of $G$ is abelian. By assumption $\mathfrak {g}$ acts on $X$ and $\cO_X\otimes _k \mathfrak{g}$ has a natural structure of Lie algebroid described in \cite[1.2.2]{Beilinson-Bernstein-Jantzen}. The image of the anchor map $\alpha:\cO_X\otimes _k \mathfrak{g} \to T_X$ coincides with $T_X(-\log D)$ (see, e.g.,  \cite[p.~478]{Winkelmann2004}). The Lie bracket
on  $\cO_X\otimes _k \mathfrak{g}$ is defined by requiring
$$[1\otimes \delta, 1\otimes \delta']= 1\otimes [\delta, \delta']=0$$
for $\delta, \delta'\in \mathfrak{g}$ and extending it to $\cO_X\otimes _k \mathfrak{g}$ by the Leibniz rule.

Let $\sigma$ be a non-degenerate alternating form on a $2r$-dimensional $k$-vector space $V$ and let  $\mathfrak{h}=V\oplus k\cdot z$ be a Lie $k$-algebra of dimension $(2r+1)$ with the Lie bracket determined by
$$[v,w]=\left\{ \begin{array}{cl}
\sigma (v,w)\cdot z &\hbox{if $v, w\in V$,}\\
0 & \hbox{if $v\in V$ and $w=z$.}\\
\end{array}
\right.$$ 
Assume that  $r\le n\le 2r$ and let $W\subset V$ be an isotropic $k$-linear subspace of dimension $(2r-n)$.
Then $I=W\oplus k\cdot z\subset \mathfrak{h}$ is an ideal  and the quotient  $\mathfrak{h}/I$ is an abelian Lie algebra  of dimension $n$, so it is isomorphic to $\mathfrak{g}$. 
Let us fix a surjective map of Lie algebras $\mathfrak{h}\to \mathfrak{g}$.
Then $\LL=  \cO_{X}\otimes _k\mathfrak{h}$ has a canonical structure of a Lie algebroid 
for which the extension
$$0\to\cO_{X}\otimes _kI\to \LL\to  \cO_{X}\otimes _k\mathfrak{g}\simeq  T_X(-\log D)\to 0$$
is abelian (one can again use \cite[1.2.2]{Beilinson-Bernstein-Jantzen}).  A contact structure on $\LL$ is defined by the map $\LL\to \cL=\cO_{X}$ induced	by the projection $\mathfrak{g}\to k\cdot z$. In this way for every $n/2\le r\le n$  we obtain a rank $(2r+1)$ semi-contact structure on $(X, D)$. Note that 
the above extension is central only for $n=2r$.
\end{Example}

The following example shows that Corollary \ref{semi-contact-Kodaira} fails if  $(X,D, \LL, \cL )$ is only generically semi-contact.

\begin{Example}
Let $X=C_1\times C_2\times C_3$, where the \(C_i\) are smooth projective curves of genus at least $2$, and take $D=0$. 
For $i=1,2,3$ choose nonconstant rational functions on $C_i$, and let $x_i$ be the induced rational functions on $X$. 
Let us set
$$\eta=x_1\,dx_2+dx_3.$$
Then $d\eta=dx_1\wedge dx_2$
and hence $$\eta\wedge d\eta
=dx_1\wedge dx_2\wedge dx_3\neq 0.$$
After choosing an effective divisor $E$ containing all poles of $\eta$, $\eta$ defines a global section of $\Omega_X^1\otimes\mathcal O_X(E)$.
It therefore gives a generically $\cO_X(E)$-valued contact structure on $\LL:=T_X$. But
 $K_X$ is ample and $\dim X=\rk \LL$.
\end{Example}

\subsection{Lines on semi-contact snc pairs}

Let us first note the following easy lemma:

\begin{Lemma}\label{rank+onP^1}
	Let 
	$$0\to \cO_{\PP^1}\to \cE_1 \to \cE_2\to 0$$
	be an exact sequence of vector bundles on $\PP^1$. Then one of the following holds:
	\begin{enumerate}
		\item the sequence does not split and $\rk ^{+}\cE_2=\rk ^{+}\cE_1-1$, or
		\item  the sequence splits and $\rk ^{+}\cE_2=\rk ^{+}\cE_1$.
	\end{enumerate}
\end{Lemma}

\begin{proof}
	If the sequence splits then the assertion is clear. If it does not split then  $\rk ^{+}\cE>0$. In this case, 
	after dualizing the above sequence,	the required equality follows from  \cite[Lemma 2.7]{KPSW}.
\end{proof}

By passing to the dual sequence and induction on $m$ we get the following corollary:

\begin{Corollary}\label{rank+onP^1-2}
	Let 
	$$0\to \cE_2 \to \cE_1\to \cO_{\PP^1} ^{\oplus m}\to 0$$
	be an exact sequence of vector bundles on $\PP^1$. Then
	$\rk ^{+}\cE_2^*\le \rk ^{+}\cE _1^*$
and equality holds if and only if the  sequence splits.
\end{Corollary}

\medskip

\begin{Lemma}\label{rank+Lie-algebroids}
	Let $(\LL, \cL)$ be a rank $(2r+1)$ special contact Lie algebroid on a projective snc pair $(X,D)$. Let $f: \PP^1\to C\subset X$ be a normalization of a rational curve with $\deg f^*\cL>0$.
	Then we have $\rk ^{+}f^*T_{X}(-\log D)\le r+1$. Moreover, if $C\cap \Supp D=\emptyset$ then
	$\rk ^{+}f^*T_{X}\le r$.	
\end{Lemma}

\begin{proof} 		
	Let us set $\cB=  \cA_{\LL}(\cL)$ and
	consider the $\LL$-Atiyah extension of $\cL$
	$$0\to \Omega_{\LL}^1\longrightarrow \cB \longrightarrow \cO_X\to 0.$$
	By Proposition \ref{general-LeBrun-lemma}, $\cD_{\LL}^{\le 1}(\cL )=\cB^*$ carries an $\cL$-valued symplectic structure and hence we have an isomorphism $\cB\simeq \cL^{-1}\otimes \cB ^*$. So for some integers $a_1\ge ...\ge a_{r+1}\ge 0$ we have
	$$f^*\cB \simeq\bigoplus _{i=1}^{r+1}(\cO_{\PP^1}(a_i)\oplus \cO_{\PP^1}(-a_i-d)),$$
	where $d=\deg f^*\cL>0$.
	In particular, $\rk ^{+} f^*(\cB^*)= r+1$.
	If the sequence 
	$$0\to f^*\Omega_{\LL}^1\to f^*\cB \to f^*\cO_X=\cO_{\PP^1}\to 0 $$
	does not split then by Lemma \ref{rank+onP^1} we have $\rk ^{+}f^*\LL=r$.
	If the above sequence splits then $\rk ^{+}f^*\LL =r+1$.
	
	We also have a short exact sequence
		$$0\to f^*\Omega_{X}^1(\log D)\to f^*\Omega_{\LL}^1\to \cO_{\PP^1}^{\oplus (2r+1-\dim X)}\to 0$$	
obtained by pulling back the dual to the sequence defining a special contact Lie algebroid.   
So using Corollary \ref{rank+onP^1-2} we get $\rk ^{+}f^*T_{X}(-\log D)\le r+1$ with equality if and only if both the above sequences split. 
	
	Now assume that   $C\cap \Supp D=\emptyset$  so that $f^*\Omega_{X}^1\simeq f^*\Omega_{X}^1(\log D)$. By construction of the Atiyah extensions, we have a commutative diagram
	$$\xymatrix{
		0\ar[r]&f^*\Omega_{X}^1\ar[r]\ar[d]&f^*(\cA_{\cL})\ar[r]\ar[d] &\cO_{\PP^1}\ar@{=}[d]\ar[r]&0\\
		0\ar[r]&f^*\Omega_{\LL}^1\ar[d]\ar[r]&f^*\cB\ar[r]\ar[d] &\cO_{\PP^1}\ar[r]&0\\
		&	\cO_{\PP^1}^{\oplus m}\ar@{=}[r]&	\cO_{\PP^1}^{\oplus m}&&&\\
	}$$
	If  $\rk ^{+}f^*T_{X}=\rk ^{+}f^*T_{X}(-\log D)=r+1$ then we have an induced splitting 
	$$f^*(\cA_{\cL} )\to f^*\cB\to f^*\Omega_{\LL}^1\to f^*\Omega_{X}^1$$
	of the first row in the above diagram. But $\deg f^*\cL\ne 0$, so the commutative diagram
	$$\xymatrix{
		0\ar[r]&f^*\Omega_{X}^1\ar[r]\ar[d]&f^*(\cA_{\cL })\ar[r]\ar[d] &\cO_{\PP^1}\ar@{=}[d]\ar[r]&0\\
		0\ar[r]&\Omega_{\PP^1}^1\ar[r]&\cA _{f^*\cL} \ar[r] &\cO_{\PP^1}\ar[r]&0\\
	}$$
	shows that this sequence does not split. Therefore 	$\rk ^{+}f^*T_{X}\le r$.
\end{proof}

\medskip

\begin{Lemma}\label{existence-of-lines}
Let $(\LL, \cL)$ be a rank $(2r+1)$ special contact Lie algebroid on a projective snc pair $(X,D)$. Assume that $K_X+D$ is not nef and let $\varphi$ be a good log contraction of $(X, D)$. 
Then either $(X, D, \LL, \cL)\simeq (\PP^{2r+1}, 0, T_{\PP^{2r+1}}, \cO_{\PP^{2r+1}} (2))$ or there exists a $\varphi$-exceptional rational curve $C$ with $\cL \cdot C=1$.
\end{Lemma}

\begin{proof} 
Without loss of generality we can assume that $\varphi$ is elementary and it contracts some $(K_X+D)$-extremal ray $R$. If $(X, D)\simeq (\PP^{2r+1}, 0)$ then $( \LL, \cL)\simeq (T_{\PP^{2r+1}}, \cO_{\PP^{2r+1}} (2))$. Otherwise, 
Lemma \ref{bound-on-length} implies that $ l_{D}(\varphi)\le 2r+1$. Since $-(K_X+D)=(r+1)c_1(\cL)$, $l_{D}(\varphi)$ is divisible by $r+1$. Therefore $ l_{D}(\varphi)=r+1$ and if $C$ is an extremal rational curve with $[C]\in R$ then $\cL\cdot C=1$. 
\end{proof}

If  $(X,D)$ is a contact projective  snc pair then  Corollary \ref{Touzet} implies that $K_X+D$ is not nef. So we have the following corollary.

\begin{Corollary}
If  $(X,D, \cL )$ is a contact projective  snc pair and $(X, D, \cL)\not \simeq (\PP^{2r+1}, 0, \cO_{\PP^{2r+1}} (2))$ then $K_X+D$ is not nef and for every good log contraction $\varphi$  there exists an extremal rational curve $C$ with $\cL \cdot C=1$.
\end{Corollary}

\section{Log contractions of semi-contact snc pairs}

This section contains proof of a generalization of Theorem \ref{main1}.
By Corollary \ref{normality-of-image}  good log contractions of semi-contact varieties induce good log contractions on the strata of stratification induced by $D$. This allows us to reduce to studying good log contractions, whose exceptional locus is not contained in the boundary. We call such contractions \emph{very good}.

In this section we use results of Section \ref{section:rational-curves} to describe the general fiber and we prove that good log contractions  of semi-contact varieties are equidimensional. This allows us to conclude that our contractions are projective bundles (we also need to prove this for every stratum of the stratification induced by the boundary divisor). The main result of the section is Theorem \ref{structure-of-log-contractions-special-Lie-alg}. Its proof requires further study of the geometry using various Atiyah extensions.

\medskip

From now on in this section we fix the following notation. Let $(\LL, \cL)$ be a special contact Lie algebroid of rank $(2r+1)$ on a projective snc pair $(Z,D^Z)$. Let us  write  $n=\dim Z$ and $U_Z=Z\backslash \Supp D^Z$.
Assume that  $(K_Z+D^Z)$ is not nef so that $(Z,D^Z)$ admits a good log contraction $\varphi_Z: Z\to S$. Let $F$ be the $(K_Z+D^Z)$-extremal face  corresponding to $\varphi_Z$. Till the end of the section we  assume that the exceptional locus $\Exc (\varphi _Z)$ is not contained in the support of $D^Z$, i.e., that $\varphi_Z$ is very good.  We also assume that $(Z, D^Z, \LL, \cL)\not \simeq (\PP^{2r+1}, 0, T_{\PP^{2r+1}}, \cO_{\PP^{2r+1}} (2))$.

\subsection{General structure of log contractions}

Let $C$ be an extremal rational curve with $[C]\in F$ and let $f: \PP^1 \to C$ be the normalization of $C$. We can choose $C$ so that $C\cap U_Z\ne \emptyset$. By Proposition \ref{Keel-McKernan} we can also assume that $f^*D^Z$ is supported in at most one point so that $f^{-1}(U_Z)=\PP^1$ or $f^{-1}(U_Z)\simeq \AA^1$.

If $H_F$ is a supporting divisor of $F$ then for large $m$ the divisor $mH_F+c_1(\cL)$ is ample and by Lemma \ref{existence-of-lines} the degree of $C$ with respect to this divisor is $1$. In particular, any closed subset of 
the connected component of $\Hom (\PP^1, Z)$ containing $[f]$ is unsplit  (a reducible degeneration would give a curve whose class has $(mH_F+c_1(\cL))$-degree strictly smaller than $1$).

\medskip

Let $V$ be an irreducible component of $\Hom (\PP^1, Z)$, which contains $[f]$.

\begin{Lemma}\label{locus-dimension-P^1}
	If $f^{-1}(U_Z)=\PP^1$ then ${\rm locus} (V)=Z$  and $\dim  {\rm locus} (V,x)=r$ for any $x\in U_Z$. Moreover, we have $\rk ^{+}f^*T_{Z}=r$.
\end{Lemma}

\begin{proof}
	Since the family $V$ is unsplit, Proposition \ref{easy-log-Ionescu-Wisniewski} gives
	\[ \dim  {\rm locus} (V)+ \dim  {\rm locus} (V,x)+1\ge  \dim V \ge -K_{Z}\cdot C+\dim Z\ge -(K_{Z}+D^Z)\cdot C+\dim Z =r+1+n
	\]
	But $ \dim  {\rm locus} (V)\le n$,  so $\dim  {\rm locus} (V,x)\ge r$ for all $x\in Z$.
	Lemma \ref{rank+Lie-algebroids} implies that
	$$r\le \dim  {\rm locus} (V,x)\le \rk ^{+}f^*T_{Z}\le r.$$
	So $\dim  {\rm locus} (V,x)=r$ for all $x\in U_Z$. This implies that $ \dim  {\rm locus} (V)=n$ and hence
	${\rm locus} (V)=Z$.
\end{proof}

Let $W$ be an irreducible component of $\Hom (\PP^1, Z; f^*D^Z\subset D^Z)$ containing $[f]$.

\begin{Lemma}\label{locus-dimension}
	If $f^{-1}(U_Z)\simeq \AA^1$ then ${\rm locus} (W)=Z$ and $\dim  {\rm locus} _{U_Z}(W,x)=r+1$ for any $x\in U_Z$. 
\end{Lemma}

\begin{proof}
By Proposition \ref{log-Ionescu-Wisniewski} we have
\begin{equation*}
	\dim  {\rm locus} _{U_Z}(W)+ \dim  {\rm locus} _{U_Z}(W,x)\ge \dim W\ge -(K_{Z}+D^Z)\cdot C+\dim Z=r+1+n.
\end{equation*}
By  Lemma \ref{rank+Lie-algebroids} for any $g: \PP^1\to X$ such that $[g]\in W$, we have	$\rk ^{+}g^*T_{Z}(-\log D^Z) \le r+1$.	So by Proposition \ref{log-bound-on-rank}   we have
$\dim  {\rm locus} _{U_Z}(W,x) \le r+1.$  Since $ \dim  {\rm locus} _{U_Z}(W)\le n$, the above inequality 
implies that $ \dim  {\rm locus} _{U_Z}(W)=n$ and $\dim  {\rm locus} _{U_Z}(W,x)=r+1$.	
\end{proof}

\begin{Remark}
	Note that if $f^{-1}(U_Z)=\PP^1$ then $f^*D^Z=0$ and $W=V$. So it is always sufficient to consider $W$. We keep the above notation to show differences and similarities to the usual non-logarithmic case (and anyway we need to consider two cases).
\end{Remark}	

\medskip 

The following corollary is an analogue of  \cite[Lemma 2.10]{KPSW}.

\begin{Corollary}\label{fiber-type}
The contraction $\varphi_Z:Z\to S$ is of fiber type, i.e., $\dim S<\dim Z$.	
\end{Corollary}

\begin{proof}
Since $C$ is an extremal rational curve, every element of $V$ gives a curve contracted by $\varphi_Z$.
But	Lemmas \ref{locus-dimension-P^1} and \ref{locus-dimension} imply that ${\rm locus} (V)=Z$, so $\varphi_Z$ cannot be birational. 	
\end{proof}

The following proposition is an analogue of  \cite[Proposition 2.11]{KPSW} in our setting. The strategy of its proof is similar but  we have more cases than just one.

\begin{Proposition}\label{general-fiber-in-general}
	Assume that  with $\dim S>0$.  Then there exists an open subset $U\subset S$ such that for some $0\le s\le r$ and for any closed point $y\in U$ we have  $(Z_y, D_y^Z)\simeq (\PP^{r+s},H_1+...+H_s)$, where $H_1,...,H_s$ are  hyperplanes in a general position. Moreover, we have $\cL |_{Z_y}\simeq \cO_{\PP^{r+s} }(1)$.
\end{Proposition}

\begin{proof} 
By generic smoothness in characteristic zero, there exists an open subset $U\subset S$ such that $\varphi|_{Z_U}: Z_U\to U$ is smooth and 
	$D_U^Z\subset Z_U$ is a divisor with strict normal crossings  relative to $U$. In particular, for every $y\in U$ 	we have a short exact sequence
	$$0\to T_{Z_y}(-\log D^Z_y)\to  T_{Z}(-\log D^Z)|_{Z_y}\to \varphi^*T_{S}|_{Z_y}\simeq \cO_{Z_y}^{\dim S}\to 0.$$
	Let us consider the map $\alpha:\Omega_Z^1(\log D^Z) \otimes \cL\to T_{Z}(-\log D^Z) $ given by composition
	$$\Omega_Z^1(\log D^Z)\otimes \cL\hookrightarrow \Omega_{\LL}^1\otimes \cL \twoheadrightarrow \cF^* \otimes \cL  \mathop{\longrightarrow}^{\simeq} \cF  \hookrightarrow \LL \twoheadrightarrow T_{Z}(-\log D^Z).$$
	This yields the following commutative diagram:
	$$\xymatrix{
		&&T_{Z_y}(-\log D_y^Z)\ar[d]\\
		(\varphi^*\Omega_{S}\otimes \cL)|_{Z_y}\simeq (\cL|_{Z_y})^{\oplus \dim S}\ar[r]^-{\gamma}\ar[rrd]_{0}\ar@{-->}[rru]^{\beta}
		&(\Omega_Z^1(\log D^Z)\otimes \cL )|_{Z_y} \ar[r]^-{\alpha |_{Z_y}}& T_{Z}(-\log D^Z)|_{Z_y}\ar[d]\\
		&&\cO_{Z_y}^{\dim S}\\
	}
	$$
	The lower arrow is zero because $\dim Z_y>0$ by Corollary \ref{fiber-type} and $\cL|_{Z_y}$ is ample as $-(K_Z+D^Z)$ is $\varphi|_Z$-ample.
	 This implies existence of $\beta$.
	The short exact sequence
	$$0\to\cO_{Z}^{\oplus (2r+1-n)}\to \LL \to  T_{Z}(-\log D^Z)\to 0
	$$
	shows that the kernel of $\delta: (\cF^* \otimes \cL )|_{Z_y} \stackrel{\simeq}{\longrightarrow} \cF  |_{Z_y}\hookrightarrow \LL |_{Z_y}\twoheadrightarrow T_{Z}(-\log D^Z)|_{Z_y}$ is contained in $\cO_{Z_y}^{\oplus  (2r+1-n)}$. On the other hand, the short exact sequence
	$$0\to \cO_{Z_y}\to (\Omega_{\LL}^1\otimes \cL )|_{Z_y} \to (\cF^*\otimes \cL)|_{Z_y}\to 0$$
	shows that the kernel of $\alpha|_{Z_y}$ is contained in an extension of $\ker \delta$ by $\cO_{Z_y}$. So $\ker \alpha|_{Z_y}$ is also contained in an extension of $\cO_{{Z_y}}^{\oplus  (2r+1-n)}$ by $\cO_{Z_y}$. 
	Since $\gamma$ is injective this implies that the map $\beta$ is non-zero. So there exists a non-zero map $\cL|_{Z_y}\to T_{Z_y}(-\log D_y^Z)$. By Theorem \ref{log-Wahl} we see that 
	for some $m$ we have $Z_y\simeq \PP^{m}$ and $\cL|_{Z_y}\simeq \cO_{\PP^m }(1)$. Moreover, $D_{y}^Z$
	is a sum of $0\le s\le m$ hyperplanes in a general position. 
	Let $C\subset Z_y$ be a line. Then 
	$$r+1=-(K_X+D)\cdot C=-(K_{Z_y}+D_y^Z)\cdot C_y=m+1-s,$$
	so $m=r+s$. 	 
\end{proof}

\medspace

\begin{Remark}
	If $\LL=T_Z(-\log D^Z)$ and $s=0$, i.e., $(Z_y, D_y^Z)\simeq (\PP^{r}, 0)$, then the composition $$T_{Z_y}=T_{Z_y}(-\log D_y^Z)\to T_Z(-\log D^Z)|_{Z_y}\to \cL |_{Z_y}$$ vanishes and $T_{Z_y}\subset \cF|_{Z_y}$. Indeed, since $\Hom (T_{\PP^r},\cO_{\PP^r}(1) ) =H^0(\PP^r, \Omega^1 _{\PP^r}(1))=0$, there are no non-zero maps $T_{Z_y}\to \cL |_{Z_y}$. This gives a quick proof of the last part of \cite[Proposition 2.11]{KPSW}.
\end{Remark}

\begin{Remark}
Proposition \ref{general-fiber-in-general} shows that 
a posteriori a general fiber  $(Z_y, D_y^Z)$ of $\varphi _Z$ is semi-contact with $b_2(Z_y)=1$ (see Proposition \ref{dimension-contact-Lie-algebroid} and Example \ref{example-P-{r+s}}).
\end{Remark}

\medspace

Let $D^Z=\sum_{i\in I} D_i^Z$ be the decomposition of $D^Z$ into irreducible components. For a subset  $J\subset I$
and a connected  component $Z'$ of $D_J^Z=\bigcap_{i\in J}D_i^Z$ we write $D^{Z'}$ for the restriction of the divisor 
$\sum _{i\in I\backslash J}D_i^Z$ to $Z'$. If $i\in I\backslash J$ then we also write
$D^{Z'}_i$ for the restriction of $D_i^Z$ to $Z'$. We also set  $U_{Z'}:=Z'\backslash \Supp D^{Z'}$.

Assume that  $\dim S>0$. Let $I^h\subset I$ be the subset corresponding to the components dominating $S$ and let
$I^v:=I\backslash I^h$ be its complement. Then we set $D^h:=\sum_{i\in I^h} D_i^Z$ and $D^v:=\sum_{i\in I^v} D_i^Z$.
In the following $D^h$ is called the \emph{horizontal part} of $D^Z$ and $D^v$ is called the \emph{vertical part} of $D^Z$. 

By Proposition \ref{general-fiber-in-general} we have $0\le s=|I^h|\le r$.
Let $U_S$ be the complement of the image of $D^v$ in $S$.
For a subset $J\subset I^h$, we set $U_J'=D^h_J\cap \varphi_Z^{-1} (U_S)$ (note that $D^h_J$ is irreducible). Clearly, $\varphi_Z(D^h_J)=S$ and $\varphi _Z |_{U_J'}$ is a projective morphism with a general fiber isomorphic to $\PP ^{r+s-|J|}$. We need the following result on the dimension of its fibers.

\begin{Proposition}\label{equidimensionality-on-open}
If $\dim S>0$  then	 for any subset $J\subset I^h$ (possibly empty), all fibers of 
$\varphi _Z |_{U_J'}$ have dimension $r+|I^h\backslash J|$.
\end{Proposition}
 
 \begin{proof}	
 	Since $-(K_Z+D^Z)=(r+1)c_1(\cL)$ is $\varphi_Z$-ample, for sufficiently ample line bundle $\cA$ on $T$  the line bundle $\cL'=\cL \otimes \varphi _Z ^*\cA$ is ample. The proof is by induction on the number of elements in $I^h\backslash J$.
 	
 	\medskip
 	
 	\emph{Step 1.} Assume that $|I^h\backslash J|=0$, i.e., $J= I^h$. 
 	
 	Let us consider $Z'=D^h_{I^h}$ (if $I^h=\emptyset$ then we set $Z'=Z$). In this case Proposition \ref{general-fiber-in-general} implies that for general $y\in S$ we have $(Z_y', D_y^{Z'})\simeq (\PP^{r}, 0)$ and $\cL |_{Z'_y}\simeq \cO_{\PP^{r} }(1)$.
 	Let us  consider  the Chow variety $\Chow _{r,1}(Z')$  of $r$-dimensional cycles on $Z'$ of degree $1$ with respect to $\cL'$ and let $Q$ be an irreducible component of $\Chow _{r,1}(Z')$ containing $ Z'_y$ for general $y\in S$. Since $N_{Z'_y/Z'}\simeq \cO_{\PP^r}^{\dim S}$, we have $h^1(N_{Z'_y/Z'})=0$
 	and $\dim Q= h^0(N_{Z'_y/Z'})=\dim S$. So $Q$ is projective of dimension $\dim S$, and there exists a universal family of cycles $\tau: \cU \to Q$ with the evaluation map $\ev : \cU\to Z'$. 
 	 Note that the canonical rational map $\psi: S\dashrightarrow Q$, mapping a general point $y\in S$ to the cycle class
 	$[Z'_y]\in Q$, is dominant and birational. So $\ev$ is also dominant and birational, and either it is an isomorphism or it has a positive dimensional fiber. In the first case, $(\varphi_Z|_{Z'}: Z'\to S)\simeq (\tau: \cU\to Q)$, so  $\varphi_Z|_{Z'}$ is equidimensional. In the second case, there exists $x_0\in Z'$  such that $\dim (\ev ^{-1}(x_0))>0$. Then by the defining property of $\Chow_{r,1}(Z')$, $\dim (\ev( \tau ^{-1}(\tau ( \ev ^{-1}(x_0)))))>r$. 
 	
 	Let $V_{\cU}\subset \Hom (\PP^1, \cU)$ be an irreducible component containing some $f: \PP^1\to \cU$ for which the image is a line in a general fiber  $Z'_y=\tau^{-1}([Z'_y])\simeq \PP^r$ and $f$ is an isomorphism onto its image. For $q\in Q$ the fiber $\tau^{-1}(q)$ represents an $r$-dimensional cycle of degree $1$ with respect to $\cL'$. Since $\ev |_{\tau^{-1}(q)}$ is finite,  $(\ev)^*\cL'$ is ample on every fiber of $\tau$ and hence it is $\tau$-ample. Therefore $V_{\cU}$ parametrizes birational maps to curves of degree $1$ with respect to some ample line bundle and hence the family  $V_{\cU}$ is unsplit. 	 Moreover, the natural map $ \tilde \ev: \Hom (\PP^1, \cU)\to  \Hom (\PP^1,Z')$ maps $V_{\cU}$ into an irreducible component $V$ of $ \Hom (\PP^1, Z')$ containing lines in a general fiber $Z'_y$ of $\varphi_Z|_{Z'}$.
 	Since  $(\varphi_Z|_{Z'}: Z'\to S)$ and $ (\tau: \cU\to Q)$ are isomorphic over dense open subsets where $\psi$ is an isomorphism,  we have ${\rm locus} (V_{\cU}, u)=\tau ^{-1}(\tau (u))\simeq \PP^r $ for general $u\in \cU$.
 	 Since $V_{\cU}$ is unsplit, the equality  ${\rm locus} (V_{\cU}, u)=\tau ^{-1}(\tau (u))$ holds for all $u\in \cU$. This implies that  ${\rm locus} (V, x_0)\supset \ev ( \tau ^{-1}(\tau ( \ev  ^{-1}(x_0))))$ has dimension $\ge (r+1)$.

 Let us assume that $x_0\in U_{Z'}$. If $g: \PP^1\to Z'$ is a map from the family $V$ and $x_0\in C:=g(\PP^1)$ then $C\subset U_{Z'}$ as $C$ is not contained in the support of $D^{Z'}$ and $D^{Z'}\cdot C=0$. But then Lemma \ref{locus-dimension-P^1} implies that  $\dim {\rm locus} (V, x_0)= r$, a contradiction. It follows that 
 $\varphi_{Z} (x_0)\not \in U_S$ and all fibers of $\varphi_Z|_{U_{I^h}'}: U_{I^h}'\to U_S$ have dimension $r$.
 	
 	\medskip
 		 
 	\emph{Step 2.}  Assume that $J\subsetneq I^h$ and the assertion holds for all subsets properly containing $J$.
 	
 	Let us set  $m=|I^h\backslash J|$ and consider $Z'=D^h_J$. By Proposition \ref{general-fiber-in-general} for general $y\in S$ we have $(Z_y, D_y^Z)\simeq (\PP^{r+s}, H_1+...+H_s)$ and $\cL |_{Z_y}\simeq \cO_{\PP^{r+s} }(1)$. This implies that
  for general $y\in S$ we have $(Z_y', D_y^{Z'})\simeq (\PP^{r+m}, H_1+...+H_m)$ and $\cL |_{Z_y'}\simeq \cO_{\PP^{r+m} }(1)$.	
 Let us consider the Chow variety $\Chow _{r+m,1}(Z')$ of $(r+m)$-dimensional cycles on $Z'$ of degree $1$ with respect to $\cL'$. Let  $Q$ be an irreducible component of $\Chow _{r+m,1}(Z')$ containing general fibers $Z'_y$ of $\varphi_Z|_{Z'}$. As in Step 1 this component is projective of dimension equal to the dimension of $S$ and there exists a universal family of pairs of cycles $\tau: \cU \to Q$ with the evaluation map $\ev : \cU\to Z'$. The map $\ev$ is birational and either it is an isomorphism or it has a positive dimensional fiber. In the first case, $(\varphi _Z|_{Z'}: Z'\to S)\simeq (\tau: \cU\to Q)$, so  $\varphi_Z|_{Z'}$ is equidimensional. 
 	In the second case, there exists $x_0\in Z'$  such that $\dim (\ev ^{-1}(x_0))>0$. Then by the defining property of $\Chow_{r+m,1}(Z')$, $\dim (\ev( \tau ^{-1}(\tau( \ev ^{-1}(x_0)))))>r+m$.

 Let us choose some $f: \PP^1\to \cU$ such that $f$ is an isomorphism onto  a line in a general fiber  $Z'_y=\tau^{-1}([Z'_y])\simeq \PP^{r+m}$, not contained in the support of $D^{Z'}_y$. By $\tilde f$ we denote the composition of $f$ with $\ev$. Let us fix some $i_0\in I^h\backslash J$ and let us set $D_{\cU}:=\ev ^* D^{Z'}_{i_0}$.
 	Let $W_{\cU}\subset \Hom (\PP^1, \cU; f^*D_{\cU}\subset D_{\cU})$ be an irreducible component containing $[f]$. the same arguments as in Step 1 show the family $W_{\cU}$ is unsplit (since it parametrizes birational morphisms to curves of degree $1$ with respect to some ample line bundle) and the natural map $ \tilde \ev: \Hom (\PP^1, \cU; f^*D_{\cU}\subset D_{\cU})\to  \Hom (\PP^1, {Z'}; \tilde f^*D^{Z'}_{i_0}\subset D^{Z'}_{i_0})$ maps $W_{\cU}$ into an irreducible component $W$ of $ \Hom (\PP^1, {Z'}; \tilde f^*D^{Z'}_{i_0}\subset D^{Z'}_{i_0})$ containing $[\tilde f]$.
 	For general $u\in \cU$,  ${\rm locus} (W_{\cU}, u)=\tau ^{-1}(\tau (u))\simeq \PP^{r+m}$. Since $W_{\cU}$ is unsplit, the equality  ${\rm locus} (W_{\cU}, u)=\tau ^{-1}(\tau (u))$ holds for all $u\in \cU$. This implies that  ${\rm locus} (W, x_0)\supset \ev( \tau ^{-1}(\tau( \ev ^{-1}(x_0))))$, so
 	\begin{equation}\label{dim-est-1}
 \dim {\rm locus} (W, x_0)	\ge (r+m+1).
 	\end{equation}
By Lemma \ref{special-Lie-algebroid-structure-v2}, $(Z', D^{Z'})$ is semi-contact of rank $(2r+1)$.
So by Lemma \ref{rank+Lie-algebroids} for any $g: \PP^1\to {Z'}$ with $[g]\in W$ we have
$$\rk ^{+}g^*T_{Z'}(-\log D^{Z'})\le r+1.$$
If $g(\PP^1)\cap \Supp D^v=\emptyset$ then pulling back the short exact sequence
 	$$0\to T_{Z'}(-\log (\sum _{i\in I^h\backslash J} D^{Z'}_i))\longrightarrow T_{Z'}(-\log D^{Z'}_{i_0})\longrightarrow \bigoplus _{i\in I^h\backslash (J\cup \{ i_0\})} \cO_{D^{Z'}_i} (D^{Z'}_i)\to 0$$
by $g$ we get an exact sequence
$$ g^*T_{Z'}(-\log D^{Z'})\longrightarrow g^*T_{Z'}(-\log D^{Z'}_{i_0})\longrightarrow \cO_{\PP^1}(1)^{\oplus (m-1)}\to 0.$$
Then
$$\rk ^{+}g^*T_{Z'}(-\log D^{Z'}_{i_0})\le (r+1)+(m-1)= r+m.$$ 
Proposition \ref{log-bound-on-rank} implies that 
if $x_0\not \in D^{{Z'}}_{i_0}\cup \Supp D^h$ then $\dim {\rm locus} (W, x_0)\le r+m$,
 a contradiction. If $x_0\in  D^{Z'}_{i_0}\backslash \Supp D^h$ then by (\ref{dim-est-1})
 $$\dim {\rm locus} (W, x_0)\cap D^{Z'}_{i_0}\ge \dim {\rm locus} (W, x_0)-1\ge r+m.$$
 But ${\rm locus} (W, x_0)$  is contained in a fiber of $\varphi_Z$ as it is covered by extremal curves passing through one point.  So ${\rm locus} (W, x_0)\cap D^h$ is contained in $(\varphi_Z |_{{Z'}\cap D^h_{i_0}})^{-1}(\varphi_Z (x_0))$, which
 by the induction assumption has dimension $(r+m-1)$, a contradiction.
It follows that 
$\varphi_{Z} (x_0)\not \in U_S$ and all fibers of 
$\varphi _Z |_{U_J'}$ have dimension $(r+m)$.
\end{proof}

\subsection{Classification of good elementary log contractions}

The following proposition is an analogue of the first part of \cite[Theorem 2.12]{KPSW}. As in \cite{KPSW} the result depends on Fujita's theorem  \cite[Lemma 2.12]{Fujita1985}, but the proof is much more involved as we need to analyse the stratification of $Z$ induced by the components of $D^Z$. We keep the notation from the previous subsection.

\begin{Proposition}\label{rough-log-contractions}
	Assume that  $\varphi _Z$ is elementary and $\dim S>0$. Then
	$S$ is smooth and  for any subset $J\subset I^h$ (possibly empty), $\varphi _Z|_{D^h_J}: D^h_J\to S$ is a projective bundle. Moreover,  there exists a simple normal crossing divisor $B^S$ on $S$ such that irreducible components of $D^v$ are preimages of irreducible components of $B^S$.
\end{Proposition}

\begin{proof} 	
Since $\varphi_Z$ is elementary, all extremal rational curves contracted by $\varphi_Z$ are proportional. 
By  Proposition \ref{general-fiber-in-general}, for general $y\in S$ we have $(Z_y, D_y^Z)\simeq (\PP^{r+s},H_1+...+H_s)$, so if $C$ is an extremal rational curve then  in $N_1(Z/S)$ the class of $C$ coincides with the class of a line in $ Z_y\simeq\PP^{r+s}$.
Therefore $C\cdot D_i^Z=1$ for $i\in I^h$ and $C\cdot D_i^Z=0$ for all $i\in I^v$.
In particular, curves contracted by $\varphi_Z$ do not properly intersect any irreducible component of $D^v$.

\begin{Lemma}\label{vertical-saturation} 
	Let $K\subset I^v$ and let $W$ be a connected component of $D_K^Z$.
	Then $W$ is saturated with respect to $\varphi_Z$, that is,
	\[
	\varphi_Z^{-1}\bigl(\varphi_Z(W)\bigr)=W
	\]
	set-theoretically.
	The same assertions hold after restriction to a stratum. More precisely,
	let $Z'=D_J^h$, let $C\subset Z'$ be an extremal rational curve
	contracted by $\varphi_{Z'}:=\varphi_Z|_{Z'}: Z'\to S$, let
	$K=\{i\in I^v:C\subset D_i^Z\}$,
	and let $Z''$ be the connected component of $D_K^{Z'}$ containing $C$. Then
\begin{enumerate}
	\item  $\varphi_{Z'}$ is a good log contraction of   $(Z', D^{Z'})$, 
	\item   $S'':=\varphi_Z(Z'')$ is normal and $\varphi_{Z''}:=\varphi_Z|_{Z''}:Z''\to S''$ is a good log contraction of   $(Z'', D^{Z''})$,
	\item  for every $i\in I^v\backslash K$, every
	nonempty component of $D_i^Z|_{Z''}$ is vertical with respect to
	$\varphi_{Z''}$,
	\item  if
	$B_i^S:=\varphi_Z(D_i^Z)$, then $S''\not\subset B_i^S$ and
	\[
	\operatorname{Supp}\bigl(D_i^Z|_{Z''}\bigr)
	=
	\varphi_{Z''}^{-1}\bigl(S''\cap B_i^S\bigr)
	\]
	set-theoretically.
\end{enumerate}	
\end{Lemma}

\begin{proof}
	 Let $F$ be a fiber of
	$\varphi_Z$ meeting $W$. We claim that $F\subset D_i^Z$ for every
	$i\in K$. Otherwise, since $F$ is connected and projective, we could
	find an irreducible curve in $F$ properly intersecting $D_i^Z$. Since this is
	impossible we have $F\subset D_K^Z$. Since
	$F$ is connected and meets the connected component $W$, we have
	$F\subset W$. It follows that
	\[
	\varphi_Z^{-1}\bigl(\varphi_Z(W)\bigr)=W
	\]
	set-theoretically.
	
	Points (1) and (2) follow from Corollary \ref{normality-of-image}. The same argument as above applies to the contraction
	$\varphi_{Z'}:Z'\to S$. Let $i\in I^v\backslash K$, the saturation
	statement gives
	\[
	\operatorname{Supp}\bigl(D_i^Z|_{Z''}\bigr)
	=
	Z''\cap\varphi_Z^{-1}(B_i^S)
	=
	\varphi_{Z''}^{-1}(S''\cap B_i^S).
	\]
	If $S''\subset B_i^S$, then the saturation of $D_i^Z$ would imply
	$Z''\subset D_i^Z$, and hence $C\subset D_i^Z$, contradicting
	$i\notin K$. Therefore $S''\cap B_i^S$ is a proper closed subset of
	$S''$, and every component of $D_i^Z|_{Z''}$ is vertical with respect
	to $\varphi_{Z''}$.
\end{proof}

\medskip

As in the proof of Proposition \ref{equidimensionality-on-open} the proof is by induction on the number of elements in $I^h\backslash J$.

\medskip

\emph{Step 1.} Let us first consider the case $J=I^h$. Let us set
$Z':=D^h_{I^h}$.  By
Proposition \ref{general-fiber-in-general}, for general $y\in S$ we
have $(Z'_y,D_y^{Z'})\simeq(\PP^r,0)$ and
$\cL|_{Z'_y}\simeq\cO_{\PP^r}(1)$.
Let $F$ be an arbitrary fiber of $\varphi_{Z'}$, taken with the reduced structure. Since the general
fiber has dimension $r$, we have
$\dim F\geq r$. Choose a projective curve $C'$ contained in an
irreducible component of $F$ of maximal dimension. By Proposition
\ref{cone-theorem-complement}, we can find an extremal rational curve
$C$ contracted by $\varphi_{Z'}$, lying in the same fiber and in the
same stratum of the stratification of $Z'$ induced by the components
of $D^{Z'}$. Let $f:\PP^1\to C$ be the normalization of $C$. Let $K$
be the set of all $j\in I^v$ for which $C\subset D_j^Z$, and let
$Z''$ be the connected component of $D_K^{Z'}$ containing $C$.

By Lemma \ref{vertical-saturation}, $S'':=\varphi_Z(Z'')$ is normal
and $\varphi_{Z''}:=\varphi_Z|_{Z''}:Z''\to S''$ is a good log
contraction of $(Z'',D^{Z''})$. Moreover, all components of
$D^{Z''}$ are vertical with respect to $\varphi_{Z''}$. By Lemma
\ref{special-Lie-algebroid-structure-v2}, $(Z'',D^{Z''})$ is
semi-contact of rank $2r+1$, and hence Corollary \ref{fiber-type}
gives $\dim S''<\dim Z''$.

By the definition of $K$, the curve $C$ is not contained in
$\operatorname{Supp}D^{Z''}$. Lemma \ref{vertical-saturation} therefore
implies that $C\subset U_{Z''}:=\varphi_{Z''}^{-1}(U_{S''})$. If
$\dim S''>0$, Proposition \ref{equidimensionality-on-open} shows that
the fiber of $\varphi_{Z''}$ containing $C$ has dimension $r$. On the
other hand, Lemma \ref{vertical-saturation} shows that the fiber $F$
is contained in $Z''$. Consequently, $F$ is the fiber of
$\varphi_{Z''}$ containing $C$, and hence $\dim F=r$.

\begin{Lemma}	\label{contraction-of-strata-to-points}
	If $\dim S''=0$ then $(Z'',D^{Z''}) \simeq (\PP^r,0)$.
\end{Lemma}

\begin{proof}
	Vanishing $D^{Z''}=0$ is clear. Otherwise, 
	Lemma \ref{vertical-saturation} gives
	$Z''=\varphi_{Z''}^{-1}(\varphi_{Z''} (D^{Z''}))=D^{Z''}$, 	which is impossible because $D^{Z''}$ is a divisor in $Z''$.
	
	By construction $-K_{Z''}=(r+1)c_1(\cL|_{Z''})$ is ample, so $Z''$ is Fano. 
	If	 $b_2(Z'')>1$ then by \cite[Theorem B]{Wisniewski1990} $Z''$ is isomorphic to $\PP^r\times \PP^r$. 
	Note that $$\dim Z''\le \dim Z'-1\le \dim Z-1\le
	2r.$$	Since $\dim Z''=2r$, we have $Z'=Z$  and $Z''$ is a divisor in $Z$.
	By Lemma \ref{vertical-saturation}, $Z''$ is the geometric fiber of 
	$\varphi_Z$ (taken with the reduced structure) over the point $\varphi_Z(Z'')$. Every curve contained in
	$Z''$ is contracted by $\varphi_Z$, so the conormal line bundle of $Z''$ has degree zero 
	on both rulings of $\PP^r\times\PP^r$. Hence
	\[
	N^*_{Z''/Z}\simeq\cO_{Z''}.
	\]
	
	Since $H^1(Z'',\cO_{Z''})=0$, the assumptions of Castelnuovo's
	criterion \cite[Proposition 3.8]{Andreatta-Wisniewski1997} are
	satisfied. Therefore $\dim S=h^0(Z'',N^*_{Z''/Z})=1$.
	On the other hand, the general fiber of $\varphi_Z$ has dimension $r$,
	so $\dim S=\dim Z-r=r+1$, a contradiction.
	
	It follows that $b_2(Z'')=1$. Since $\dim Z''\le 2r$ and $D^{Z''}=0$,  Proposition \ref{dimension-contact-Lie-algebroid} implies that $Z''\simeq \PP^r$.
\end{proof}

If $\dim S''=0$ then Lemma \ref{vertical-saturation} shows that $Z''=F$. So the above lemma implies that also in this case $\dim F=r$.  Thus every fiber of $\varphi_{Z'}$ has dimension $r$. Fujita's theorem
\cite[Lemma 2.12]{Fujita1985} now implies that $S$ is smooth and
$\varphi_{Z'}$ is a projective bundle.

\medskip

For every $i\in I^v$, let $B_i^S:=\varphi_Z(D_i^Z)$ (with the reduced structure).
For a subset $K\subset I^v$, we denote by $B_K^S$ the scheme-theoretic
intersection $\bigcap_{i\in K}B_i^S$.
Lemma \ref{vertical-saturation} implies that
$\varphi_{Z'}^{-1}(B_i^S)=D_i^{Z'}$ set-theoretically. Since
$\varphi_{Z'}$ is smooth and $B_i^S$ is reduced, its inverse image is
reduced. Thus this equality is scheme-theoretic. Consequently,
inverse images commute with scheme-theoretic intersections and hence
$\varphi_{Z'}^{-1}(B_K^S)=D_K^{Z'}$
scheme-theoretically. The pair $(Z',D^{Z'})$ is snc, so for every $K\subset I^v$ the intersection $D_K^{Z'}$ is smooth and, whenever
it is nonempty, it has codimension $|K|$ in $Z'$. Since
$\varphi_{Z'}$ is smooth this implies that $B_K^S$ is smooth and
$\operatorname{codim}_S B_K^S
=
\operatorname{codim}_{Z'}D_K^{Z'}
=
|K|.$
In particular, every $B_i^S$ is a smooth irreducible divisor, and
every nonempty intersection $B_K^S$ is smooth of the expected
codimension. Therefore
$B^S:=\sum_{i\in I^v}B_i^S$ is a simple normal crossing divisor on
$S$.

\medskip

\emph{Step 2.} Assume that $J\subsetneq I^h$ and that the assertion
holds for every subset properly containing $J$. Let
$m:=|I^h\backslash J|$ and set $Z':=D_J^h$. By Lemma
\ref{vertical-saturation},
$\varphi_{Z'}:=\varphi_Z|_{Z'}:Z'\to S$ is a good log contraction of
$(Z',D^{Z'})$. By Proposition \ref{general-fiber-in-general}, for
general $y\in S$ we have
$(Z'_y,D_y^{Z'})\simeq(\PP^{r+m},H_1+\cdots+H_m)$ and
$\cL|_{Z'_y}\simeq\cO_{\PP^{r+m}}(1)$.

Let $F$ be an arbitrary fiber of $\varphi_{Z'}$ and let $F_0$ be an
irreducible component of $F$ of maximal dimension. Since the general
fiber has dimension $r+m$, we have 
$\dim F_0\geq r+m$. If $F_0$ were contained in a horizontal
component $D_{i_0}^{Z'}$ of $D^{Z'}$, then the induction assumption
applied to $J\cup\{i_0\}$ would give $\dim F_0\leq r+m-1$, a
contradiction. Thus $F_0$ is not contained in any horizontal component
of $D^{Z'}$.

Choose a projective curve $C'\subset F_0$ which is not contained in
any horizontal component of $D^{Z'}$. By Proposition
\ref{cone-theorem-complement}, we can find an extremal rational curve
$C$ contracted by $\varphi_{Z'}$, lying in the same fiber and in the
same stratum as $C'$. In particular, $C$ is not contained in any
horizontal component of $D^{Z'}$. Let $f:\PP^1\to C$ be the
normalization of $C$, let $K$ be the set of all $j\in I^v$ for which
$C\subset D_j^Z$, and let $Z''$ be the connected component of
$D_K^{Z'}$ containing $C$.

By Lemma \ref{vertical-saturation}, $S'':=\varphi_Z(Z'')$ is normal
and $\varphi_{Z''}:=\varphi_Z|_{Z''}:Z''\to S''$ is a good log
contraction of $(Z'',D^{Z''})$. By Lemma
\ref{special-Lie-algebroid-structure-v2}, $(Z'',D^{Z''})$ is
semi-contact, and Corollary \ref{fiber-type} gives
$\dim S''<\dim Z''$.

Let $D^{Z'',v}$ be the vertical part of $D^{Z''}$ with respect to
$\varphi_{Z''}$. By Lemma \ref{vertical-saturation}, every component
of $D^{Z''}$ induced by a vertical component of $D^Z$ is vertical
with respect to $\varphi_{Z''}$. The components indexed by
$I^h\backslash J$ are horizontal, so the horizontal part
$D^{Z'',h}$ has exactly $m$ irreducible components.

Lemma \ref{vertical-saturation} shows that $F\subset Z''$ and that
$C\subset U_{Z''}:=\varphi_{Z''}^{-1}(U_{S''})$. If $\dim S''>0$,
Proposition \ref{equidimensionality-on-open} shows that the fiber of 
$\varphi_{Z''}$ containing $C$ has dimension $r+m$. Since this fiber
is $F$, we obtain $\dim F=r+m$.

\begin{Lemma}	\label{contraction-of-strata-to-points-v2}
	If $\dim S''=0$ then $(Z'', D^{Z''})\simeq (\PP^{r+m}, H_1+...+H_m)$. 
\end{Lemma}

\begin{proof}	
	Since $\dim S>0$, the set $K$ is nonempty. Therefore
	$\dim Z''\leq\dim Z'-1\leq\dim Z-1\leq 2r$.
	Suppose that $\rho (Z'')>1$. Clearly, $(Z'', D^{Z''})$ is log Fano and by Lemma
	\ref{existence-of-lines}, the log Fano pseudoindex satisfies
	\[
	\imath(Z'',D^{Z''})=r+1>\frac{\dim Z''}{2}.
	\]
	By \cite[Theorem 4.3]{Fujita-Log-Fano}, this would imply $\dim Z''=2\imath(Z'',D^{Z''})-1=2r+1$,
	contradicting $\dim Z''\leq2r$. Therefore $\rho (Z'')=1$. Then equality
	$-K_{Z''}=(r+1)c_1(\cL|_{Z''})+D^{Z''}$
	shows that $-K_{Z''}$ is ample. Thus $Z''$ is Fano and
	$\Pic Z''\simeq\ZZ$. Moreover, if $C$ is the extremal rational curve
	used in the construction of $Z''$, then $\cL\cdot C=1$. Therefore
	$\cL|_{Z''}$ is the ample generator of $\Pic Z''$.
	Since 
	$\cL\cdot C=1$ and  $D_i^{Z'}\cdot C=0$ for every $i\in K$, it follows that 
	$ \cO_{Z''}(-D_i^{Z'})\simeq\cO_{Z''}$
	for every $i\in K$.

	Lemma \ref{vertical-saturation} also shows that $Z''$ is the geometric fiber  of $\varphi_{Z'}$ over the point
	$\varphi_{Z'}(Z'')$. Since the boundary has simple normal crossings,
	\[
	N^*_{Z''/Z'}
	=
	I_{Z''}/I_{Z''}^2
	\simeq
	\bigoplus_{i\in K}\cO_{Z''}(-D_i^{Z'})
	\simeq
	\cO_{Z''}^{\oplus(\dim Z'-\dim Z'')}.
	\]
	Since $Z''$ is Fano, Kodaira vanishing gives
	$H^1(Z'',\cO_{Z''})=0$. Hence all the assumptions of Castelnuovo's
	criterion \cite[Proposition 3.8]{Andreatta-Wisniewski1997} are
	satisfied, and
	\[\dim S	=	h^0(Z'',N^*_{Z''/Z'})=	\dim Z'-\dim Z''.	\]
	On the other hand, the general fiber of $\varphi_{Z'}$ has dimension
	$r+m$, so
	$	\dim Z''=\dim Z'-\dim S=r+m$.
	Let us write $D_i^{Z''}=a_ic_1(\cL|_{Z''})$ for $i\in I^h\backslash J$ and some positive integers $a_i$. Then 
	$$-K_{Z''}=\left(r+1+\sum_{i\in I^h\backslash J} a_i\right)c_1(\cL|_{Z''}).$$
	Since there are $m$ integers $a_i\geq1$, we have
	\[
	r+1+\sum_{i\in I^h\backslash J}a_i
	\geq r+m+1
	=
	\dim Z''+1.
	\]
	The Kobayashi--Ochiai theorem implies that equality holds. Therefore
	$a_i=1$ for every $i\in I^h\backslash J$ and
	$ (Z'',\cL|_{Z''})
	\simeq
	(\PP^{r+m},\cO_{\PP^{r+m}}(1)).$
	The components of $D^{Z''}$ are consequently hyperplanes, and hence
	$(Z'',D^{Z''})
	\simeq
	(\PP^{r+m},H_1+\cdots+H_m)$.
\end{proof}

If $\dim S''=0$ then Lemma \ref{vertical-saturation} shows that $Z''=F$.
So the above lemma implies that $\dim F=r+m$.
Thus all fibers of $\varphi_{Z'}$ have dimension $r+m$. Then Fujita's
theorem \cite[Lemma 2.12]{Fujita1985} implies that $\varphi_{Z'}$ is a projective bundle.

This completes the induction. In particular, for $J=\emptyset$ the
map $\varphi_Z:Z\to S$ is a projective bundle and hence is smooth.
Lemma \ref{vertical-saturation} then gives, scheme-theoretically,
$D_i^Z=\varphi_Z^{-1}(B_i^S)$ for every $i\in I^v$. The divisor
$B^S$ is the simple normal crossing divisor constructed in Step~1.
\end{proof}
 
\medskip

\begin{Example}\label{construction}  
Let us fix some $1\le s\le r$, a projective snc pair $(S, B^S)$.
Assume that $1\le \dim S\le (r-s+1)$. Let $\underline{\cM}=(\cM_1, ...,\cM_s)$ be a collection of 
 $s$ line bundles  on $S$ and let
$$0\to \cO_S^{\oplus s}\to \cC = \cD_{T_{S}(-\log B^S)} ^{\le 1}(\underline{\cM})\to T_{S}(-\log B^S)\to 0$$
be the  extension defined in Subsection \ref{Atiyah-classes}, dual to the logarithmic multi-Atiyah extension of $\cM$
(so that $\cC= (\cB_{\underline{\cM}})^*$).
Let us consider extensions 	
$$0\to \cO_S^{\oplus (r+1-s-\dim S)}\to \cG \to \cC \to 0$$
(so that $\cG$ has rank $r+1$) and
$$0\to \bigoplus _{i=1}^s\cM_i\mathop{\longrightarrow}^{\oplus \gamma _i} \cE \longrightarrow\cG  \to 0.$$
Let $\cG_i$ be the cokernel of $\gamma_i: \cM_i\to \cE$. 
Then on $Z:=\PP_S(\cE)$ we can consider $s$ horizontal divisors $D^h_i:=\PP(\cG _i)$. If $\varphi_Z: Z\to S$ is the canonical projection then the divisor $D^Z:=\sum _{i=1}^sD^h_i+ \varphi_Z^{-1}(B^S)$ has simple normal crossings.
Note also that  $\dim Z=\dim S+s+r$. In particular, if $\dim Z=2r+1$ then $\cG=\cC$.
\end{Example}

\medskip

\begin{Theorem}\label{structure-of-log-contractions-special-Lie-alg}
Let  $(Z, D^Z, \cL)$ be a rank $(2r+1)$ projective semi-contact snc pair.  Let  $\varphi_Z: Z\to S$ be a very good elementary log contraction with  $\dim S>0$. Then $S$ is smooth, $\cE=\varphi_{Z*}\cL$ is locally free, $(Z, \cL)\simeq (\PP(\cE), \cO_{\PP(\cE)}(1))$, and there exists a simple normal crossing divisor $B^S$ on $S$ such that one of the following cases occur:
\begin{enumerate}
		\item  $\rk \cE=r+1$, $D^Z=\varphi_Z^{-1}(B^S)$ and we have a short exact sequence
		$$0\to \cO_S^{\oplus (r+1-\dim S)}\to \cE \to T_{S}(-\log B^S) \to 0.$$
		\item  $\rk \cE=r+s+1$ for some $1\le s\le r$, and $\cE$ and $D^Z$ are as in Example \ref{construction}.
	\end{enumerate}
\end{Theorem}
 
\begin{proof}
Proposition  \ref{rough-log-contractions} implies that $S$ is smooth and $\varphi_Z$ is a projective bundle.
So Proposition \ref{general-fiber-in-general} implies that $\cE=\varphi_{Z*}\cL$ is locally free and $(Z, \cL)\simeq (\PP(\cE), \cO_{\PP(\cE)}(1))$. It also implies that there exists $0\le s\le r$ such that $\cE$ has rank $(r+s+1)$ and $D^Z$ has exactly $s$ horizontal components $D^h_1,...,D_s^h$.  The relative Euler exact sequence for $\varphi_Z$ shows that
$\omega_{Z/S}=\varphi_Z^*(\det \cE)\otimes \cL^{-(r+s+1)}$ and hence
$\omega_{Z}\otimes \cL^{r+s+1}=\varphi_Z^*(\omega_S\otimes \det \cE).$
Since  $\omega_Z(D^Z)\simeq \cL^{-(r+1)}$, we get 
$$\cL^s(-D^Z)\simeq \varphi_Z^*(\omega_S\otimes \det \cE).$$
By Proposition  \ref{rough-log-contractions} we also know that there exists a simple normal crossing divisor $B^S$ on $S$ such that 
$\varphi_Z^{-1}(B^S)$ is the vertical part of $D^Z$.
In the rest of the proof we set $n:=\dim Z$ and $q:=2r+1-n$. 
Since $n=\dim S+r+s$, we have $q=r+1-s-\dim S$.

Let $(\LL, \cL)$ be a special contact Lie algebroid on $(Z,D^Z)$.
Let us set $\cB=  \cA_{\LL}(\cL)$ and consider the $\LL$-Atiyah extension of $\cL$
$$0\to \Omega_{\LL}^1\longrightarrow \cB \longrightarrow \cO_Z\to 0.$$
We also consider the logarithmic Atiyah bundle \(\cB_{\cL}\) of \(\cL\) on
\((Z,D^Z)\) and  the relative logarithmic Atiyah bundle \(\cB_{\cL/S}\) of \(\cL\) on
\((Z,D^Z)/(S, B^S)\) (see Subsection \ref{subsection-Atiyah-diff}).
 By Lemma \ref{seq:A-B-D^h} we have short exact sequences
 \begin{equation}\label{new-seq-1}
 	0\longrightarrow \varphi_Z^*\Omega_S^1(\log B^S)\longrightarrow \cB_{\cL}\longrightarrow \cB_{\cL/S}\longrightarrow 0
 \end{equation}
 and
 \begin{equation}\label{new-seq-2}
 	0\longrightarrow \varphi_Z^*\cE(-1)\longrightarrow \cB_{\cL/S}\longrightarrow \bigoplus_{i=1}^s\cO_{D_i^h}\longrightarrow 0,
 \end{equation}
 where the direct sum in the second sequence is zero when \(s=0\). Since \(\varphi_Z\) and every restriction \(\varphi_Z|_{D_i^h}:D_i^h\to S\) are projective bundles, the projection formula and standard vanishings  give
 \[
 R^1\varphi_{Z*}\bigl(\varphi_Z^*\cE(-1)\bigr)=\cE\otimes  R^1\varphi_{Z*}\bigl(\cO_Z(-1)\bigr)=0
 \qquad\text{and}\qquad
 R^1\varphi_{Z*}\cO_{D_i^h}=0.
 \]
 Therefore \(R^1\varphi_{Z*}\cB_{\cL/S}=0\). 
 Similarly,
$R^1\varphi_{Z*}\bigl(\varphi_Z^*\Omega_S^1(\log B^S)\bigr)\simeq \Omega_S^1(\log B^S)\otimes R^1\varphi_{Z*}\cO_Z=0$,
 so the first exact sequence gives
 \[
 R^1\varphi_{Z*}\cB_{\cL}=0.
 \]

Let us set $\cF_1:=\varphi_{Z*}\cB_{\cL}$ and $\cF_2:=\varphi_{Z*}(\cB_{\cL}^*(-1))$.
By Lemma \ref{Atiyah-bundles-sequence} we have the short exact sequence
$$0\to \cB_{\cL}\to \cB\to  \cO_Z^{\oplus q} \to 0.$$
So pushing down the above sequence gives
\[ 0\longrightarrow \cF_1 \longrightarrow \varphi_{Z*}\cB \longrightarrow
\cO_S^{\oplus q} \longrightarrow 0 .\]
Pushing down the dual of above sequence twisted by $\cO_Z(-1)$ gives  $\varphi_{Z*}(\cB^*(-1))\simeq\cF_2$.
 
Let us recall that by Proposition \ref{general-LeBrun-lemma} existence  of the contact structure on $(\LL, \cL)$ implies existence of an $\cL$-valued symplectic form on $\cB^*$. In particular, we have an isomorphism 
 $\cB\simeq \cB^*(-1).$ Therefore $\cF_2\simeq \varphi_{Z*}\cB $ and we have a short exact sequence
\begin{equation}\label{F_1-F_2}
0\longrightarrow
\cF_1 \longrightarrow \cF_2 \longrightarrow \cO_S^{\oplus q} \longrightarrow 0.
\end{equation}
It remains only to identify \(\cF_1\) and \(\cF_2\).  To do that we consider the following two cases.
 
\medskip 

\emph{Case 1: }{$s=0$.}
 
In this case sequences (\ref{new-seq-1}) and (\ref{new-seq-2}) reduce to  
 \begin{equation*}
 	0\longrightarrow \varphi_Z^*\Omega_S^1(\log B^S)\longrightarrow \cB_{\cL}\longrightarrow \varphi_Z^*\cE(-1)\longrightarrow 0.
 \end{equation*}
Pushing it down to $S$, we obtain a canonical isomorphism $\Omega^1_S(\log B^S)\simeq \cF_1$.
Pushing down the dual of above sequence twisted by $\cO_Z(-1)$ gives  $\cE^*\simeq \cF_2$. 
 So dualizing  (\ref{F_1-F_2}) we get the short exact sequence
 $$0\to \cO_S^{\oplus q}\to \cE \to T_S(-\log B^S)\to 0.$$
 Since $q=r+1-\dim S$,  this gives the first case in the theorem.

\medskip

\emph{Case 2:} {$1\le s\le r$. } 

Let us set $\cM:=\omega_S\otimes \det \cE(B^S)$ and $D^h=\sum _{i=1}^sD^h_i$. Then $\cL^s(-D^h)\simeq \varphi_Z^*\cM$.

Let us set $J=\{1,...,s\}$ and $\cG= \varphi_{Z*}(\cL|_{D^h_J})$. Note that $\cG$ is locally free of rank $(r+1)$ and  $(D^h_J, \cL|_{D^h_J})\simeq (\PP (\cG), \cO_{\PP(\cG )}(1))$.   
By Subsection \ref{hyperplane-subbundles} there exist line bundles $\cM_i$ such that $\cL(-D^h_i)\simeq \varphi_Z^*\cM _i$ for $i=1,...,s$. Let $\gamma_i: \cM_i\to \cE=\varphi_{Z*}\cL$ be  the canonical injective map  corresponding to  $\varphi_Z^*\cM _i\to \cL$ and let $\cG_i$ be the cokernel of $\gamma_i$. Then $\cG_i$
is locally free and the projection $\cE\to \cG_i$ defines the embedding $D^h_i\simeq \PP(\cG_i)\hookrightarrow Z\simeq \PP(\cE)$.

\begin{Lemma}\label{seq-M-E-G}
	We have a short exact sequence
	$$0\to \bigoplus _{i\in J}\cM_i\mathop{\longrightarrow}^{\oplus \gamma _i} \cE \longrightarrow\cG  \to 0.$$
\end{Lemma}

\begin{proof}
	Note that if $H_1,...,H_s$ are hyperplanes in general position in $\PP^{r+s}$ then
	$H^1(J_{H_1\cap ...\cap H_s}(1))=0.$
	Hence $$R^1\varphi_{Z*}(J_{D^h_J}\otimes \cL )=0.$$ 
	Similarly, vanishings
	$H^0(J_{H_1\cap ...\cap H_s}^2(1))=H^1(J_{H_1\cap ...\cap H_s}^2(1))=0$
	imply that $$\varphi_{Z*}(J^2_{D_J^h}\otimes\cL)=R^1\varphi_{Z*}(J^2_{D_J^h}\otimes\cL)=0.$$
	
	We claim that the required sequence is the push-forward of
	the short exact sequence
	$$0\to J_{D^h_J}\otimes \cL \to \cL\to \cL|_{D^h_J}\to 0.$$
	By the above we only need to check that $\varphi_{Z*}(J_{D^h_J}\otimes \cL )\simeq \bigoplus _{i\in J}\cM_i$
	and the push-forward map is given by $\oplus _{i\in J} \gamma _i$.
	
	The push-forward of the sequence
	$$0\to  J_{D^h_J}^2\otimes \cL\to  J_{D^h_J}\otimes \cL\to  J_{D^h_J}/ J_{D^h_J}^2 \otimes \cL \to 0$$
	gives an isomorphism  $\varphi_{Z*}(J_{D^h_J}\otimes \cL )\to \varphi_{Z*}(J_{D^h_J}/ J_{D^h_J}^2\otimes \cL )$.
	We also have a surjective map
	$$   \bigoplus _{i\in J}\cO _Z(-D^h_i)\otimes \cL\simeq \varphi_Z^*\left( \bigoplus _{i\in J}\cM_i \right)\to J_{D^h_J}\otimes \cL ,$$
	which after tensoring with $\cO_{D^h_J}$ gives an isomorphism
	$$  ( \varphi_{Z}|_{D^h_J})^*\left( \bigoplus _{i\in J}\cM_i \right)\to J_{D^h_J}/ J_{D^h_J}^2\otimes \cL .$$
	Taking the pushforward gives the required assertions.
\end{proof}

The following lemma is a generalization of Lemma \ref{relative-cotangent-bundle}. It also generalizes \cite[Proposition 2.10]{Dolgachev-Kapranov} that corresponds to the case when $S$ is a point 
(our approach gives also a canonical sequence explaining the isomorphism).

\begin{Lemma}\label{relative-cotangent-bundle-v2} 
	The relative evaluation map $\varphi_{Z}^*\varphi_{Z*} \cB_{\cL/S}\to \cB_{\cL/S}$ is injective and its cokernel is isomorphic to $\varphi_Z^*\cG(-1)$.
	Moreover, $\varphi_{Z*} \cB_{\cL/S}$ is canonically isomorphic to $ \cO_S^{\oplus s}$.
\end{Lemma} 

\begin{proof}  
	The relative Atiyah bundle $\cA_{\cL/S}$ is isomorphic to $\varphi _Z^*\cE (-1)$ (this follows from the relative Euler sequence). So by Lemma \ref{seq:A-B-D^h} we have the following short exact sequence:
	$$0\to\varphi _Z^*\cE (-1)\to \cB_{\cL/S}\to {\bigoplus _{i=1}^s}\cO_{D^h_i}\to 0.$$
	Pushing it down we see that $\varphi_{Z*} \cB_{\cL/S}\to \varphi _{Z*}\left( {\bigoplus_{i=1}^s}\cO_{D^h_i}\right)=\cO_S^{\oplus s}$
	is an isomorphism. So we have a well-defined map
	$${\bigoplus_{i=1}^s}\cO_{Z}(-D^h_i)\to {\bigoplus_{i=1}^s}\cO_{Z}\simeq \varphi_{Z}^*\varphi_{Z*}\cB_{\cL/S}.$$
	On the other hand, the twisted pull-back of  $\oplus \gamma _i$ defines the map
	$${\bigoplus_{i=1}^s}\cO_{Z}(-D^h_i)\simeq \varphi_Z^*\left( \bigoplus _{i\in J}\cM_i \right) (-1)\to \varphi _Z^*\cE (-1),$$
	whose cokernel is isomorphic to $\varphi_Z^*\cG(-1)$  by Lemma \ref{seq-M-E-G}.
	Using the sum of these maps and the relative evaluation we can form the sequence
	$$ 0 \to {\bigoplus_{i=1}^s}\cO_{Z}(-D^h_i)	\to \varphi_{Z}^*\varphi_{Z*} \cB_{\cL/S}\oplus \varphi _Z^*\cE (-1)\to \cB_{\cL/S}\to 0.$$
	It is exact as can be easily checked fiberwise on the fibers of $\varphi_Z$. Therefore
	we have  the induced short exact sequence
	$$0\to \varphi_{Z}^*\varphi_{Z*} \cB_{\cL/S}\to \cB_{\cL/S}\to \varphi_Z^*\cG(-1)\to 0.$$
\end{proof}

The following lemma and its proof are analogous to that of Lemma \ref{exact-sequence-without-contact}.   

\begin{Lemma}\label{exact-sequence-without-contact-v2}	 
	Let $\underline{\cM}=(\cM_1,...,\cM_s)$ and let  
	$$0\to  \Omega^1_S(\log B^S)\to \cB_{\underline{\cM}} \to \cO_S^{\oplus s}\to 0$$
	be the logarithmic multi-Atiyah extension of $\underline{\cM}$ (see Subsection \ref{subsection-Atiyah-diff}). 
	Then we have a canonical short exact sequence
	$$0\to \varphi_Z^* \cB_{\underline{\cM}}\to \cB_{\cL}\to \varphi_Z^*\cG(-1)\to 0.$$
	In particular, $\varphi_{Z*}\cB_{\cL}\simeq \cB_{\underline{\cM}}$ and the first map can be identified with the relative evaluation map $\varphi_{Z}^*\varphi_{Z*}\cB_{\cL}\to \cB_{\cL}$.
\end{Lemma} 

\begin{proof}
	Since  $\cL(-D^h_i)\simeq \varphi_Z^*\cM _i$, as in the proof of Lemma \ref{exact-sequence-without-contact}
	we have a commutative diagram
	$$\xymatrix{
		0\ar[r]& \varphi_Z^*\Omega^1_S(\log B^S)\ar[r]\ar[d]^{d\varphi_Z}&  \varphi_Z^*\cB  _{\cM _i}\ar[d]^{\alpha_i}\ar[r]&  \varphi_Z^* \cO_S\ar[r]\ar@{=}[d]&0\\
		0\ar[r]&\Omega_Z^1(\log D^Z)\ar[r]& \cB _{\cL}\ar[r]& \cO_Z\ar[r]&0\\
	}$$
	This gives rise to the commutative diagram
	$$\xymatrix{
		&(\Omega^1_S(\log B^S))^{\oplus s }\ar@/^1pc/[rrr]\ar[r]\ar[d]&  {\bigoplus}_{i=1}^s\cB  _{\cM _i}\ar[d]\ar[rr]_{\bigoplus _{i=1}^s\varphi _{Z*}\alpha_i}&& (\varphi _{Z*}\cB_{\cL})^{\oplus s}\ar[d]\\
		&\Omega^1_S(\log B^S)\ar[r]\ar@/_1pc/[rrr]& \cB_{\underline{\cM}}\ar@{.>}[rr]&& \varphi _{Z*}\cB_{\cL}\\
	}$$
	in which the first and third vertical arrows are codiagonals and
	the dotted arrow exists by the universal property of the pushout (note that the top curved arrow does not depend on $\alpha_i$ and it is equal to the direct sum of $s$ copies of the bottom curved arrow).
	So by construction and Lemma \ref{seq:A-B-D^h} we have a commutative diagram
	$$\xymatrix{
		0\ar[r]& \varphi_Z^*\Omega^1_S(\log B^S)\ar[r]\ar@{=}[d]&  \varphi_Z^*\cB  _{\underline{\cM}}\ar[d]\ar[r]&  \varphi_Z^* (\cO_S^{\oplus s})\ar[r]\ar[d]&0\\
		0\ar[r]&\varphi_Z^*\Omega^1_S(\log B^S) \ar[r]& \cB _{{\cL}}\ar[r]&  \cB _{{\cL}/S}\ar[r]&0\\
	}$$
	in which the last  vertical arrow can be identified with the  relative evaluation map $\varphi_{Z}^*\varphi_{Z*} \cB_{\cL/S}\to \cB_{\cL/S}$. So by  Lemma \ref{relative-cotangent-bundle-v2}
	we get a short exact sequence
	$$0\to \varphi_Z^* \cB_{\underline{\cM}}\to \cB_{\cL}\to \varphi_Z^*\cG(-1)\to 0.$$
	Taking pushforward we get  $ \cB_{\underline{\cM}}\simeq \varphi_{Z*}\cB_{\cL}$, which implies also the last part of the lemma.
\end{proof}

Lemma \ref{exact-sequence-without-contact-v2} shows that $$\cF_1=\varphi_{Z*}\cB_{\cL}\simeq \cB_{\underline{\cM}}.$$ 
Dualizing the sequence from Lemma \ref{exact-sequence-without-contact-v2}, twisting it by $\cO_Z(-1)$ and pushing it down to $S$, we obtain 
$$\cF_2=\varphi_{Z*}(\cB_{\cL}^*(-1))\simeq \cG^*.$$ Therefore, dualizing sequence~(\ref{F_1-F_2}), we obtain a short exact sequence $$0\to \cO_S^{\oplus q}\to \cG\to \cB_{\underline{\cM}}^*\to 0.$$ By definition, $\cB_{\underline{\cM}}^*=\cD_{T_S(-\log B^S)}^{\leq 1}(\underline{\cM})=\cC$, and $q=r+1-s-\dim S$. Thus we obtain the first extension from Example \ref{construction}. Lemma \ref{seq-M-E-G} gives the second extension, and the description $D^Z=\sum_{i=1}^sD_i^h+\varphi_Z^{-1}(B^S)$, where $D_i^h=\PP(\cG_i)$, identifies the boundary with the divisor from that example. This proves the second case and finishes the proof of Theorem \ref{structure-of-log-contractions-special-Lie-alg}.
\end{proof}

\medskip

As a corollary we get the following classification of semi-contact varieties:

\begin{Corollary}\label{classification-semi-contact-varieties}
	Let $X$ be a smooth projective semi-contact variety of rank $(2r+1)$.
	Then one of the following occurs:
	\begin{enumerate}
		\item $X$ is contact and $\dim X=2r+1$,
		\item $X\simeq\PP^{r}$,
	\item $K_X\sim 0$ and $r\le \dim X\le 2r$,
		\item $r+1\leq\dim X\leq 2r$ and there exist a smooth projective variety $S$ with $1\leq\dim S\leq r$ and a vector bundle $\cE$ of rank $r+1$ on $S$ fitting into a short
		exact sequence
		\[
		0\longrightarrow
		\cO_S^{\oplus(r+1-\dim S)}
		\longrightarrow
		\cE
		\longrightarrow
		T_S
		\longrightarrow0
		\]
		such that
		$X\simeq\PP_S(\cE)$.
	\end{enumerate}
\end{Corollary}

\begin{proof}
	Let $(\LL,\cL)$ be a special contact Lie algebroid of rank $(2r+1)$ on $X$.
	
If $K_X$ is pseudoeffective then  Corollary \ref{semi-contact-Kodaira} 
implies that 	$K_X\sim 0$ and	$\dim X<2r+1$. On the other hand, Lemma \ref{rank-estimate}  gives $r\leq\dim X$,
so we are in case~\textup{(3)}.
	
Now assume that $K_X$ is not pseudoeffective and hence $K_X$ is	not nef.
	
If $b_2(X)=1$ then Proposition \ref{dimension-contact-Lie-algebroid}, applied with $D=0$, gives two possibilities. In the first one, $\dim X=2r+1$ and $X$ is contact, which is case~\textup{(1)}. In the second one, for some $0\leq s\leq r$, $X\simeq\PP^{r+s}$
and $D$ is an arrangement of $s$ hyperplanes in general position. Since $D=0$, we must have $s=0$, and hence $ X\simeq\PP^r$.
Thus case~\textup{(2)} occurs.
	
If $b_2(X)\geq2$ then by the cone theorem, there exists an elementary contraction $\varphi\colon X\longrightarrow S$
with $\dim S>0$. Since the boundary is empty, we can apply Theorem \ref{structure-of-log-contractions-special-Lie-alg}.
Case~\textup{(2)} of that theorem cannot occur, because in that case the boundary contains at least one horizontal divisor. Therefore we are in case~\textup{(1)}. It follows that $S$ is smooth, $\cE:=\varphi_*\cL$ is locally free of rank $r+1$,	$X\simeq\PP_S(\cE) $
and there is a short exact sequence
	\[
	0\longrightarrow
	\cO_S^{\oplus(r+1-\dim S)}
	\longrightarrow
	\cE
	\longrightarrow
	T_S
	\longrightarrow0.
	\]
In particular, $\dim S\leq r+1$. If $\dim S=r+1$, then $\dim X=2r+1$ and case~\textup{(1)} occurs.
Otherwise, we have  $1\leq\dim S\leq r$	and hence $r+1\leq\dim X\leq 2r$. This is precisely case~\textup{(4)}.
\end{proof}

 Clearly, case (1) occurs. Case (2) occurs by Example \ref{example-P-{r+s}}. Example \ref{example-semi-abelian-semi-contact} 
 shows that case (3) is realized. To see that case (4) also occurs let us prove the following lemma (generalizing example (3) from Proposition \ref{semi-contact-log-Fano}).
 
 \begin{Lemma}\label{example-semi-contact}
 For every $1\leq d\leq r$, the variety
 	$\PP_{\PP^d}\left(
 	\cO_{\PP^d}(1)^{\oplus(d+1)}
 	\oplus
 	\cO_{\PP^d}^{\oplus(r-d)}
 	\right)$
 	is  semi-contact of rank $(2r+1)$.
 \end{Lemma}

 \begin{proof}
 If $d=r$ then 	$\PP_{\PP^d}\left(\cO_{\PP^d}(1)^{\oplus(d+1)}\oplus
\cO_{\PP^d}^{\oplus(r-d)} \right)\simeq \PP^r\times\PP^r$. This variety is  semi-contact of rank $(2r+1)$ by Proposition \ref{semi-contact-log-Fano} that will be proven later.
 So in the following we assume that $1\leq d\leq r-1$. Let us set
$k:=r+1-d\geq2$ and $S:=\PP^d.	$
 	Choose line bundles  $\cM_1:=\cO_{\PP^d}(1)$,
 	$\cM_2=\cdots=\cM_k:=\cO_{\PP^d}$
 and consider the smooth projective variety
 	\[
 	q\colon
 	Y:=
 	\PP_S(\cO_S\oplus\cM_1)
 	\times_S\cdots\times_S
 	\PP_S(\cO_S\oplus\cM_k)
 	\longrightarrow S.
 	\]
 Note that	$\dim Y=\dim S+k=r+1$.

 For every $1\leq i\leq k$, let $S\longrightarrow\PP_S(\cO_S\oplus\cM_i)$
 be the section corresponding to the quotient $\cO_S\oplus\cM_i\longrightarrow\cM_i$,
 	and let $B_i\subset Y$ be its inverse image under the projection onto the
 	$i$-th factor. Let us set
 	\[
 	B:=\sum_{i=1}^kB_i
 	\qquad\text{and}\qquad
 	\Sigma:=\bigcap_{i=1}^kB_i.
 	\]
 	Then $(Y,B)$ is an snc pair and the restriction of $q$ to $\Sigma$ induces an
 	isomorphism $q|_{\Sigma}\colon\Sigma\xrightarrow{\sim}S$.
 	We use this isomorphism to identify $\Sigma$ with $S$.
 	
 	The restriction of the logarithmic tangent bundle to $S$ fits into the
 	logarithmic normal exact sequence
 	\[
 	0\longrightarrow\cO_S^{\oplus k}
 	\longrightarrow
 	T_Y(-\log B)|_S
 	\longrightarrow
 	T_S
 	\longrightarrow0.
 	\]
 	The class of this extension is
 	\[
 	\bigl(
 	c_1(\cM_1),\ldots,c_1(\cM_k)
 	\bigr)
 	=
 	\bigl(
 	c_1(\cO_{\PP^d}(1)),0,\ldots,0
 	\bigr)
 	\in
 	H^1(S,\Omega_S^1)^{\oplus k}.
 	\]
 	Therefore
 	\[
 	\cE:=T_Y(-\log B)|_S
 	\simeq
 	\cO_{\PP^d}(1)^{\oplus(d+1)}
 	\oplus
 	\cO_{\PP^d}^{\oplus(r-d)}.
 	\]
 	More precisely, the nontrivial component of the above extension is the
 	Euler sequence, whereas the remaining $r-d$ components are trivial. 
 	
 By Proposition \ref{standard-contact-structure} the snc pair $(\tilde X, \tilde D):=\bigl(\PP\bigl(T_Y(-\log B)\bigr) ,\tilde\pi^{-1}(B)\bigr)$
 is contact. Let $\tilde\pi\colon\tilde X\longrightarrow Y$ be the canonical projection. The inverse image of the stratum $S$ is
 	\[
 	X:=\tilde\pi^{-1}(S)
 	\simeq
 	\PP_S\bigl(T_Y(-\log B)|_S\bigr)
 	=
 	\PP_S(\cE)\simeq \PP_{\PP^d}\left(
 	\cO_{\PP^d}(1)^{\oplus(d+1)}
 	\oplus
 	\cO_{\PP^d}^{\oplus(r-d)}
 	\right).
 	\]
 	Since $B_1,\ldots,B_k$ are all the irreducible components of $B$, the
 	boundary induced on $X$ is empty. By restriction of the contact Lie
 	algebroid on $(\tilde X,\tilde D)$ to this stratum, $X$ is a smooth
 	projective semi-contact variety of rank $(2r+1)$.
 \end{proof}
 
 \section{Constructions of contact snc pairs}\label{examples-contact-log-pairs}
 
 In this section we give some examples of contact  snc pairs. It is well known that a Mori contraction of a projective contact complex manifold is of fiber type (see \cite[Lemma 2.10]{KPSW}). In fact, any contact smooth projective variety with $b_2\ge 2$ is the projectivization of the tangent bundle over a smooth projective  variety (see \cite[Theorem 2.12]{KPSW}). Similarly, the projectivization of the logarithmic tangent bundle over a projective snc pair with a natural simple normal crossing divisor carries a canonical contact structure. Below we describe this construction in more detail and we call it the standard contact log structure.  However, in the logarithmic case  there also exist more complicated contact structures that live on projectivizations of some special vector bundles on lower dimensional varieties.
 
 In this section we do not need to assume that varieties are projective.

 \subsection{Standard contact log structure}

 \begin{Proposition}\label{standard-contact-structure}
 	Let $r\geq 1$ and let $(Y,B)$ be an snc pair, where $\dim Y=r+1$. Let $\pi: X=\PP( T_Y(-\log B))\to Y$
 	and let $D=\pi^{-1}(B)$. Then the snc pair $(X,D)$ carries a canonical contact structure. Moreover, $\pi$ is a good elementary log contraction of $(X,D)$.
 \end{Proposition}
 
 \begin{proof}
 	Let us set $\cO_X(1):=\cO_{\PP( T_Y(-\log B))}(1)$. Let us consider the snc pair consisting of the total space of the logarithmic cotangent bundle 	
 	$\tilde \pi : \tilde X=\TT^*Y(-\log B)=\VV (T_Y(-\log B))\to Y$
 	with the divisor $\tilde D:=\tilde{\pi}^{-1}(B)$. Let
 	$
 	\tilde \eta\in
 	H^0\left(
 	\tilde X,\Omega_{\tilde X}^1(\log\tilde D)
 	\right)
 	$
 	be the canonical logarithmic Liouville form. Intrinsically, if
 	$\alpha\in\Omega_Y^1(\log B)_y$ is regarded as a point of
 	$\tilde X$, then
 	\[
 	\tilde \eta_\alpha(v)=\alpha(d\tilde\pi(v)).
 	\]
 	Then $\tilde \omega:=d\tilde \eta$ is a canonical logarithmic symplectic form on
 	$(\tilde X,\tilde D)$.
 	Indeed, the vertical tangent bundle $T_{\tilde X/Y}$ of $\tilde\pi$ is canonically isomorphic to
 	\(
 	\tilde\pi^*\Omega_Y^1(\log B).
 	\)
 	If $v_\beta$ is the vertical tangent vector corresponding to
 	$\beta\in\Omega_Y^1(\log B)_y$, and $w$ is a logarithmic tangent
 	vector whose image in $T_Y(-\log B)_y$ is $u$, then
 	\[
 	\tilde \omega(v_\beta,w)=\beta(u),
 	\]
 	up to the sign determined by the convention for the Liouville form.
 	Therefore the perfect pairing
 	\[
 	\Omega_Y^1(\log B)\otimes T_Y(-\log B)\longrightarrow\cO_Y
 	\]
 	shows that $\tilde \omega$ is nondegenerate.
 	
 	Let $\tilde X^\circ$ be the complement of the zero section in $\tilde X$ and let $p:\tilde X^\circ\to X$
 	denote the canonical projection. Note that $p$ is the principal $\mathbb G_m$-bundle associated with $\cO_X(-1)$. If $h_a$ denotes multiplication by $a\in\mathbb G_m$ in the fibers of $\tilde X\to Y$, then
 	\[
 	h_a^*\tilde\eta=a\tilde\eta
 	\qquad\text{and}\qquad
 	h_a^*\tilde\omega=a\tilde\omega.
 	\]
 	Thus $\tilde\omega|_{\tilde X^\circ}$ is $\mathbb G_m$-equivariant. So Corollary \ref{general-LeBrun-lemma-v3} gives a canonical $\cO_X(1)$-valued contact structure on $T_X(-\log D)$.	
 	The last assertion concerning $\pi$ follows from the fact that $\cO_X(-(K_X+D))\simeq\cO_X(r+1)$ is $\pi$-ample and $N^1(X/Y)$ is generated by $\cO_X(1)$.
 \end{proof}
 
 \medskip
 
 \begin{Remark}
 	One can easily see that the canonical contact structure $T_X(-\log D)\to\cO_X(1)$
 	is obtained as the composition of $d\pi:T_X(-\log D)\to\pi^*T_Y(-\log B)$ with the tautological quotient map $\pi^*T_Y(-\log B)\xrightarrow{\theta}\cO_X(1)$.
 \end{Remark}

 \subsection{Non-standard contact log structures}

 \begin{Proposition}\label{non-standard-contact-structure}
 	Let us fix positive integers $1\leq s\leq r$.
 	Let $(Y,B)$ be an snc pair, where $\dim Y=r-s+1$, and let $\cM_1,\ldots,\cM_s$ be line bundles on $Y$. Let
  $\cC$ be the canonical extension of $T_Y(-\log B)$ by $\cO_Y^{\oplus s}$ with extension class
 	\(
 	-\bigoplus_{i=1}^s c_1^{\log}(\cM_i).
 	\)
 	Let $\cE=\bigoplus_{i=1}^s\cM_i\oplus\cC$ and let $\pi:X=\PP(\cE)\to Y$. Let $D^h_j$ be the divisor corresponding to the projection
 	\(
 	\cE\to\bigoplus_{i\ne j}\cM_i\oplus\cC
 	\)
 	and let $D=\sum_{j=1}^sD^h_j+\pi^*B$. Then the snc pair $(X,D)$ carries a canonical contact structure. Moreover, $\pi$ is a good elementary log contraction of $(X,D)$.
 \end{Proposition}
 
 \begin{proof}
 	Let us set $\underline{\cM}:=\bigoplus_{i=1}^{s}\cM_i$ and consider the total space
 	$
 	q:Z:=\mathbb V_Y(\underline{\cM})\longrightarrow Y.
 	$
 	Thus the fibers of $Z$ are the duals of the fibers of $\underline{\cM}$. For every $i$, let $H_i\subset Z$ be the zero divisor of the $i$-th fiber coordinate, and put $A:=\sum_{i=1}^{s}H_i+q^{-1}(B)$.
 	Then $(Z,A)$ is an snc pair and
 	\(
 	\dim Z=\dim Y+s=r+1.
 	\)
 	
 	We first identify the logarithmic tangent bundle of $(Z,A)$. The logarithmic relative tangent sequence is
 	\[
 	0\longrightarrow\cO_Z^{\oplus s}
 	\longrightarrow T_Z(-\log A)
 	\longrightarrow q^*T_Y(-\log B)
 	\longrightarrow0,
 	\]
 	where the $i$-th standard basis vector in $\cO_Z^{\oplus s}$ corresponds to the Euler vector field $\epsilon_i$ along the $i$-th line-bundle factor. By definition, we have
 	\[
 	\cC=
 	\cD_{T_Y(-\log B)}^{\leq1}(\cM_1)
 	\times_{T_Y(-\log B)}\cdots\times_{T_Y(-\log B)}
 	\cD_{T_Y(-\log B)}^{\leq1}(\cM_s).
 	\]
 	Thus $\cC$ acts naturally on the total space $Z$: an element of $\cD_{T_Y(-\log B)}^{\leq1}(\cM_i)$ acts on $\cM_i^{-1}$ and hence induces a fiberwise linear vector field on $\mathbb V_Y(\cM_i)$, while the common symbol gives its projection to $T_Y(-\log B)$. These vector fields preserve the divisors $H_i$ and $q^{-1}(B)$. Therefore we obtain a canonical morphism
 	\( 	q^*\cC\to T_Z(-\log A),	\)
 	which fits into the commutative diagram
 	\[
 	\begin{CD}
 		0@>>>q^*\cO_Y^{\oplus s}@>>>q^*\cC
 		@>>>q^*T_Y(-\log B)@>>>0\\
 		@.@V{1_i\mapsto-\epsilon_i}V{\simeq}V@VVV@V{\mathrm{id}}VV@.\\
 		0@>>>\cO_Z^{\oplus s}@>>>T_Z(-\log A)
 		@>>>q^*T_Y(-\log B)@>>>0.
 	\end{CD}
 	\]
 	Indeed, the identity operator on $\cM_i^{-1}$ induces minus the Euler vector field $\epsilon_i$, because $\epsilon_i$ acts as multiplication by $-1$ on functions of weight $-1$. Hence we have a canonical isomorphism
 	\[
 	q^*\cC\simeq T_Z(-\log A).
 	\]
 	
 	Now let us consider \( \tilde X:=\mathbb V_Y(\cE).\) 
 	The direct-sum decomposition $\cE=\underline{\cM}\oplus\cC$ gives a projection $\rho:\tilde X\to Z$, and the preceding identification gives canonical isomorphisms
 	\[
 	\tilde X
 	\simeq\mathbb V_Z(q^*\cC)
 	\simeq\mathbb V_Z\bigl(T_Z(-\log A)\bigr),
 	\]
 	i.e., $\tilde X$ is canonically the total space of the logarithmic cotangent bundle of $(Z,A)$ (here we use the assumption that $\cE$ is the direct sum $\underline{\cM}\oplus\cC$ and not just an extension of $\cC$ by $\underline{\cM}$). Let
 	\(
 	\tilde D:=\rho^{-1}(A)
 	\)
 	and let
 	\(
 	\tilde\eta\in
 	H^0\left(
 	\tilde X,\Omega_{\tilde X}^1(\log\tilde D)
 	\right)
 	\)
 	be the canonical logarithmic Liouville form. As in the proof of Proposition \ref{standard-contact-structure}, its exterior derivative
 	\(
 	\tilde\omega:=d\tilde\eta
 	\)
 	is the canonical logarithmic symplectic form on the logarithmic cotangent bundle of $(Z,A)$.
 	
 	We claim that $\tilde\omega$ is homogeneous of weight $1$. Indeed, let $m_a:Z\to Z$ be scalar multiplication by $a\in\mathbb G_m$ in the fibers of $Z=\mathbb V_Y(\underline{\cM})$. Under the identification of $\tilde X$ with the logarithmic cotangent bundle of $Z$, the scalar action on $\tilde X=\mathbb V_Y(\cE)$ is the composition of the logarithmic cotangent lift of $m_a$ and multiplication by $a$ in the cotangent fibers. Cotangent lifts preserve the logarithmic Liouville form, whereas fiberwise multiplication by $a$ multiplies it by $a$. Hence
 	\[
 	a^*\tilde\eta=a\tilde\eta
 	\qquad\text{and}\qquad
 	a^*\tilde\omega=a\tilde\omega.
 	\]
 	
 	Let $\tilde X^\circ$ be the complement of the zero section in $\tilde X$ and let
 	\(
 	p:\tilde X^\circ\to X
 	\)
 	denote the canonical projection. Then $p$ is the principal $\mathbb G_m$-bundle associated with $\cO_X(-1)$ and
 	\(
 	\tilde D|_{\tilde X^\circ}=p^{-1}(D).
 	\)
 	Therefore Corollary \ref{general-LeBrun-lemma-v3} applied to the $\GG_m$-equivariant logarithmic symplectic form $\tilde\omega|_{\tilde X^\circ}$ induces a canonical $\cO_X(1)$-valued contact structure on $T_X(-\log D)$.
 	The last part comes from the fact that $\cO_X(-(K_X+D))\simeq\cO_X(r+1)$ is $\pi$-ample and $N^1(X/Y)$ is generated by $\cO_X(1)$.
 \end{proof}
 
 \medspace 

 The following two examples show that in Theorem \ref{main1} a non-split extension 
 \[
 0\longrightarrow \cM\longrightarrow \cE\longrightarrow \cC
 \longrightarrow0
 \]
 can define a logarithmic contact pair, but an arbitrary non-split extension need not lead to such a structure.
 The problem with generalization of construction from the proof of Proposition \ref{non-standard-contact-structure} is  that $\rho:\mathbb V_Y(\cE)\to Z$
 is in general only a torsor under $\mathbb V_Z(q^*\cC)
 \simeq \TT_Z^*(\log A).$
 So it has no canonical identification with the logarithmic cotangent bundle
 and therefore carries no canonical Liouville form. Local Liouville forms
 may fail to glue, and the obstruction to gluing is controlled by the
 extension cocycle. This explains why some non-split extensions still give
 rise to logarithmic contact structures, whereas others do not.

 \begin{Example}[\emph{A non-split extension defining a log contact snc pair}]\label{AI1}

 Let $Y$ be a complex elliptic curve, let $B=0$, and set $r=s=1$ and $\cM=\cO_Y$. Then
 	\(
 	\cC\simeq \cO_Y\oplus T_Y\simeq\cO_Y^{\oplus2}.
 	\)
 	Choose a nonzero class
 	\(
 	\xi\in H^1(Y,\cO_Y)
 	\) 	and let us represent it by a \v{C}ech 	cocycle $\{f_{ij}\}$.
 Let us consider 
  a (non-split) extension
 \begin{equation}\label{eq:contact-extension}
 	0\longrightarrow\cM\longrightarrow\cE
 	\longrightarrow\cC\longrightarrow0.
 \end{equation}	 
 defined by the class $\beta\in H^1(Y,\cC^*)  =H^1(Y,\cO_Y)\oplus H^1(Y,\Omega_Y^1)$
 of the cocycle $(f_{ij},df_{ij})$. The class of $(df_{ij})$ in $H^1(Y,\Omega_Y^1)$ vanishes by degeneration of the Hodge-to-de Rham spectral sequence and hence $\beta=(\xi ,0)$ but the choice of the cocycle helps to obtain convenient local fiber coordinates.
 
 Let $\tilde\pi:\tilde X:=\VV_Y(\cE)\longrightarrow Y.$ Let us choose an affine open covering $\{U_i\}$ so that on each $U_i$ we have a local coordinate $y_i$ and a local splitting of~\eqref{eq:contact-extension}. Let
 \( \frac{\partial}{\partial y_i} \) be the local generator of $T_Y|_{U_i}$. Let
 	\( (y_i;z_i,t_i,x_i) \)
 	be the local coordinates on \(\tilde\pi^{-1}(U_i)\), where \(z_i\), \(t_i\) are the fiber coordinates dual to the local generators of \(\cM\), \(\cO_Y\subset\cC\), 
 	 and $x_i$ corresponds to
 	$\frac{\partial}{\partial y_i} \in T_Y\subset\cC$.
 In particular, the divisor \( \tilde D:=\VV_Y(\cC)\subset\tilde X \) is locally given by $z_i=0$.
 Then  we have the overlap formulas
 	\begin{equation}\label{eq:contact-overlap}
 		z_i=z_j,\qquad
 		t_i=t_j+f_{ij}z_j,\qquad
 		x_i=\frac{\partial y_j}{\partial y_i}x_j+\frac{\partial f_{ij}}{\partial y_i}z_j.
 	\end{equation}
 	In particular,
 	\[
 	x_i\,dy_i=x_j\,dy_j+z_j\,df_{ij}.
 	\]
 	
 	On $\tilde\pi^{-1}(U_i)$ define
 	\begin{equation}\label{eq:contact-local-eta}
 		\eta_i:=t_i\frac{dz_i}{z_i}-dt_i+x_i\,dy_i.
 	\end{equation}
 	Using~\eqref{eq:contact-overlap}, we have
 	\[
 		\eta_i
 		=(t_j+f_{ij}z_j)\frac{dz_j}{z_j}
 		-d(t_j+f_{ij}z_j)+x_j\,dy_j+z_j\,df_{ij}
 		=t_j\frac{dz_j}{z_j}-dt_j+x_j\,dy_j
 		=\eta_j.
 	\]
 	Thus the local forms~\eqref{eq:contact-local-eta} glue to a global
 	logarithmic one-form
 	\(
 	\eta\in H^0\bigl(\tilde X,\Omega_{\tilde X}^1(\log\tilde D)\bigr).
 	\)
 	Let us set
 	\(\tilde\omega:=d\eta. \)
 	On $\tilde\pi^{-1}(U_i)$ we have
 	\[
 	\tilde\omega
 	=-\frac{1}{z_i}dz_i\wedge dt_i+dx_i\wedge dy_i
 	\]
 	and hence
 	\[
 	\tilde\omega^{\wedge2}
 	=-\frac{2}{z_i}
 	dz_i\wedge dt_i\wedge dx_i\wedge dy_i.
 	\]
 	Therefore $\tilde\omega$ is a logarithmic symplectic form on
 	$(\tilde X,\tilde D)$.
 	
 	Let $\tilde X^0$ be the complement of the zero section in $\tilde X$, and
 	let
 	\(
 	p:\tilde X^0\longrightarrow X:=\PP_Y(\cE)
 	\)
 	be the canonical projection. The scalar $\GG_m$-action is locally given by
 	\[
 	u\cdot(y_i;z_i,t_i,x_i)
 	=(y_i;uz_i,ut_i,ux_i).
 	\]
 So $\tilde\omega|_{\tilde X^0}$ is $\GG_m$-equivariant and hence by Corollary \ref{general-LeBrun-lemma-v3}
it gives rise to an $\cO_X(1)$-valued contact structure on $(X,D)$, where $ D:=\PP_Y(\cC)$ (cf. the proof of Proposition \ref{standard-contact-structure}). Since we use local coordinates, it is easy to write down the corresponding contact form explicitly.  Namely, the restriction of $\eta$ to $\tilde X^0$ has weight $1$ and it is horizontal	for the $\GG_m$-action. Hence it descends to an $\cO_X(1)$-valued logarithmic	one-form $\theta$ on $(X,D)$, where $ D:=\PP_Y(\cC)$.	Since \( \tilde\omega|_{\tilde X^0}=d(p^*\theta)\)
 	is nondegenerate, $\theta$ is a contact form on $(X,D)$.
 	Finally, note that $\pi:X\to Y$ is an elementary log contraction and its restriction to
 	$X\setminus D$ is not an isomorphism.
 \end{Example}
 
 \begin{Example}[\emph{A non-split extension which does not give rise to a log contact pair}]\label{AI2}
 	Let $Y$ be a smooth complex projective curve of genus $g\geq2$, let $B=0$,
 	and again set $r=s=1$ and $\cM=\cO_Y$. Then
 	\(
 	\cC\simeq\cO_Y\oplus T_Y.
 	\)
 	Choose a nonzero class
 	\(
 	0\neq\zeta\in
 	\operatorname{Ext}^1(T_Y,\cO_Y)
 	=H^1(Y,\Omega_Y^1),
 	\)
 	and let
 	\begin{equation}\label{eq:V-extension}
 		0\longrightarrow\cO_Y\longrightarrow\cV
 		\longrightarrow T_Y\longrightarrow0
 	\end{equation}
 	be the corresponding non-split extension. Set
 	\(
 	\cE:=\cO_Y\oplus\cV,
 	\)
 	and embed $\cM=\cO_Y$ as the subbundle $\cO_Y\subset\cV$. This gives a
 	non-split extension
 	\begin{equation}\label{eq:noncontact-extension}
 		0\longrightarrow\cM\longrightarrow\cE
 		\longrightarrow\cO_Y\oplus T_Y=\cC\longrightarrow0.
 	\end{equation}
 Set $X:=\PP_Y(\cE)$, $D:=\PP_Y(\cC)$, and let $\pi:X\to Y$ be the projection. We claim that the snc pair $(X,D)$
 does not carry any contact structure.
 	
Let \( \tilde\pi:\tilde X:=\VV_Y(\cE)\to Y, \) \( \tilde D:=\VV_Y(\cC), \)
and choose an affine open covering $\{U_i\}$ with local coordinates $y_i$
 	and local splittings of~\eqref{eq:V-extension}. Represent $\zeta$ by a
 	\v{C}ech cocycle
 	\(
 	\gamma_{ij}\,dy_i\in\Gamma(U_{ij},\Omega_Y^1).
 	\)
 	Let \(	(y_i;z_i,t_i,x_i) \)
 	be the corresponding coordinates on $\tilde\pi^{-1}(U_i)$, where $z_i$
 	corresponds to $\cM\subset\cV$, $t_i$ to the direct summand
 	$\cO_Y\subset\cE$, and $x_i$ to a local lift of
 	$\frac{\partial}{\partial y_i} $ to $\cV$. On overlaps we have
 	\begin{equation}\label{eq:noncontact-overlap}
 		z_i=z_j,\qquad
 		t_i=t_j,\qquad
 		x_i=\frac{\partial y_j}{\partial y_i}x_j+\gamma_{ij}z_j.
 	\end{equation}
 	The divisor $\tilde D$ is locally given by $z_i=0$.
 	
 	Since $\det\cE\simeq T_Y$ and
 	$\cO_X(D)\simeq\cO_X(1)$, we have	$\cO_X(-(K_X+D))\simeq\cO_X(2)$.
Assume  $(X,D)$ carries a contact structure $\theta: T_X(-\log D)\to \cL$. Then
 	$\cL^{\otimes2}\simeq\cO_X(2)$,
  	and consequently $\cL\simeq\cO_X(1)\otimes\pi^*\cN$
 	for some $\cN\in\operatorname{Pic}(Y)[2]$.
 	Using
 	\[
 	\Omega_{X/Y}^1(\log D)\simeq
 	\pi^*\cC\otimes\cO_X(-1),
 	\]
 	the relative logarithmic cotangent sequence gives
 	\begin{equation}\label{eq:relative-log-sequence}
 		0\longrightarrow
 		\pi^*(\Omega_Y^1\otimes\cN)\otimes\cO_X(1)
 		\longrightarrow
 		\Omega_X^1(\log D)\otimes\cL
 		\longrightarrow
 		\pi^*(\cC\otimes\cN)
 		\longrightarrow0.
 	\end{equation}
 	Therefore the contact form $\theta^*$ induces \( \bar\omega\in H^0(Y,\cC\otimes\cN).\)
 	
 	Suppose first that $\cN\neq\cO_Y$. Then  $H^0(Y,\cN)=0$ and $H^0(Y,T_Y\otimes\cN)=0$, so $\bar\omega=0$. It follows that the rank-two integrable bundle	$T_{X/Y}(-\log D)$ coincides with the kernel of the contact form $\theta$, a contradiction. 
 	
 	We can therefore assume that $\cN=\cO_Y$. Since $H^0(Y,T_Y)=0$, the relative
 	component $\bar\omega$ of a contact form is a nonzero multiple of
 	\((1,0)\in H^0(Y,\cO_Y\oplus T_Y).\)
 	After rescaling, we can assume that it is equal to $(1,0)$. Let
 	$p:\tilde X^0\longrightarrow X$
 	be the canonical $\GG_m$-torsor. On $\tilde\pi^{-1}(U_i)$, the pullback of
 	the contact form must then have the form
 	\begin{equation}\label{eq:general-local-contact-form}
 		\eta_i
 		=t_i\frac{dz_i}{z_i}-dt_i
 		+(a_i z_i+b_i t_i+c_i x_i)\,dy_i
 	\end{equation}
 	for some regular functions $a_i,b_i,c_i$ on $U_i$.
 	Taking the exterior derivative of~\eqref{eq:general-local-contact-form}, we
 	obtain
 	\[
 	d\eta_i
 	=-\frac{1}{z_i}dz_i\wedge dt_i
 	+(a_i\,dz_i+b_i\,dt_i+c_i\,dx_i)\wedge dy_i,
 	\]
 	because all terms containing the wedge product of two differentials from the
 	one-dimensional base vanish. Hence
 	\begin{equation}\label{eq:general-top-wedge}
 		(d\eta_i)^{\wedge2}
 		=-\frac{2c_i}{z_i}
 		dz_i\wedge dt_i\wedge dx_i\wedge dy_i.
 	\end{equation}
 	The contact condition therefore forces every $c_i$ to be invertible.
 	
 	The horizontal terms in~\eqref{eq:general-local-contact-form} glue to a
 	section
 	\[
 	h\in H^0(Y,\Omega_Y^1\otimes\cE)
 	=\operatorname{Hom}(T_Y,\cE).
 	\]
 	In the coordinates above,
 	\[
 	h\left(\frac{\partial}{\partial y_i}\right)
 	=a_i z_i+b_i t_i+c_i x_i.
 	\]
 	The functions $c_i$ are the local expressions of the endomorphism
 	\[
 	T_Y\xrightarrow{h}\cE
 	\longrightarrow\cC
 	\longrightarrow T_Y.
 	\]
 	By~\eqref{eq:general-top-wedge}, this endomorphism is an isomorphism.
 	Projecting $\cE=\cO_Y\oplus\cV$ onto $\cV$ therefore produces a splitting
 	of~\eqref{eq:V-extension}, contrary to $\zeta\neq0$. Thus $(X,D)$ admits no
 	logarithmic contact structure.
 \end{Example}

\section{Semi-contact snc pairs with two log contractions}

In this section we describe structure of semi-contact snc pairs with two log contractions. To do so we first need to classify semi-contact log Fano varieties.

\subsection{Semi-contact log Fano varieties}

Using Theorem \ref{structure-of-log-contractions-special-Lie-alg} one can classify semi-contact snc pairs with two log contractions. 
Let us first note the following generalization of LeBrun--Salamon's result \cite[Corollary 4.2]{LeBrun-Salamon}.

\begin{Proposition}\label{semi-contact-log-Fano}
	Let $(Z, D^Z)$ be a rank $(2r+1)$ semi-contact log Fano variety with $b_2(Z)\ge 2$. Then one of the following occurs:
	\begin{enumerate} 
		\item  $(Z, D^Z)\simeq (\PP (T_{\PP^{r+1}}), 0)$.
		\item For some $m\ge 1$, $Z \simeq\PP ( \cO_{\PP^{r}}(1)^{\oplus (r+1)}\oplus \cO_{\PP^{r}}(m))$ and $D^Z$ corresponds to the projection
		$$\cO_{\PP^{r}}(1)^{\oplus (r+1)}\oplus \cO_{\PP^{r}}(m)\to \cO_{\PP^{r}}(1)^{\oplus (r+1)}.$$
		\item $(Z, D^Z)\simeq (\PP^r\times \PP^r, 0)$.
	\end{enumerate}
	Moreover, each of the above pairs is a semi-contact log Fano variety  (with exactly two good elementary log contractions).
\end{Proposition}

\begin{proof}
	If $D^Z\ne 0$ then by Lemma \ref{existence-of-lines}, \cite[Theorem 4.3]{Fujita-Log-Fano} shows that  $(Z,D^Z)$ is as in (2). The main problem is to show that such pair $(Z,D^Z)$ has a contact structure.
	Let us set $S=\PP^r$ and $\cM:= \cO_{\PP^{r}}(m)$. Since $m\ne 0$, the Atiyah extension of $\cM$ is non-trivial and hence 	$\cA_{\cM}= \cO_{\PP^r} (-1)^{\oplus (r+1)}$. Let 
	$\cE= \cO_{\PP^{r}}^{\oplus (r+1)} (1)\oplus \cO_{\PP^{r}}(m)$  and let us consider $\pi: Z=\PP (\cE) \to S$ with the divisor $D^Z=\PP( \cO_{\PP^{r}}^{\oplus (r+1)} (1))$
	corresponding to the projection  $\cE \to \cO_{\PP^{r}}^{\oplus (r+1)} (1)$.
	Then $(Z,D^Z)$ has a contact structure by
	Proposition \ref{non-standard-contact-structure}. 
	
	Now assume that $D^Z=0$. If $\dim Z=2r+1$ then $Z\simeq \PP(T_{{\PP^{r+1}}})$ by \cite[Corollary 4.2]{LeBrun-Salamon}. If $\dim Z\le 2r $ then by \cite[Theorem B]{Wisniewski1990} we get $Z\simeq \PP^r\times \PP^r$. This variety is semi-contact, since it is the boundary in the contact snc pair $(\PP^{r}\times {\PP^{r+1}}, \PP^r\times \PP^r)$ (or in any variety from case (2)).
	
Note that in all the cases we have $b_2(Z)= 2$. So, since $(Z, D^Z)$ is log Fano, it has exactly two good elementary log contractions.
\end{proof}

\medskip

\begin{Example}\label{example-log-Fano}
	It is useful to explicitly describe contractions of contact snc pairs in case (3). $(Z,D^Z)$ admits exactly two good elementary log contractions, one of which is given by the projectivization $\varphi_1$, so it is of fiber  type.  
	
	If $m=1$ then $Z=\PP^{r}\times {\PP^{r+1}}$ and $D^Z=\PP^r\times \PP^r$ is a divisor of type $(0,1)$. The second log contraction is the projection onto the second factor and it realizes $Z$ as $\PP (T_{\PP^{r+1}} (-\log H))$ for the hyperplane $H\subset \PP^{r+1}$, which is the image of $D^Z$.
	
	If $m>1$ then  the second log contraction is the blow up of the cone  in $\PP^{\binom{m+r}{m}+r}$
	with vertex $\PP^{r}$ over the $m$-th Veronese embedding $\PP^{r}\hookrightarrow \PP^{{\binom{m+r}{m}}-1}$.
	This log contraction is divisorial and it contracts $D^Z=\PP^r\times \PP^r$ to the vertex of the cone  (the restriction to $D^Z$ is projection onto the second factor of the product). Note that in this case  $Z$ is not isomorphic to $\PP( T_Y(-\log B))$ for any $(Y,B)$.
\end{Example}

\subsection{Semi-contact snc pairs with two contractions}

The following theorem gives a generalization of \cite[Proposition 2.13]{KPSW} to the semi-contact case.

\begin{Theorem}\label{two-contractions}
	Let $(Z,D^Z)$ be a rank $(2r+1)$ semi-contact projective snc pair that admits two	 good elementary log contractions
	$\varphi_1: Z\to S_1$  and $\varphi_2: Z\to S_2$. Assume that $\varphi _1$ is very good and $(Z, D^Z)$ is not log Fano.
	Then there exist
	\begin{enumerate}
		\item an $(r+1)$-dimensional projective variety $S$ smooth outside a closed point $s_0\in S$ at which
		 $\hat \cO_{S, s_0}$ is isomorphic to the completion of the cone over the $m$-th Veronese embedding of $\PP^r$ at the vertex for some $m\ge 1$, and
		\item a projective simple normal crossing divisor $B$ on $S\backslash \{s_0\}$,	
	\end{enumerate}
	such that 
	if $\psi_1: S_1\to S$ is the blow up of $S$ at $s_0$ with exceptional divisor $E$ then
	$$(\varphi_1: Z\to S_1)\simeq (\PP (T_{S_1}(-\log (\psi_1^{-1}(B)+E)))\to S_1)$$
	and $D^Z=\varphi _1^{-1}(\psi_1^{-1}(B)+E).$
	Moreover, each $(S, B)$ as above gives a semi-contact projective snc pair $(Z, D^Z)$ with two good elementary log contractions.
\end{Theorem}

\begin{proof}
	Let $\varphi_3: Z\to S$ be the good log contraction of $(Z, D^Z)$ contracting the extremal face spanned by the rays corresponding to $\varphi_1$ and  $\varphi_2$. There exist morphisms $\psi_1: S_1\to S$ and $\psi_2: S_2\to S$ such that  $\varphi _3=\psi _1\circ \varphi_1=\psi _2\circ \varphi_2$.
	Theorem \ref{structure-of-log-contractions-special-Lie-alg} describes the structure of $\varphi _1$ and it implies that $\varphi _1$ is a projective bundle with fibers of dimension $\ge r$. 
	If $\dim S=0$ then $(Z, D^Z)$ is log Fano, contradicting our assumption. So we can assume that $\dim S>0$. Then Proposition \ref{general-fiber-in-general} implies that a general fiber of $\varphi_3$ is a projective space. So a general fiber of $\psi _1$ is trivial, i.e., $\psi_1$ is birational.
	If $\varphi _2$ is of fiber type then the same argument shows that $\psi_2$ is birational. But then 
	general fibers of $\varphi _1$ and $\varphi_2$ coincide, which contradicts our assumption that they contract different extremal rays. So $\varphi_2$ is birational. If $\Exc (\varphi _2)$ is not contained in $\Supp D^Z$ then $\varphi _2$ is very good. Then Theorem \ref{structure-of-log-contractions-special-Lie-alg} implies that 
	 $\varphi _2$ is of fiber type, a contradiction. Therefore $\Exc (\varphi _2)\subset \Supp D^Z$.  For any stratum $Z'$ of the stratification induced by $D^Z$, the divisor $-(K_{Z'}+D^{Z'})=-(K_Z+D^Z)|_{Z'}$ is $\varphi _2|_{Z'}$-ample. Let us choose a stratum $Z'$ for which $\varphi _2|_{Z'}$ is not an isomorphism on $U_{Z'}$. Then  $(Z',D^{Z'})$ is semi-contact with two different good log contractions $\varphi _1|_{Z'}:Z'\to S_1'=\varphi _1 (Z')$ and  $\varphi _2|_{Z'}:Z'\to S_2'=\varphi _2 (Z')$. In fact, both of them are projective bundles by 
	Theorem \ref{structure-of-log-contractions-special-Lie-alg}. Let us consider the log contraction
	$\varphi _3|_{Z'}: Z'\to S'=\varphi _3 (Z')$. Since there are no curves contracted by both $\varphi_1|_{Z'}$ and $\varphi_2|_{Z'}$ and the fibers of these morphisms have dimension $\ge r$, we have $\dim Z'=2r$ and $\dim S'=0$. But then we are in the previous case. So $(Z', D^{Z'})\simeq (\PP^r\times \PP^r, 0)$, and  $\varphi_1|_{Z'}$ and $\varphi_2|_{Z'}$ correspond to projections onto both factors of the product. It follows that $\varphi_2$ is divisorial and hence by \cite[Proposition 5.1.6]{Kawamata-Matsuda-Matsuki} we have $Z'=\Exc (\varphi _2)$.	
	
	If $\dim S_1=r$ then $S_1=S_1'\simeq \PP^r$ is contracted by $\psi _1$ to a point as $\varphi _3$ contracts $Z'$ to a point. This contradicts our assumption that $\dim S>0$.
	So $\dim S_1=r+1$ and by  Theorem \ref{structure-of-log-contractions-special-Lie-alg} there exists a simple normal crossing divisor $B_1$ on $S_1$ such that $Z\simeq \PP (T_{S_1} (-\log B_1))$ and $D^Z=\varphi _1^{-1}(B_1)$. Moreover, since $D^{Z'}=0$, $S_1'$ is a connected component of $B_1$. 
	In the following we write $E$ for $S_1'$.
	By construction $\psi _1(E)$ is a point
	that we denote by $s_0$. Since
	$\dim N_1(S_1/S)=1$ (this follows from \cite[Lemma 3.2.5]{Kawamata-Matsuda-Matsuki}) and $-K_{E}=-(K_{S_1}+B_1)|_{E}$ is ample, the divisor $-(K_{S_1}+B_1)$ is $\psi_1$-ample. This shows that $\psi_1$ is a  good elementary log contraction of $(S_1,B_1)$. Therefore by \cite[Proposition 5.1.6]{Kawamata-Matsuda-Matsuki}, $E$ is equal to the exceptional locus of  $\psi_1$. We also know that
	$$T_{E} (-\log E)|_{E}=T_{S_1} (-\log B_1)|_{E}\simeq \cO _{\PP^r} (1)^{\oplus (r+1)}.$$
	In particular, the short exact sequence
	$$0\to \cO_{E}\to T_{S_1} (-\log E)|_{E}\to T_{E}\to 0$$
	does not split. The class of this extension in $\Ext ^1 _{E} (T_{E}, \cO_{E})= H^1(\Omega_{E}^1)$
	is equal to $c_1(\cO_{E}({E}))$, so the normal bundle $\cO_{E}({E})$ is non-trivial. 
	Since a subvariety with ample normal bundle can not be contracted to a point by a non-constant morphism,
	$\cO_{E}({E})$  is not ample and hence $\cO_{E}(-E)\simeq \cO _{\PP^r} (m)$ for some $m>0$.
	Now  Mori's version of the Grauert--Hironaka--Rossi theorem (see \cite[Lemma 3.33]{Mori1982}) implies that
	the formal completion of $S_1$ along $E$ is isomorphic to the formal completion of the total space of  $\cO_{\PP^r}(-m)$ along the zero section. Moreover, the contraction $\psi_1$ in a neighbourhood of $s_0\in S$  looks formally like contraction of the zero section. Thus $\hat \cO_{S, s_0}$
	is isomorphic to the completion of the cone over the $m$-th Veronese embedding of $\PP^r$ at the vertex.
	This shows that $(Z, D^Z)$ is of the described form.

	Now let us assume that $(Z,D^Z)$ is an snc pair coming from some $(S,B)$ as described in  the theorem.
	The fact that  $(Z,D^Z)$ is semi-contact follows from Proposition \ref{standard-contact-structure}. Since $c_1(\cO _E(E))\ne 0$, the sequence 
	$$0\to \cO_{E}\to T_{S_1} (-\log (\psi_1^{-1}(B)+E))|_{E}\to T_{E}\to 0$$
	is non-split and hence $T_{S_1} (-\log (\psi_1^{-1}(B)+E))|_{E}\simeq \cO_{\PP^{r}}(1)^{\oplus (r+1)}$. Therefore
	$Z'=\varphi_1^{-1}(E)\simeq \PP^r\times \PP^r$.
	Since $C=\PP^1 \times \{x_0\}$ is not contracted by $\varphi_1$ and $(K_Z+D^Z)\cdot C= K_{Z'}\cdot C=-(r+1)<0$,
	the pair $(Z, D^Z)$ admits another good log contraction.
\end{proof} 

\medskip

\begin{Remark}
	If $E$ is a divisor on a smooth projective variety and  $(E, \cO _E(-E))\simeq (\PP^r, \cO_{\PP^r}(m))$ for some $m>0$ then by Artin-Grauert--Fujiki's theorem $E$ can be contracted to a point. However, the contraction exists only in the category of either algebraic or analytic spaces. This is why to describe $(Z, D^Z)$ in the above theorem we start with a projective log pair $(S, B)$ and not with the snc pair $(S_1, B_1)$.
\end{Remark}

\begin{Remark}
	Note that in the above theorem $(S^{\an},s_0)$ is analytically isomorphic to the germ of the cone over the $m$-th Veronese embedding of $\PP^r$ at the vertex. This immediately follows from Grauert's criterion \cite{Grauert1962}, whose assumptions are trivially satisfied. The only reason for the above formulation of Theorem \ref{two-contractions} is avoiding analytic methods when dealing with algebraic varieties.
\end{Remark}

\begin{Example}\label{blow-up-works}
	Here we analyse the map $\psi_2: S_2\to S$ induced by the second contraction in Theorem \ref{two-contractions}. 
	By Theorem \ref{structure-of-log-contractions-special-Lie-alg} the restriction of the second good elementary log contraction $\varphi_2: Z\to S_2$ to $Z'=\Exc (\varphi _2)$ is isomorphic to the projection of $\PP^r\times \PP^r$ onto the second factor. Therefore all reduced  fibers of $\psi_2$ are isomorphic to $\PP^r$. Note also that over $S\backslash \{s_0\}$, $\psi_2$ is isomorphic to $\PP(T_S (-\log B))$. Let $\cL _2$ be the unique reflexive extension of the line bundle $\cO_{\PP(T_S (-\log B) |_{S\backslash \{s_0\}})}(1)$. By \cite[Corollary 1.7]{Hartshorne-Stable-reflexive-sheaves}, the sheaf $\psi_{2*}\cL_2$ is reflexive.
	So $\psi_{2*}\cL_2\simeq T_S (-\log B)$ as both sheaves are reflexive and isomorphic over $S\backslash \{s_0\}$. Since $S$ is klt, the sheaf $T_S (-\log B)$ (isomorphic to $T_S$ at $s_0$) is not locally free at $s_0$ unless $S$ is smooth at $s_0$ (see \cite[Theorem 6.1]{GKKP-Extension-of-diff-forms}), which is equivalent to $m=1$. So  $\psi_2$ is a projective bundle if and only if $m=1$. For $m>1$ this shows that Fujita's theorem \cite[Lemma 2.12]{Fujita1985} fails for even mildly singular varieties (note that $S_2$ is klt as the pair $(S_2, \varphi_{2*}(D^Z))$ is dlt). The problem here is that the fiber of $\psi_2$ over $s_0$ is non-reduced.
	
	For $m=1$ we have $(S_2\to S)\simeq (\PP (T_S(-\log B))\to S)$ and this example shows that $\varphi_2$ is the blow up of the fiber  of this projectivization over $s_0$.  
\end{Example}

\medskip

Proposition \ref{semi-contact-log-Fano} and Theorem \ref{two-contractions} imply the following corollary:

\begin{Corollary}\label{two-contractions-dim}
	Let $(Z,D^Z)$ be a semi-contact projective snc pair that admits two	good elementary log contractions
	$\varphi_1$  and $\varphi_2$. Assume that $\Exc(\varphi _1)\not \subset \Supp D^Z$.
	Then either $(Z, D^Z)$ is contact or  $(Z, D^Z)\simeq (\PP^r\times \PP^r, 0)$.
\end{Corollary}

\begin{proof}
	Assume first that $(Z,D^{Z})$ is log Fano. Since $(Z,D^{Z})$ admits two distinct elementary log contractions,
we have $b_{2}(Z)\geq 2$. So we can apply Proposition \ref{semi-contact-log-Fano}. In the first two cases the pair $(Z, D^Z)$ is contact and in the third one we have  $(Z, D^Z)\simeq (\PP^r\times \PP^r, 0)$.

Now let us assume that $(Z,D^{Z})$ is not log Fano. Then by Theorem \ref{two-contractions} 
there exist a projective variety $T_1$ and a simple normal crossing divisor $B_1$ (equal to $\psi_{1}^{-1}(B)+E$
from the theorem)  such that 
	\[
	(\varphi_{1}\colon Z\longrightarrow T_{1})
	\simeq
	\left(
	\mathbb{P}\bigl(
	T_{T_{1}}(-\log B_1
	\bigr)
	\longrightarrow T_{1}
	\right)
	\]
	and
	$D^{Z}=	\varphi_{1}^{-1}(B_1).$
	Therefore $(Z,D^{Z})$ is a standard contact snc pair of the form considered
	in Proposition \ref{standard-contact-structure}.
\end{proof}

\appendix

\section{Appendix: Semi-contact snc pairs in small dimensions}

In this section we try to classify low dimensional semi-contact snc pairs.
In the whole section $(X, D)$ is a semi-contact projective snc pair
with a special contact Lie algebroid $(\LL, \cL)$  of rank $(2r+1)$.
Let us first note the following simple corollary of our results:

\begin{Lemma}\label{semi-contact-dim-r}
If $\dim X=r$ and $(K_X+D)$ is not nef then $(X, D)\simeq (\PP^r, 0)$. Moreover, $\LL\simeq \cO_{\PP^{r}}^{\oplus r}\oplus \cO_{\PP^{r}}(1)^{\oplus (r+1)}$ or $\LL\simeq \cO_{\PP^r}^{\oplus (r+1)}\oplus T_{\PP^r}$. 
\end{Lemma}

\begin{proof}
	If $b_2(X)=1$ then the assertion follows from Proposition \ref{dimension-contact-Lie-algebroid}. So we can assume that  $b_2(X)>1$. Then $D=0$ again by  Proposition \ref{dimension-contact-Lie-algebroid}. 
Since $K_X$ is not nef, the cone theorem provides an elementary  $K_X$-negative contraction
$\varphi: X\to Y$ with $\dim Y>0$.	Then Theorem \ref{structure-of-log-contractions-special-Lie-alg} implies that 
$X\simeq \PP_Y(\cE)$ for a rank $r+1$ vector bundle $\cE$ on $Y$. Hence $\dim X=\dim Y+r>r$, a contradiction.
\end{proof}

\medskip 

\subsection{Dimension $1$}

In case $\dim X=1$ the following lemma gives a full classification of semi-contact snc pairs.

\begin{Lemma}\label{semi-contact-dim-1}
If $\dim X=1$ then $r=1$ and one of the following occurs:
	\begin{enumerate}
		\item $X\simeq \PP^1$ and $D$ consists of $2$ points,
		\item $X\simeq \PP^1$ and $D=0$,
		\item $X$ is an elliptic curve and $D=0$.
	\end{enumerate}
\end{Lemma}

\begin{proof}
Equality $r=1$ follows from Lemma \ref{rank-estimate}. Lemma \ref{special-contact-Kodaira} implies that $h^0(\cL)>0$ and hence $\deg (K_X+D)\le 0$. Since  $-(K_X+D)=(r+1) c_1(\cL)$, we see that $(X,D)$ is as above.

If $X\simeq \PP^1$ and $D=0$ then we are in ituation of Lemma \ref{semi-contact-dim-r}.

If  $X \simeq \PP^1$ and $D\ne 0$ then $\cL\simeq \cO_{X}$ and hence $\LL \simeq \cO_{X}^{\oplus 3}$.
Similarly, if $X$ is an elliptic curve then  $\cL\simeq \cO_{X}$.
In both cases existence a semi-contact structure on  $\LL =\cO_{X}^{\oplus 3}$ follows from Example \ref{example-semi-abelian-semi-contact}.
\end{proof}

\subsection{Dimension $2$}

Assume that $\dim X=2$. Then by Lemma \ref{rank-estimate} we have either $r=1$ or $r=2$.

If  $K_X+D$ is pseudoeffective then $K_X+D\sim 0$. By Example \ref{example-semi-abelian-semi-contact}  such semi-contact structures  really occur both for $r=1$ and $r=2$ on compactifications of semi-abelian surfaces (as $X$ one can even take an abelian surface with $D=0$).

If $K_X+D$ is not  pseudoeffective and $r=2$ then we can use Lemma \ref{semi-contact-dim-r}. So in the following we can assume that $K_X+D$ is not  pseudoeffective and $r=1$.

If $b_2(X)=1$ then Proposition \ref{dimension-contact-Lie-algebroid} implies that $X\simeq \PP^2$. $D$ is a hyperplane and  $\LL\simeq \cO_{\PP^{2}}\oplus \cO_{\PP^{2}}(1)^{\oplus 2}$.

If  $b_2(X)>1$ then $(X,D)$ admits a good elementary log contraction $\varphi: X\to Y$ with $\dim Y>0$. 

If $\varphi$ is of fiber type then  by Theorem \ref{structure-of-log-contractions-special-Lie-alg}, $r=1$, $Y$ is a  smooth projective curve and there exists a reduced divisor $B$ on $Y$ such that $X\simeq \PP(\cE)$, $D=\varphi^{-1} (B)$ and
we have
$$0\to \cO_Y\to \cE \to T_{Y}(-B)\to 0.$$
 
If $\varphi$ is divisorial then $Z=\Exc (\varphi)\subset \Supp D$. So $Z$ is an irreducible component of $D$ and $(Z, D^Z)$  is semi-contact. But $-(K_X+D)$ is $\varphi$-ample, so $-(K_Z+D^Z)=-(K_X+D)|_Z$ is ample. Hence $(Z, D^Z)\simeq (\PP^1, 0)$ by Lemma \ref{semi-contact-dim-1}. This shows that the image $D^Y$ of the divisor $(D-Z)$ is disjoint from the point $y_0=\varphi(Z)$. Note that if we write
$$K_X+D=\varphi^* (K_Y+D^Y)+\alpha Z$$
then $\alpha=-2/Z^2>0$. So $K_Y+D^Y$ is not pseudoeffective. Unfortunately, $Y$ can be singular at $y_0$ so we cannot immediately run an inductive procedure.

\subsection{Dimension $3$}

Assume that $\dim X=3$. Then by Lemma \ref{rank-estimate} we have $1\le r\le 3$.

If $K_X+D$ is pseudoeffective then by Corollary \ref{semi-contact-Kodaira} we have $r>1$ and
$K_X+D\sim 0$. By Example \ref{example-semi-abelian-semi-contact}
such semi-contact structures occur  for both $r=2$ and $r=3$ on compactifications of semi-abelian $3$-folds.

\begin{Lemma}
If $K_X+D$ is not pseudoeffective and $b_2(X)=1$ then $X\simeq \PP^3$ and one of the following cases occurs:
\begin{enumerate}
	\item $r=1$, $D=0$ and  $\LL=T_{\PP^3}$,
	\item $r=1$, $D$ is a sum of $2$  hyperplanes, and $\LL=T_{\PP^3}(-\log D)\simeq \cO_{\PP^{3}}\oplus \cO_{\PP^{3}}(1)^{\oplus 2} $,
	\item $r=2$,  $D$ is a hyperplane and  $\LL\simeq \cO_{\PP^{3}}^{\oplus 2}\oplus \cO_{\PP^{3}}(1)^{\oplus 3}$,
	\item $r=3$,  $D=0$ and  $\LL\simeq \cO_{\PP^{3}}^{\oplus 3}\oplus \cO_{\PP^{3}}(1)^{\oplus 4}$,
	\item $r=3$,  $D=0$ and  $\LL\simeq \cO_{\PP^{3}}^{\oplus 4}\oplus T_{\PP^{3}}$.
\end{enumerate}
\end{Lemma}

\begin{proof}
If $r\ge 2$ the assertion follows from Proposition \ref{dimension-contact-Lie-algebroid}.
So we can assume that $r=1$. Then $(X, D)$ is a contact log Fano variety with $-(K_X+D)=2c_1(\cL)$. 
If $D=0$ then it is well-known that $X\simeq \PP^3$ (see \cite[2.1]{Buczynski-Kapustka-Kapustka} for the history of this assertion). So in the following we assume that $D\ne 0$.
Since $D=a c_1(\cL)$ for some $a>0$, the Kobayashi--Ochiai  
theorem (see \cite[Chapter V, Theorem 1.11]{Kollar-Rational}) implies that $X\simeq \PP^3$ or $X$ is isomorphic to a $3$-dimensional quadric $Q^3$ in $\PP^4$. In the first case $a=2$ and $D$ is either a sum of two hyperplanes or $D$ is a smooth quadric. In the first case $(X,D)$ is contact.

If $D$ is a smooth quadric then $\cL \simeq \cO_{\PP^3} (1)$ and existence of a contact structure implies existence of a surjective map $T_{X}(-\log D)\to \cL$ with kernel $\cF$.
In particular, since  $\cF \otimes \cL^{-1}$ is a rank $2$ vector bundle, we have
$$\int_{X}c_3((T_{X}(-\log D))\otimes \cL ^{-1})=\int_Xc_3 (\cF \otimes \cL^{-1})=0.$$
But a short computation using the short exact sequence
$$0\to T_{\PP^3}(-\log D) \to T_{\PP^3}\to \cO_D(D)\to 0$$
shows that $\int_{\PP^3}c_3((T_{\PP^3}(-\log D))\otimes  \cL^{-1})=-1\ne 0$, a contradiction.

If $X\simeq Q^3$ then $a=1$ and hence $D$ is isomorphic to a smooth $2$-dimensional quadric   in $\PP^3$. The proof is this case is similar as above. Namely, if $\cF$ is the kernel of $T_{X}(-\log D)\to \cL$ then
$$\int_{X}c_3((T_{X}(-\log D))\otimes \cL ^{-1})=\int_Xc_3 (\cF \otimes \cL^{-1})=0.$$
Now short exact sequences
$$0\to T_{Q^3}(-\log D) \to T_{Q^3}\to \cO_D(D)\to 0$$
and 
$$0\to T_{Q^3}\to T_{\PP^4}|_{Q^3}\to \cO_{Q^3} (2)\to 0 $$
allow us to compute $\int_{Q^3}c_3((T_{Q^3}(-\log D))\otimes  \cO_{Q^3}(-1))=-2$, a contradiction.
\end{proof}

Now let us assume that $b_2 (X)>1$ and $K_X+D$ is not pseudoeffective. Then the pair $(X,D)$ admits a good elementary log contraction $\varphi: X\to Y$ with $\dim Y>0$. Note that in this case $r=1$ or $r=2$ by Lemma \ref{semi-contact-dim-r}.

If $\varphi$ is of fiber type then by Theorem \ref{structure-of-log-contractions-special-Lie-alg}, $X\simeq \PP(\cE)$ for some vector bundle $\cE$ and there exists a simple normal crossing divisor $B$ on $Y$ 
such that  $D^v=\varphi^{-1}(B)$ and one of the following happens:
\begin{enumerate}
	\item $r=1$, $\dim Y=1$ and there exists a line bundle $\cM$ on $Y$ such that
	$$0\to \cM\to \cE\to \cD ^{\le 1}_{T_{Y}(-\log B)}(\cM)\to 0$$
	and $D^h\simeq \PP (\cD ^{\le 1}_{T_{Y}(-\log B)}(\cM))$,
	\item  $r=1$, $\dim Y=2$, $D^h=0$ and $\cE\simeq T_{Y}(-\log B)$,
	\item $r=2$, $\dim Y=1$, $D^h=0$ and
	 $$0\to \cO_Y^{\oplus 2}\to \cE \to T_{Y}(-\log B)\to 0.$$
\end{enumerate}

If $\varphi$ is of divisorial type then the exceptional divisor $Z=\Exc (\varphi)$ is contained in $D$, so it is an irreducible component of $D$. It follows that $(Z, D^Z)$ is semi-contact and $\varphi|_Z$ is a good log contraction of fiber type. If  $\dim \varphi(Z)=0$ then by classification in the $2$-dimensional case and by Proposition \ref{semi-contact-log-Fano} we have the following possibilities:
\begin{enumerate}
	\item $r=2$, $Z\simeq \PP^2$ and $D^Z=0$,
	\item $r=1$,  $Z\simeq \PP^2$ and $D^Z$ is a line (so $D$ has another irreducible component),
	\item $r=1$, $Z\simeq \PP^1\times \PP^1$ and $D^Z=0$.
\end{enumerate}
If $\dim \varphi(Z)=1$ then $\varphi|_Z$ is elementary (e.g., by Proposition \ref{semi-contact-log-Fano} and Theorem \ref{two-contractions}). So
as in the $2$-dimensional case $r=1$, $T=\varphi(Z)$ is a smooth projective curve and there exists a reduced divisor $B$ on $T$ such that $Z\simeq \PP(\cE)$ (over $T$), $D^Z=(\varphi|_Z)^{-1} (B)$ and
we have
$$0\to \cO_T\to \cE \to T_{T}(-B)\to 0.$$

If $\varphi$ is small then $\Exc (\varphi)\subset \Supp D$. Let $Z$ be an irreducible component of $D$,
which contains a curve contracted by $\varphi$. The pair $(Z, D^Z)$ is semi-contact and $\varphi|_Z$ is a good
divisorial log contraction of this pair (possibly non-elementary).  
Then $\varphi|_Z$ factors through an elementary divisorial contraction of $(Z, D^Z)$. So from the $2$-dimensional case this we know that  $\Exc (\varphi|_Z)$ is isomorphic to $\PP^1$ and it does not intersect any other component of $D$. Hence  $\Exc (\varphi)$ is a disjoint union of numerically equivalent $\PP^1$'s.
Each connected component $C\simeq \PP^1$ of $\Exc (\varphi)$ is an intersection of two irreducible components $Z$ and $Z'$ of $D$.
Note that the conormal bundle $N^*_{C/X}$ is isomorphic to $\cO_C(-Z)\oplus \cO_C(-Z')\simeq \cO_{\PP^1}(-a)\oplus \cO_{\PP^1}(-b)$, where $a=-(C.C)_Z>0$ and $b=-(C.C)_{Z'}>0$.

\section{Appendix: Compactification of a $1$-jet bundle to a b-contact manifold}

In this appendix we work in the category of smooth real manifolds and $\cC^{\infty}$-functions. 
With suitable modifications, one can replace this with real analytic or complex analytic (holomorphic) objects. Before explaining analogue of our construction in the real case, let us recall some basic notions of b-geometry. Let $(M,Z)$ be a b-manifold, i.e., a smooth real manifold $M$  with the boundary hypersurface $Z$, called the {critical set}
(see \cite[2.1]{Melrose-book} for a precise definition). Then one defines a  b-tangent bundle $\leftindex^{b}{TM}$ (see \cite[Lemma 2.5]{Melrose-book}) and its dual, a b-cotangent bundle $\leftindex^{b}{T^*M}=(\leftindex^{b}{TM})^*$.

Let $(M, Z)$ be a b-manifold of dimension $(2r+1)$. In \cite{Miranda-Oms} a \emph{b-contact structure} on $(M,Z)$ is defined as a corank $1$ subbundle $ \xi$ of the b-tangent bundle (such $\xi$ is called a \emph{corank $1$ b-distribution}) given by the kernel of a globally defined b-$1$-form $\alpha\in \leftindex^{b}\Omega^1 (M)$ such that $\alpha\wedge (d\alpha)^{\wedge r}$
is nowhere vanishing as a section of $\bigwedge ^{2r+1}(\leftindex^{b}{T^*M})$. Note that in this definition  $\xi$ is necessarily coorientable, i.e., the normal bundle $\leftindex^{b}{T^*M}/\xi$ is trivial. But from our point of view 
the following definition is more natural:

\begin{Definition}
	A \emph{b-contact structure} on the b-manifold $(M,Z)$ is a corank $1$ b-distribution $\xi \subset \leftindex^{b}{TM}$ such that for every point $p\in M$ there exists an open neighbourhod  of $p$ on which $\xi$  is given by 	the kernel of a local  b-$1$-form $\alpha$  such that $\alpha\wedge (d\alpha)^{\wedge r}$ does not vanish at $p$.
\end{Definition}

It is well known that the bundle of contact elements of $M$ is the projectivized cotangent bundle
and it carries a natural contact structure (see, e.g., \cite[Appendix 4 D]{Arnold-Mathematical-methods} or \cite[1.2]{Geiges-Contact-topology}).  Note that here we use traditional in real geometry convention that the projectivized bundle is defined by lines in the fibers of the bundle and not by hyperplanes as in the main body of the paper (in this case we talked about projectivizations rather than projectivized bundles).
There is also an analogous construction of the manifold of oriented contact elements that in the case of a Riemannian manifold is isomorphic to the unit sphere bundle in $T^*M$ (see, e.g.,  \cite[Example 3.45]{McDuff-Salamon}). Both these structures come from the homogeneous symplectic structure on the cotangent bundle $T^*M$ (see \cite[A.2]{Kashiwara-Schapira} and \cite[Appendix 4 E]{Arnold-Mathematical-methods}).
Let us recall that the canonical contact structure on the projectivised cotangent bundle does not need to be coorientable
(see, e.g., \cite[Examples 1.2.4 and 2.1.11]{Geiges-Contact-topology}). This is the main reason why we extend the definition of a b-contact structure.

For any b-manifold $(M,Z)$ there is also an analogous construction of a canonical contact structure on the projectivization of the b-cotangent bundle  $\leftindex^{b}{T^*M}$. An  analogue of the construction for the unit sphere bundle in $\leftindex^{b}T^*M$ (presumably after fixing an auxiliary Riemannian metric) is described in \cite[Section 6, Example 2]{Bradell-et-al} and \cite[Example 8.4]{Miranda-Oms}.

Here we describe a different b-contact manifold, which is a fibration with compact fibers over a given smooth (possibly compact) manifold $M$. This is an easy analogue of Proposition \ref{non-standard-contact-structure} (in case of one line bundle and $B=0$). To describe the construction we need to fix a line bundle $L\to M$ (possibly trivial). It is well known that the vector bundle  $J^1(L)\to M$ of first jets of sections of $L$  carries a natural contact structure (see \cite[Example 3.44]{McDuff-Salamon} for the case of the trivial line bundle $L\to M$). Possibly the best way to describe this contact structure is to note that the canonical symplectic structure on the cotangent bundle of the complement $(L^*)^{\times}$ of the zero section in $L^*$ descends to a contact  structure on the quotient
$T^*((L^*)^{\times})/\RR^{\times}$  (this point of view on contact structures was used already in \cite{LeBrun}; cf. the proof of Proposition \ref{general-LeBrun-lemma}). Now $J^1(L)$ can be canonically identified with this contact manifold $T^*((L^*)^{\times})/\RR^{\times}$.
Our aim is to (relatively) compactify this contact manifold to a b-contact manifold. The b-compactification is given by the projectivized bundle ${\mathbf P}(J^1(L)\times \RR)\to M$ with the critical set $Z$ given by the added hypersurface at infinity.
We describe this $b$-contact structure in local coordinates.

Let $(q^1,..., q^n)$ be local coordinates for $M$ over some open subset over which we fix a local trivialization of $L$. Then 
$J^1(L)$ has local coordinates $(q^1,..., q^n, p_1,...,p_n, t)$ with standard notation, and on $J^1(L)\times \RR$ we have an additional coordinate $z$ corresponding to the second factor of the product. Then on $J^1(L)\times \RR$ with a b-manifold structure given by the critical set $z=0$,  we can define a b-$1$-form by 
$$t\frac{dz}{z}-dt+\sum p_idq^i.$$
This form is homogeneous with respect to the natural $\RR^{\times}$-action on fiber coordinates over $M$ (see the proof of Proposition \ref{non-standard-contact-structure}) and it is non-vanishing outside of the zero section. So it gives rise to a contact b-$1$-form on ${\mathbf P}(J^1(L)\times \RR)$. The proof that the form does not depend on the choice of local coordinates is left to the reader.

\medskip

\section*{Acknowledgements}

The author would like to thank the referee for the exceptionally thorough report 
and for the considerable effort devoted to reviewing the manuscript. Examples \ref{AI1} and \ref{AI2}, together with the example in Lemma \ref{example-semi-contact}, were suggested by ChatGPT (GPT-5.6 Sol) and subsequently verified and rewritten by the author. The author was partially supported by Polish National Centre (NCN) contract numbers 2021/41/B/ST1/03741 and 2025/59/B/ST1/02168.

\bibliographystyle{amsalpha}
\bibliography{References}

\newcommand{\etalchar}[1]{$^{#1}$}
\providecommand{\MR}[1]{}
\providecommand{\bysame}{\leavevmode\hbox to3em{\hrulefill}\thinspace}
\providecommand{\MR}{\relax\ifhmode\unskip\space\fi MR }
\providecommand{\MRhref}[2]{%
  \href{http://www.ams.org/mathscinet-getitem?mr=#1}{#2}
}
\providecommand{\href}[2]{#2}
\begin{thebibliography}{ORCW21}

\bibitem[{Arn}89]{Arnold-Mathematical-methods}
{Arnol{\cprime d}, V. I.}, \emph{Mathematical methods of classical mechanics},
  second ed., Graduate Texts in Mathematics, vol.~60, Springer-Verlag, New
  York, 1989, Translated from the Russian by K. Vogtmann and A. Weinstein.
  \MR{997295}

\bibitem[Ati56]{Atiyah-Krull-Schmidt}
M.~Atiyah, \emph{On the {K}rull-{S}chmidt theorem with application to sheaves},
  Bull. Soc. Math. France \textbf{84} (1956), 307--317. \MR{86358}

\bibitem[AW97]{Andreatta-Wisniewski1997}
Marco Andreatta and Jaros\l aw~A. Wi\'{s}niewski, \emph{A view on contractions
  of higher-dimensional varieties}, Algebraic geometry---{S}anta {C}ruz 1995,
  Proc. Sympos. Pure Math., vol.~62, Amer. Math. Soc., Providence, RI, 1997,
  pp.~153--183. \MR{1492522}

\bibitem[AW01]{Andreatta-Wisniewski2001}
\bysame, \emph{On manifolds whose tangent bundle contains an ample subbundle},
  Invent. Math. \textbf{146} (2001), no.~1, 209--217. \MR{1859022}

\bibitem[BB93]{Beilinson-Bernstein-Jantzen}
A.~Be\u{\i}linson and J.~Bernstein, \emph{A proof of {J}antzen conjectures}, I.
  {M}. {G}el\cprime fand {S}eminar, Adv. Soviet Math., vol.~16, Amer. Math.
  Soc., Providence, RI, 1993, pp.~1--50. \MR{1237825}

\bibitem[BDM{\etalchar{+}}19]{Bradell-et-al}
Roisin Braddell, Amadeu Delshams, Eva Miranda, C\'{e}dric Oms, and Arnau
  Planas, \emph{An invitation to singular symplectic geometry}, Int. J. Geom.
  Methods Mod. Phys. \textbf{16} (2019), no.~suppl. 1, 1940008, 16.
  \MR{3904653}

\bibitem[BDPP13]{BDPP2013}
S\'{e}bastien Boucksom, Jean-Pierre Demailly, Mihai P\u{a}un, and Thomas
  Peternell, \emph{The pseudo-effective cone of a compact {K}\"{a}hler manifold
  and varieties of negative {K}odaira dimension}, J. Algebraic Geom.
  \textbf{22} (2013), no.~2, 201--248. \MR{3019449}

\bibitem[Bea98]{Beauville-Fano-contact}
Arnaud Beauville, \emph{Fano contact manifolds and nilpotent orbits}, Comment.
  Math. Helv. \textbf{73} (1998), no.~4, 566--583. \MR{1639888}

\bibitem[BKK22]{Buczynski-Kapustka-Kapustka}
Jaros{\l}aw Buczy\'{n}ski, Grzegorz Kapustka, and Micha{\l} Kapustka,
  \emph{Special lines on contact manifolds}, Ann. Inst. Fourier (Grenoble)
  \textbf{72} (2022), no.~5, 1859--1909. \MR{4485842}

\bibitem[BW22]{Buczynski-Wisniewski}
Jaros\l~aw Buczy\'{n}ski and Jaros\l aw~A. Wi\'{s}niewski, \emph{Algebraic
  torus actions on contact manifolds}, J. Differential Geom. \textbf{121}
  (2022), no.~2, 227--289, With an appendix by Andrzej Weber. \MR{4466670}

\bibitem[CG18]{Cavalcanti-Gualtieri2018}
Gil~R. Cavalcanti and Marco Gualtieri, \emph{Stable generalized complex
  structures}, Proc. Lond. Math. Soc. (3) \textbf{116} (2018), no.~5,
  1075--1111. \MR{3805052}

\bibitem[CMSB02]{Cho-Miyaoka-Shepherd-Barron}
Koji Cho, Yoichi Miyaoka, and N.~I. Shepherd-Barron, \emph{Characterizations of
  projective space and applications to complex symplectic manifolds}, Higher
  dimensional birational geometry ({K}yoto, 1997), Adv. Stud. Pure Math.,
  vol.~35, Math. Soc. Japan, Tokyo, 2002, pp.~1--88. \MR{1929792}

\bibitem[CVdB10]{Calaque-Van-den-Bergh}
Damien Calaque and Michel Van~den Bergh, \emph{Hochschild cohomology and
  {A}tiyah classes}, Adv. Math. \textbf{224} (2010), no.~5, 1839--1889.
  \MR{2646112}

\bibitem[CZ19]{Chen-Zhu}
Qile Chen and Yi~Zhu, \emph{{$\Bbb A^1$}-curves on log smooth varieties}, J.
  Reine Angew. Math. \textbf{756} (2019), 1--35. \MR{4026447}

\bibitem[Dem02]{Demailly-Contact}
Jean-Pierre Demailly, \emph{On the {F}robenius integrability of certain
  holomorphic {$p$}-forms}, Complex geometry ({G}\"{o}ttingen, 2000), Springer,
  Berlin, 2002, pp.~93--98. \MR{1922099}

\bibitem[DK93]{Dolgachev-Kapranov}
I.~Dolgachev and M.~Kapranov, \emph{Arrangements of hyperplanes and vector
  bundles on {$\bold P^n$}}, Duke Math. J. \textbf{71} (1993), no.~3, 633--664.
  \MR{1240599}

\bibitem[Dru98]{Druel-Contact}
St\'{e}phane Druel, \emph{Structures de contact sur les vari\'{e}t\'{e}s
  alg\'{e}briques de dimension 5}, C. R. Acad. Sci. Paris S\'{e}r. I Math.
  \textbf{327} (1998), no.~4, 365--368. \MR{1649963}

\bibitem[Fuj87]{Fujita1985}
Takao Fujita, \emph{On polarized manifolds whose adjoint bundles are not
  semipositive}, Algebraic geometry, {S}endai, 1985, Adv. Stud. Pure Math.,
  vol.~10, North-Holland, Amsterdam, 1987, pp.~167--178. \MR{946238}

\bibitem[Fuj14]{Fujita-Log-Fano}
Kento Fujita, \emph{Simple normal crossing {F}ano varieties and log {F}ano
  manifolds}, Nagoya Math. J. \textbf{214} (2014), 95--123. \MR{3211820}

\bibitem[Gei08]{Geiges-Contact-topology}
Hansj\"{o}rg Geiges, \emph{An introduction to contact topology}, Cambridge
  Studies in Advanced Mathematics, vol. 109, Cambridge University Press,
  Cambridge, 2008. \MR{2397738}

\bibitem[GKKP11]{GKKP-Extension-of-diff-forms}
Daniel Greb, Stefan Kebekus, S\'{a}ndor~J. Kov\'{a}cs, and Thomas Peternell,
  \emph{Differential forms on log canonical spaces}, Publ. Math. Inst. Hautes
  \'{E}tudes Sci. (2011), no.~114, 87--169. \MR{2854859}

\bibitem[GMP14]{Guillemin-Miranda-Pires}
Victor Guillemin, Eva Miranda, and Ana~Rita Pires, \emph{Symplectic and
  {P}oisson geometry on {$b$}-manifolds}, Adv. Math. \textbf{264} (2014),
  864--896. \MR{3250302}

\bibitem[Got02]{Goto2002}
Ryushi Goto, \emph{Rozansky-{W}itten invariants of log symplectic manifolds},
  Integrable systems, topology, and physics ({T}okyo, 2000), Contemp. Math.,
  vol. 309, Amer. Math. Soc., Providence, RI, 2002, pp.~69--84. \MR{1953353}

\bibitem[Got16]{Goto2016}
\bysame, \emph{Unobstructed deformations of generalized complex structures
  induced by {$C^\infty$} logarithmic symplectic structures and logarithmic
  {P}oisson structures}, Geometry and topology of manifolds, Springer Proc.
  Math. Stat., vol. 154, Springer, [Tokyo], 2016, pp.~159--183. \MR{3555982}

\bibitem[Gra62]{Grauert1962}
Hans Grauert, \emph{\"{U}ber {M}odifikationen und exzeptionelle analytische
  {M}engen}, Math. Ann. \textbf{146} (1962), 331--368. \MR{137127}

\bibitem[Har77]{Har77}
Robin Hartshorne, \emph{Algebraic geometry}, Grad. Texts Math., no.~52,
  Springer-Verlag, New York-Heidelberg, 1977. \MR{0463157}

\bibitem[Har80]{Hartshorne-Stable-reflexive-sheaves}
\bysame, \emph{Stable reflexive sheaves}, Math. Ann. \textbf{254} (1980),
  no.~2, 121--176. \MR{597077}

\bibitem[Kat89]{Kato-Logarithmic-structures}
Kazuya Kato, \emph{Logarithmic structures of {F}ontaine-{I}llusie}, Algebraic
  analysis, geometry, and number theory ({B}altimore, {MD}, 1988), Johns
  Hopkins Univ. Press, Baltimore, MD, 1989, pp.~191--224. \MR{1463703}

\bibitem[Kaw97]{Kawamata1997}
Yujiro Kawamata, \emph{On {F}ujita's freeness conjecture for {$3$}-folds and
  {$4$}-folds}, Math. Ann. \textbf{308} (1997), no.~3, 491--505. \MR{1457742}

\bibitem[Keb02]{Kebekus-CMSB}
Stefan Kebekus, \emph{Characterizing the projective space after {C}ho,
  {M}iyaoka and {S}hepherd-{B}arron}, Complex geometry ({G}\"{o}ttingen, 2000),
  Springer, Berlin, 2002, pp.~147--155. \MR{1922103}

\bibitem[KM98]{Kollar-Mori-book}
J\'{a}nos Koll\'{a}r and Shigefumi Mori, \emph{Birational geometry of algebraic
  varieties}, Cambridge Tracts in Mathematics, vol. 134, Cambridge University
  Press, Cambridge, 1998, With the collaboration of C. H. Clemens and A. Corti,
  Translated from the 1998 Japanese original. \MR{1658959}

\bibitem[KM99]{Keel-McKernan}
Se\'{a}n Keel and James McKernan, \emph{Rational curves on quasi-projective
  surfaces}, Mem. Amer. Math. Soc. \textbf{140} (1999), no.~669, viii+153.
  \MR{1610249}

\bibitem[KMM87]{Kawamata-Matsuda-Matsuki}
Yujiro Kawamata, Katsumi Matsuda, and Kenji Matsuki, \emph{Introduction to the
  minimal model problem}, Algebraic geometry, {S}endai, 1985, Adv. Stud. Pure
  Math., vol.~10, North-Holland, Amsterdam, 1987, pp.~283--360. \MR{946243}

\bibitem[Kol96]{Kollar-Rational}
J\'anos Koll\'ar, \emph{Rational curves on algebraic varieties}, Ergeb. Math.
  Grenzgeb. (3), vol.~32, Springer-Verlag, Berlin, 1996. \MR{1440180}

\bibitem[KPSW00]{KPSW}
Stefan Kebekus, Thomas Peternell, Andrew~J. Sommese, and Jaros\l aw~A.
  Wi\'{s}niewski, \emph{Projective contact manifolds}, Invent. Math.
  \textbf{142} (2000), no.~1, 1--15. \MR{1784795}

\bibitem[KS90]{Kashiwara-Schapira}
Masaki Kashiwara and Pierre Schapira, \emph{Sheaves on manifolds}, Grundlehren
  der mathematischen Wissenschaften [Fundamental Principles of Mathematical
  Sciences], vol. 292, Springer-Verlag, Berlin, 1990, With a chapter in French
  by Christian Houzel. \MR{1074006}

\bibitem[Lan14]{Langer-Lie-algebroids}
Adrian Langer, \emph{Semistable modules over {L}ie algebroids in positive
  characteristic}, Doc. Math. \textbf{19} (2014), 509--540. \MR{3218782}

\bibitem[Laz04]{Laz04}
Robert Lazarsfeld, \emph{Positivity in algebraic geometry. {I}}, Ergeb. Math.
  Grenzgeb. (3), vol.~48, Springer-Verlag, Berlin, 2004, Classical setting:
  line bundles and linear series. \MR{2095472}

\bibitem[LeB95]{LeBrun}
Claude LeBrun, \emph{Fano manifolds, contact structures, and quaternionic
  geometry}, Internat. J. Math. \textbf{6} (1995), no.~3, 419--437.
  \MR{1327157}

\bibitem[LS94]{LeBrun-Salamon}
Claude LeBrun and Simon Salamon, \emph{Strong rigidity of positive
  quaternion-{K}\"{a}hler manifolds}, Invent. Math. \textbf{118} (1994), no.~1,
  109--132. \MR{1288469}

\bibitem[Mel93]{Melrose-book}
Richard~B. Melrose, \emph{The {A}tiyah-{P}atodi-{S}inger index theorem},
  Research Notes in Mathematics, vol.~4, A K Peters, Ltd., Wellesley, MA, 1993.
  \MR{1348401}

\bibitem[MO23]{Miranda-Oms}
Eva Miranda and C\'{e}dric Oms, \emph{Contact structures with singularities:
  from local to global}, J. Geom. Phys. \textbf{192} (2023), Paper No. 104957,
  22. \MR{4629748}

\bibitem[Mor82]{Mori1982}
Shigefumi Mori, \emph{Threefolds whose canonical bundles are not numerically
  effective}, Ann. of Math. (2) \textbf{116} (1982), no.~1, 133--176.
  \MR{662120}

\bibitem[MS98]{McDuff-Salamon}
Dusa McDuff and Dietmar Salamon, \emph{Introduction to symplectic topology},
  second ed., Oxford Mathematical Monographs, The Clarendon Press, Oxford
  University Press, New York, 1998. \MR{1698616}

\bibitem[Ogu18]{Ogus-log}
Arthur Ogus, \emph{Lectures on logarithmic algebraic geometry}, Cambridge
  Studies in Advanced Mathematics, vol. 178, Cambridge University Press,
  Cambridge, 2018. \MR{3838359}

\bibitem[ORCW21]{Occhetta-Romano-SolaConde-Wisniewski}
Gianluca Occhetta, Eleonora~A. Romano, Luis E.~Sol\'{a} Conde, and Jaros\l
  aw~A. Wi\'{s}niewski, \emph{High rank torus actions on contact manifolds},
  Selecta Math. (N.S.) \textbf{27} (2021), no.~1, Paper No. 10, 33.
  \MR{4215378}

\bibitem[Pym17]{Pym2017}
Brent Pym, \emph{Elliptic singularities on log symplectic manifolds and
  {F}eigin-{O}desskii {P}oisson brackets}, Compos. Math. \textbf{153} (2017),
  no.~4, 717--744. \MR{3705241}

\bibitem[Sau89]{Saunders-Jet-book}
D.~J. Saunders, \emph{The geometry of jet bundles}, London Mathematical Society
  Lecture Note Series, vol. 142, Cambridge University Press, Cambridge, 1989.
  \MR{989588}

\bibitem[Tou16]{Touzet2016}
Fr\'{e}d\'{e}ric Touzet, \emph{On the structure of codimension 1 foliations
  with pseudoeffective conormal bundle}, Foliation theory in algebraic
  geometry, Simons Symp., Springer, Cham, 2016, pp.~157--216. \MR{3644247}

\bibitem[Wah83]{Wahl-Ample-in-tangent}
J.~M. Wahl, \emph{A cohomological characterization of {${\bf P}^{n}$}}, Invent.
  Math. \textbf{72} (1983), no.~2, 315--322. \MR{700774}

\bibitem[Win04]{Winkelmann2004}
J\"{o}rg Winkelmann, \emph{On manifolds with trivial logarithmic tangent
  bundle}, Osaka J. Math. \textbf{41} (2004), no.~2, 473--484. \MR{2069097}

\bibitem[Wi{\'{s}}90]{Wisniewski1990}
Jaros{\l}aw~A. Wi{\'{s}}niewski, \emph{On a conjecture of {Mukai}}, Manuscripta
  Math. \textbf{68} (1990), no.~2, 135--141. \MR{1063222}

\end{thebibliography}

\end{document}